\documentclass[preprint,12pt]{elsarticle}
\usepackage{epsfig}
\usepackage{amsmath, amssymb, amscd, amsthm, amsfonts}
\usepackage{lineno,hyperref,color,bm}
\usepackage{mathtools}
\usepackage{epstopdf}

\newtheorem{theorem}{Theorem}

\newtheorem{example}[theorem]{Example}
\newtheorem{remark}[theorem]{Remark}

\begin{document}

\begin{frontmatter}

\title{Fractional Buffer Layers: \\ Absorbing Boundary Conditions for Wave Propagation}

\author[SHUaddress]{Min Cai} 

\author[Brownaddress]{Ehsan Kharazmi\corref{mycorrespondingauthor}}
\cortext[mycorrespondingauthor]{Corresponding author}
\ead{ehsan\_kharazmi@brown.edu}

\author[SHUaddress]{Changpin Li} 

\author[Brownaddress,PNNLaddress]{George Em Karniadakis}

\address[SHUaddress]{Department of Mathematics, Shanghai University, 99 Shangda Road, Shanghai 200444, People's Republic of China.}

\address[Brownaddress]{Division of Applied Mathematics, Brown University, 170 Hope St, Providence, RI 02906, USA.}

\address[PNNLaddress]{Pacific Northwest National Laboratory, Richland, WA 99354, USA.}

\begin{abstract}
We develop fractional buffer layers (FBLs) to absorb propagating waves without reflection in bounded domains. Our formulation is based on variable-order spatial fractional derivatives. We select a proper variable-order function so that dissipation is induced to absorb the coming waves in the buffer layers attached to the domain. In particular, we first design proper FBsL for the one-dimensional one-way and two-way wave propagation. Then, we extend our formulation to two-dimensional problems, where we introduce a consistent variable-order fractional wave equation. In each case, we obtain the fully discretized equations by employing a spectral collocation method in space and Crank-Nicolson or Adams-Bashforth method in time. We compare our results with the perfectly matched layer (PML) method and show the effectiveness of FBL in accurately suppressing any erroneously reflected waves, including corner reflections in two-dimensional rectangular domains. FBLs can be used in conjunction with any discretization method appropriate for fractional operators describing wave propagation in bounded or truncated domains. 
\end{abstract}


%\item Removing the corner reflections in two-dimensional problems.
%\item Less number of equations in the system of equations and therefore more computationally efficient in higher dimensions.
%\item Different buffer layers with various characteristics in each propagation direction. 
%\item Flexibility to absorb anisotropic wave propagation.


\end{frontmatter}

%\linenumbers

%%%%%%%%%%%%%%%
\section{Introduction}
%%%%%%%%%%%%%%%
Waves are omnipresent in nature in diverse physical and biological phenomena. They are governed by first-order, second-order or higher order  
partial differential equations (PDES) that lead to one-way or two-way waves, giving rise to oscillating solutions that propagate through spatio-temporal domains while conserving energy in lossless media. Examples include the scalar wave equations for pressure waves in gases, Maxwell’s equations in electromagnetism, Schr\"{o}dinger’s equation in quantum mechanics, and elastic vibrations. In solving numerically the governing equations of these examples on a finite domain, it is necessary to truncate the computational domain in a certain way that does not introduce significant artifacts into the computation. Therefore, an efficient model for propagation of waves from the interior of a finite domain and absorption through its boundaries is desirable in numerous physical problems \cite{BaylissTurkel@1980}. Several approaches have been studied in the literature. Some works have focused on finding analytical boundary conditions for the differential equations and then discretizing the analytical conditions \cite[and references therein]{Xavier@2008}. Some other works are dedicated to constructing absorbing boundary conditions via directly working with approximations to the wave equation \cite[and references therein]{Higdon@1986,Higdon@1987}. A well-known approach in the literature is the perfectly matched layer (PML) method, which in theory absorbs strongly the outgoing waves from the interior of a computational domain without reflecting them back into the interior. This method was first introduced by J.P. B\'{e}renger \cite{Berenger@1994}, and  thereafter a number of works based on PML and corresponding numerical schemes were developed; see for example \cite[and references therein]{Fang@JSC2019,Huang@JCAM2020}. A comprehensive introduction to the PML for wave equations can be found in \cite{Johnson2007} and other modifications to the PML method are also proposed in \cite{Grote@2010, Kim@2019}. It has been observed, however, that most of these PML based methods may still suffer from some reflections in the discretized wave equation and specifically exhibit the issue of corner reflections in high-dimensional problems with singular domains.

We develop a new absorbing layer, namely a \textit{fractional buffer layer} (FBL), by exploiting
the flexibility and expressivity of variable-order fractional operators. It is widely recognized that fractional calculus has important applications in various scientific fields. In particular, the fractional derivatives that extend the notion of their integer-order counterparts have been shown to provide a powerful mathematical tool that can be used to describe many physical problems. especially anomalous transport \cite{atanackovic2009generalized,LiCai@2019,MetzierKlafter@PR2000,Uchaikin@2013,Zaslavsky@PR2002}. In a more general setting, variable-order fractional derivatives further extend the fixed fractional order to a spatio-temporal variable function. They can accurately describe the multi-scale behavior of systems with temporally/spatially varying properties \cite{Coimbra@AP2003, Chechkin@JPAMG2005,PatnaikHollkampHollkamp@PMPES2020,Sun@PhysA2009}.

\begin{figure}[h!]
\centering
\includegraphics[width=1\linewidth]{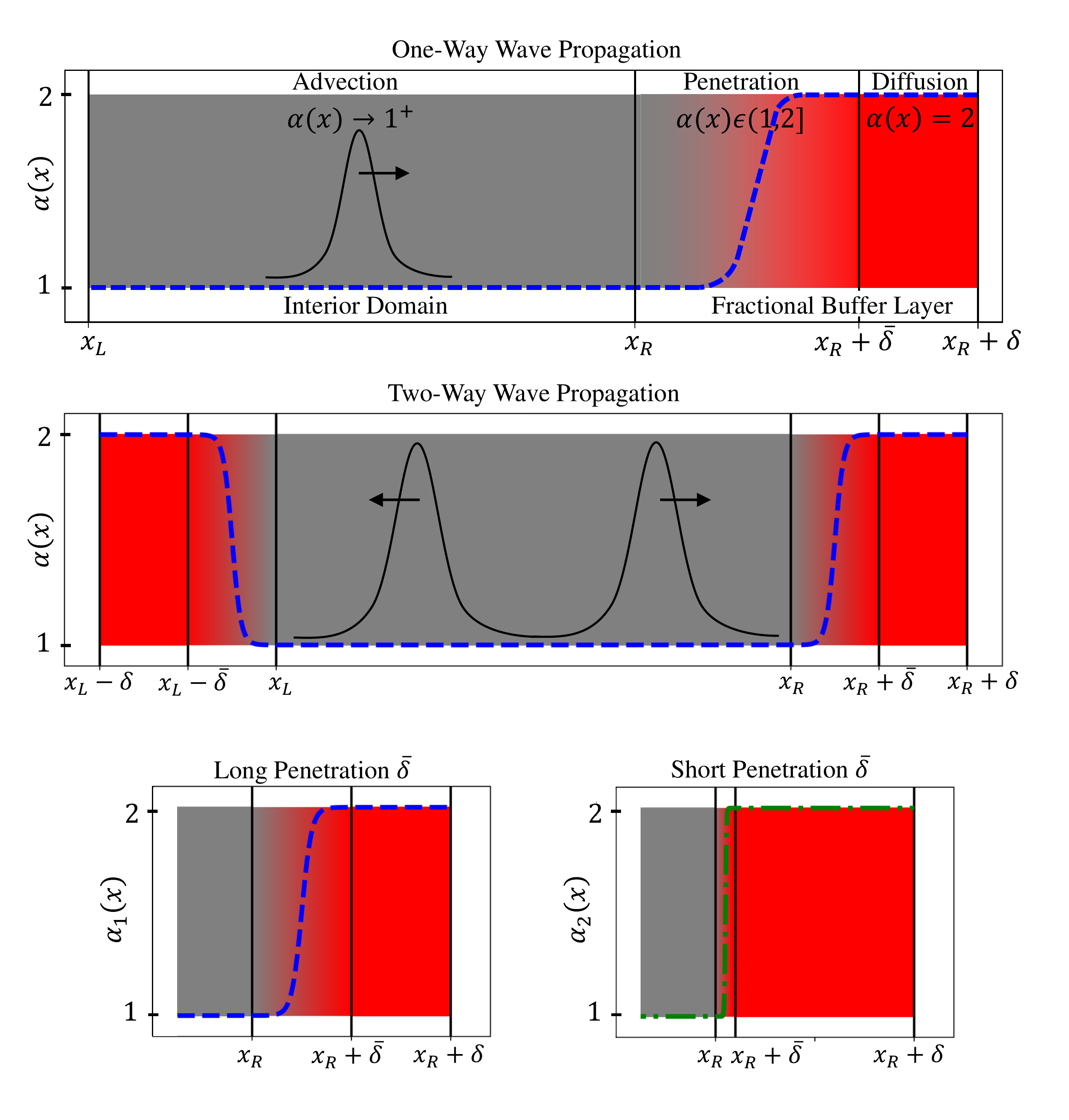}
\caption{Fractional buffer layers (FBLs) for absorbing one-dimensional waves. Top: FBL for one-way wave propagation. Middle: FBL for two-way wave propagation. Bottom: Different characterization of FBLs with long and short penetration regions. The gray shaded area shows the interior domain $[x_L, x_R]$ with advection dominating. The red shaded area shows the full absorption area with diffusion dominating.}
\label{Fig: FBL_Master}
\end{figure}

One of the interesting applications of the variable-order fractional differential operators is in modeling wave propagation in finite domains. In \cite{Zhao@JCP2015}, a variable-order time fractional differential operator is employed to control wave reflections in truncated computational domains 
% of the integer-order one-dimensional wave equation 
by switching from a wave- to a diffusion-dominated equation at the boundaries. In that approach, a priori knowledge of the time that waves reach the boundary is needed to effectively design the variable-order function. In the current paper, we consider the application of variable-order space fractional differential operators, where the variable-order function is solely a function of space variables. We formulate our method by extending the space domain and attaching a \textit{fractional buffer layer} (FBL) to its boundaries. Then, we define a proper variable-order function such that we recover the original equation inside the interior domain of interest and introduce dissipation in the buffer layer to absorb the wave. A schematic plot of one-dimensional FBLs is shown in Fig. \ref{Fig: FBL_Master}. The top panel shows the one-dimensional one-way wave, where the wave propagates to the right inside the interior domain and then penetrates the buffer but it is fully absorbed in the FBL appended to the right boundary. The middle panel shows a similar approach but for one-dimensional two-way waves that propagate in both directions, and thus we append FBLs to the both boundaries. 

The left- and right-sided variable-order fractional derivatives that are used in developing FBLs are given in the Riemann-Liouville sense as 
\begin{flalign}
\label{eq:LeftRL}
\qquad
\prescript{RL}{x_L}{\mathcal{D}}_{x}^{\alpha(x)} u(x)
=\frac{1}{\Gamma(n-\alpha(x))}\frac{\partial^n}{\partial x^n}
\int_{x_L}^{x}(x-s)^{n-1-\alpha(x)}u(s){\rm d}s,
&&
\end{flalign} 
and 
\begin{flalign}
\label{eq:RightRL}
\qquad
\prescript{RL}{x}{\mathcal{D}}_{x_R}^{\alpha(x)}u(x)
=\frac{1}{\Gamma(n-\alpha(x))}\left(-\frac{\partial}{\partial x}\right)^{n}
\int_{x}^{x_R}(s-x)^{n-1-\alpha(x)}u(s){\rm d}s,
&&
\end{flalign} 
respectively, where $n-1<\alpha(x)\leq n\in\mathbb{Z}^{+}$ and $x\in[x_L, x_R]$. 
%
% The left- and right-sided variable-order fractional derivatives in Riemann-Liouville sense are given as 
% \begin{flalign}
% \label{eq:LeftRL}
% \qquad
% \prescript{RL}{x_L}{\mathcal{D}}_{x}^{\alpha(x)} u(x)
% =\frac{1}{\Gamma(n-\alpha(x))}\frac{\partial^n}{\partial x^n}
% \int_{x_L}^{x}(x-s)^{n-1-\alpha(x)}u(s){\rm d}s,
% &&
% \end{flalign} 
% and 
% \begin{flalign}
% \label{eq:RightRL}
% \qquad
% \prescript{RL}{x}{\mathcal{D}}_{x_R}^{\alpha(x)}u(x)
% =\frac{1}{\Gamma(n-\alpha(x))}\left(-\frac{\partial}{\partial x}\right)^{n}
% \int_{x}^{x_R}(s-x)^{n-1-\alpha(x)}u(s){\rm d}s,
% &&
% \end{flalign} 
% respectively. 
We can show for any integer $n\in\mathbb{Z}^{+}$ that
\begin{align}
\label{Eq: consistency}
\lim\limits_{\alpha(x)\rightarrow n}
\left(\prescript{RL}{x_L}{\mathcal{D}}_{x}^{\alpha(x)}\right)(\cdot)
=\frac{\partial^{n}}{\partial x^{n}} (\cdot), 
\quad 
\lim\limits_{\alpha(x)\rightarrow n}
\left(\prescript{RL}{x}{\mathcal{D}}_{x_R}^{\alpha(x)}\right)(\cdot)
=\left(-\frac{\partial}{\partial x}\right)^{n} (\cdot).
\end{align} 
This consistency property of variable-order fractional derivatives provides the flexibility in our formulation to continuously switch between a first-order derivative as $\alpha(x) \rightarrow 1$ (an advection operator) and a second-order derivative as $\alpha (x) \rightarrow 2$ (a diffusion operator) over the buffer layer. 
%
% A schematic plot of one-dimensional FBLs is shown in Fig. \ref{Fig: FBL_Master}. The top panel shows the one-dimensional one-way wave, where the wave propagates to the right inside the interior domain and then penetrates and is fully absorbed in the FBL appended to the right boundary. The middle panel shows similar approach but for one-dimensional two-way waves that propagate to both directions and thus we append FBLs to the both boundaries. 
%
We discuss in detail the characterization of FBLs and the choices of variable-order function with long and short penetration regions in Section \ref{Sec:Advc}.

%%%%%%%%%%%%%%%%%%%%%%%%%%%%%%%%%%%%
%%%%%%%%%%%%%%%%%%%%%%%%%%%%%%%%%%%%%%

The advantages of our proposed FBL method compared to other approaches such as PML method are: 
\begin{itemize}
    % \item {\color{blue}Absorbing the coming waves of the advection dominant (fractional) equation without reflection.}
    \item Removing the corner reflection in two-dimensional problems.
    \item Less number of equations in the system of equations and therefore more computationally efficient in higher dimensions.
    \item Different buffer layers with various characteristics in each propagation direction. 
    \item Flexibility to absorb anisotropic wave propagation.
    
\end{itemize}

We obtain the fully discrete numerical schemes for the proposed FBLs by employing a spectral collocation method in space and second-order time integration techniques. The detailed derivation of the numerical methods is given in \ref{Sec:FullyDiscreteScheme}.
Several other numerical methods have been also developed in the literature to solve equations with variable-order fractional differential operators. Finite difference methods were developed for variable-order time and space fractional wave equations in \cite{Liu@FCAA2013, Sweilam@FCAA2012}. High-order spectral collocation and Petrov-Galerkin spectral methods have been also developed in \cite{ SamieeUnified1, SamieeUnified2, ZayernouriKarniadakis@JCP2015, ZengZhangKarniadakis@SIAM2015, ZhaoMaoKarniadakis@CMAME2019} and successfully used and analysed in several fractional differential equations of distributed-order \cite{kharazmiPGDis@SISC2017, SamieefullyDis} and random order in the context of uncertainty quantification \cite{AkhavanUQ@Arxiv2020,kharazmiUQ@JVVUQ2019}.
% , and sensitivity analysis \cite{kharazmiFSEM@JSC2018}. 
A comprehensive review of high-order numerical methods for fractional differential equations is given in \cite{KharazmiRev}. 

The rest of the paper is organized as follows. In section \ref{Sec:Advc}, we develop the FBL for one-way wave equations. We extend the formulation to the two-way wave equations in one- and two-dimensions in sections \ref{Sec:Wave_1D} and  \ref{Sec:Wave_2D}, respectively. We end the paper with some concluding remarks in section \ref{sec:Conclusion}. The Appendix contains details of the numerical discretization method used in the current work.

%%%%%%%%%%%%%%%
\section{FBL: One-dimensional one-way waves}\label{Sec:Advc}
%%%%%%%%%%%%%%%
We first consider the one-dimensional one-way wave equation in the bounded domain $\Omega = (x_L,x_R)$. The governing equation is the advection equation given as
\begin{equation}\label{eq: Advec}
\frac{\partial}{\partial t}u(x,t)
= V \frac{\partial}{\partial x}u(x,t), 
\end{equation}
subject to a given initial condition $u(x,0)=u_{0}(x,t)$ and a suitable boundary condition, where $|V|$ is the magnitude of the propagation velocity. The sign of $V$ identifies the direction of wave propagation in a sense that when $V>0$ and $V<0$ the wave propagates to the left $(-x)$ and to the right $(+x)$, respectively.

%%%%%%%%%%%%%%%
\subsection{Formulation of FBL}
%%%%%%%%%%%%%%%
We develop the FBL for the one-dimensional one-way wave equation \eqref{eq: Advec} in order to absorb the coming waves at the boundary. Our formulation is based on the variable-order space fractional diffusion equation, in which we replace the first-order spatial derivative in \eqref{eq: Advec} with a corresponding left- or right-sided variable-order fractional derivative of Riemann-Liouville sense and let the variable-order be $\alpha(x)\in(1,2]$. We see from the consistency relations \eqref{Eq: consistency} that the variable-order fractional derivative converges to the first-order derivative when $\alpha(x)\rightarrow 1$, and therefore we recover \eqref{eq: Advec} in the limiting case. Moreover, when $\alpha(x)\rightarrow2$ we obtain a diffusion equation. Therefore, the variable-order fractional derivative can be thought of as a transition from an advection to a diffusion process. This provides us the means to adjust the variable-order $\alpha(x)$ in order to design a compliant process that behaves as advection inside the interior domain of interest and then switches to diffusion to damp out the erroneous motion in the buffer layer. \\

\noindent $\bullet$ \textbf{Characterization of the buffer layer}:
To construct the FBL, we first define a buffer layer of size $\delta$ next to the boundaries and extend the computational domain to $(x_L-\delta, x_R)$ and $(x_L, x_R+\delta)$ for $V>0$ and $V<0$, respectively. We select the variable-order derivative such that we have the advection dominant in $(x_L, x_R)$ and then transition to diffusion dominant in the buffer layers. We divide the buffer layer into two regions, namely \textit{penetration (transition)} region, and \textit{fully diffusion} region. In the penetration region, we let the wave go through the buffer layer and start to transition (smoothly or abruptly) into the next fully diffusion region. Thus, we define the variable-order function such that 
\begin{equation}
\label{eq:VO fun Advec1}
\left\{
\begin{array}{lll}
\alpha(x)\rightarrow1^{+},
&x\in[x_L,x_R],
&: \text{\textit{advection}}
\\[3pt]
\alpha(x)\in(1,2],
&x\in[x_L-\bar{\delta},x_L)
\vee
(x_R,x_R+\bar{\delta}],
&: \text{\textit{penetration}}
\\[3pt]
\alpha(x)=2,
&x\in[x_L-\delta,x_L-\bar{\delta})
\vee
(x_R+\bar{\delta},x_R+\delta],
&: \text{\textit{diffusion}}
\end{array}
\right.
\end{equation} 
where $\bar{\delta}\in(0,\delta)$ is the penetration length.

We consider two main classes of the variable-order functions: \textit{i}) the step function that has a very short penetration region, and \textit{ii}) the smooth $\tanh$ function with the tuning parameter $\omega$ to adjust the characteristics of the buffer layer. Their definitions for the right-propagating advection equation are given as follows, and can be mirrored for the left-propagating case. 
\begin{itemize}
\item Step function
\begin{flalign}
\label{eq:RightStepFunction}
   \alpha(x)
   =\left\{
   \begin{array}{ll}
   1+\epsilon, &x\in[x_L,x_R+\bar{\delta}],
   \\[3pt]
   2, &x\in(x_R+\bar{\delta},x_R+\delta].
   \end{array}
   \right.
   &&
\end{flalign}

\item Smooth $\tanh$ function 
\begin{flalign}
\label{eq:RightSmoothFunction}
  \alpha(x)
  =\left\{
  \begin{array}{ll}
   1+\epsilon, &x\in[x_L,x_R],
   \\[3pt] 
    1.5+(0.5-\epsilon) \tanh(\omega(x-x_R-\frac{\bar{\delta}}{2})), &x\in(x_R, x_R+\bar{\delta}],
  \\[3pt]
   2, &x\in(x_R+\bar{\delta},x_R+\delta].
   \end{array}
  \right.
  &&
\end{flalign}
\end{itemize}

\noindent Here $\epsilon$ takes a very small positive value (we choose $\epsilon = 10^{-5}$ in numerical simulations). In the case of smooth $\tanh$ function, we choose the parameter $\omega$ based on the length of the penetration region $\bar{\delta}$ such that we do not violate condition \eqref{eq:VO fun Advec1}.\\

\noindent $\bullet$ \textbf{Left propagation}: In this case $V>0$, and the one-way wave moves to the left in the interior domain $[x_L,x_R]$. Our aim is to absorb it (let it decay) in a buffer layer $[x_L-\delta,x_L)$ next to the left boundary. Therefore, we let $(x,t)\in[x_L-\delta,x_R]\times[0,T]$ and construct the FBL by considering the left-sided variable-order space fractional derivative. Hence, 
\begin{equation}
\label{eq:FracAdvcLeftRL}
 \frac{\partial }{\partial t}u(x,t)
=|V| \prescript{RL}{x_L-\delta}{\mathcal{D}}_{x}^{\alpha(x)}u(x,t),
\end{equation}
subject to the initial condition $u(x,0)=u_{0}(x)$ and the boundary conditions $u(x_L-\delta,t)=0$ and $u(x_R,t)=u_{0}(x_R)$.\\

\noindent $\bullet$ \textbf{Right propagation}: In this case $V<0$, and the one-way wave moves to the right in the interior domain $[x_L,x_R]$. Our aim is to absorb it (let it decay) in a buffer layer $(x_R,x_R+\delta]$ next to the right boundary. Therefore, we let $(x,t)\in[x_L,x_R+\delta]\times[0,T]$ and construct the FBL by considering the right-sided variable-order space fractional derivative. Hence, 
\begin{equation}
\label{eq:FracAdvcRightRL}
 \frac{\partial }{\partial t}u(x,t)
=|V|\prescript{RL}{x}{\mathcal{D}}_{x_R+\delta}^{\alpha(x)}u(x,t),
\end{equation}
subject to the initial condition $u(x,0)=u_{0}(x)$ and the boundary conditions $u(x_L,t)=u_{0}(x_L)$ and $u(x_R+\delta,t)=0$.\\

% \noindent $\bullet$ \textbf{Choice of variable-order function $\alpha(x)$}: We consider two main classes of the variable-order functions. The step function has a very short penetration region. The smooth function has the tuning parameter to adjust the characteristics of the buffer layer. We give their definitions for the right-propagating advection equation and mirror these functions for the left-propagating advection equation. 
% \begin{itemize}
% \item Step function
% \begin{flalign}
% \label{eq:RightStepFunction}
%   \alpha(x)
%   =\left\{
%   \begin{array}{ll}
%   1+\epsilon, &x\in[x_L,x_R+\bar{\delta}],
%   \\[3pt]
%   2, &x\in(x_R+\bar{\delta},x_R+\delta].
%   \end{array}
%   \right.
%   &&
% \end{flalign}

% \item Smooth tanh function 
% \begin{flalign}
% \label{eq:RightSmoothFunction}
%   \alpha(x)
%   =\left\{
%   \begin{array}{ll}
%   1+\epsilon, &x\in[x_L,x_R],
%   \\[3pt] 
%     1.5+(0.5-\epsilon) \tanh(\omega(x-x_R-\frac{\bar{\delta}}{2})), &x\in(x_R, x_R+\bar{\delta}],
%   \\[3pt]
%   2, &x\in(x_R+\bar{\delta},x_R+\delta].
%   \end{array}
%   \right.
%   &&
% \end{flalign}
% \end{itemize}

% \noindent Here $\epsilon$ takes a very small positive value (we choose $\epsilon = 10^{-5}$ in numerical simulations). In the case of smooth tanh function, we choose the parameter $\omega$ based on the length of the penetration region $\bar{\delta}$ such that we do not violate condition \eqref{eq:VO fun Advec1}.

\begin{remark}[Length of the FBL]\label{remark:VO functions}
In numerical simulations we will demonstrate that the length of FBL, $\delta$, does not have significant effect on the absorption of the coming wave provided that the variable-order function satisfies  condition  \eqref{eq:VO fun Advec1}, and thus the system in the buffer layer becomes fully diffusion-dominated. It is, however, the penetration length $\bar\delta$ and the slope of variable order function $\omega$, which become important. For more details, we refer to the discussion in the next subsection.
\end{remark}

%%%%%%%%%%%%%%%%%%%%%%%%
\subsection{Numerical simulations}
%%%%%%%%%%%%%%%%%%%%%%%%
First, we investigate the performance of FBLs in absorbing the coming waves for the one-way wave equation \eqref{eq: Advec}. We assume that the wave is propagating to the right and therefore we consider \eqref{eq:FracAdvcRightRL} to formulate the right FBL. We employ a spectral collocation method with $P$ collocation points in space and use a second-order time integration technique with time steps $\tau$. The details of the fully discrete scheme are provided in  \ref{Sec:FullyDiscreteScheme}. We note that the analytical solution to \eqref{eq: Advec} with $V=-1$ is $u^{exact}(x,t)=u_{0}(x-t)$, which we consider it as the reference solution to compare our results in the interior domain. We define the pointwise error at any fixed time $t$ as $ |u(x,t)-u^{exact}(x,t)|$.

\begin{example}
\label{Ex: FBL one way right}
We consider the FBL for right-propagation advection equation given in  \eqref{eq:FracAdvcRightRL} with the propagation velocity $V=-1$. The interior domain is $(x_L,x_R)=(-5,5)$ with a buffer layer of length $\delta>0$ next to the right boundary; thus, $x\in [x_L,x_R+\delta]$. We assume a smooth initial condition of the form $u_{0}(x)=e^{-x^2}$ and consider different characterizations of the smooth $\tanh$ variable-order function in \eqref{eq:RightSmoothFunction} with fixed parameter $\epsilon=10^{-5}$. In particular, we consider long and short penetration regions with moderate and steep variable-order functions.
\end{example}

\begin{figure}
\centering
$t=0$\\[-13 pt]
\noindent\rule{4cm}{0.7pt}
\\
\includegraphics[clip, trim=0.5cm 0cm 1.7cm 0cm, width=0.47\linewidth]{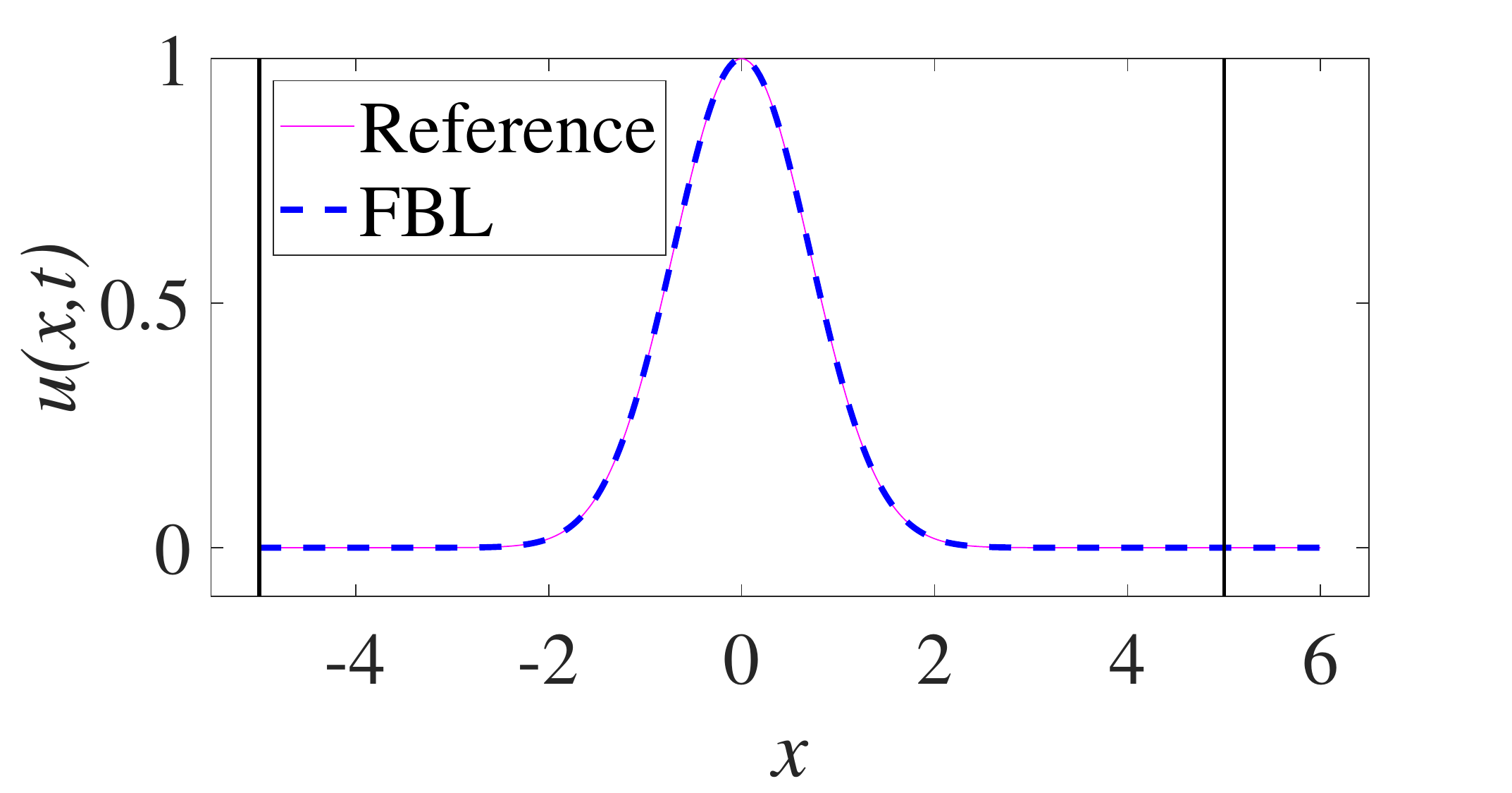}
\\[-6 pt]
$t=3$\\[-13 pt]
\noindent\rule{4cm}{0.7pt}
\\
\includegraphics[clip, trim=0.5cm 0cm 1.7cm 0cm, width=0.47\linewidth]{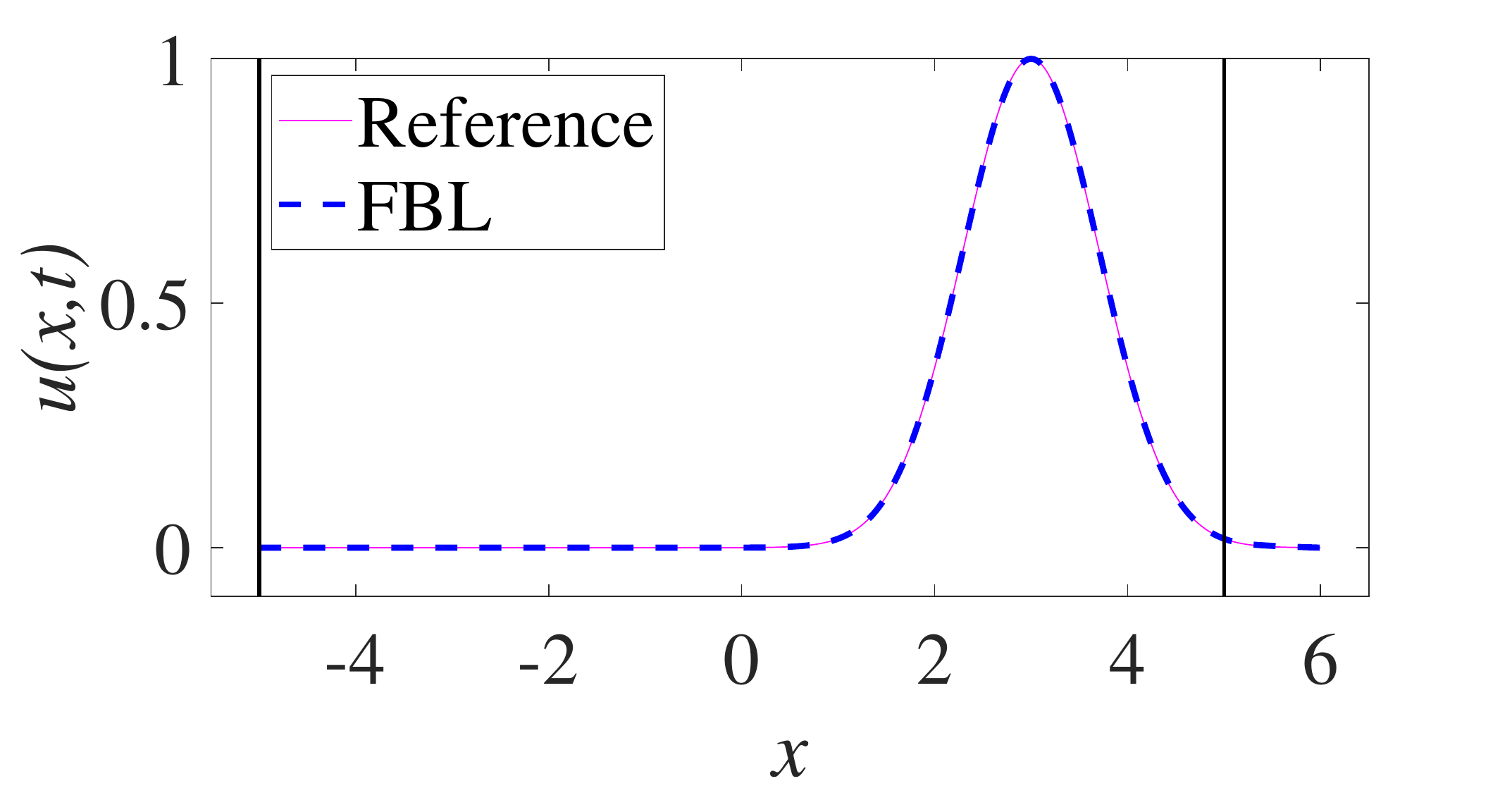}
\includegraphics[clip, trim=0.5cm 0cm 1.7cm 0cm, width=0.47\linewidth]{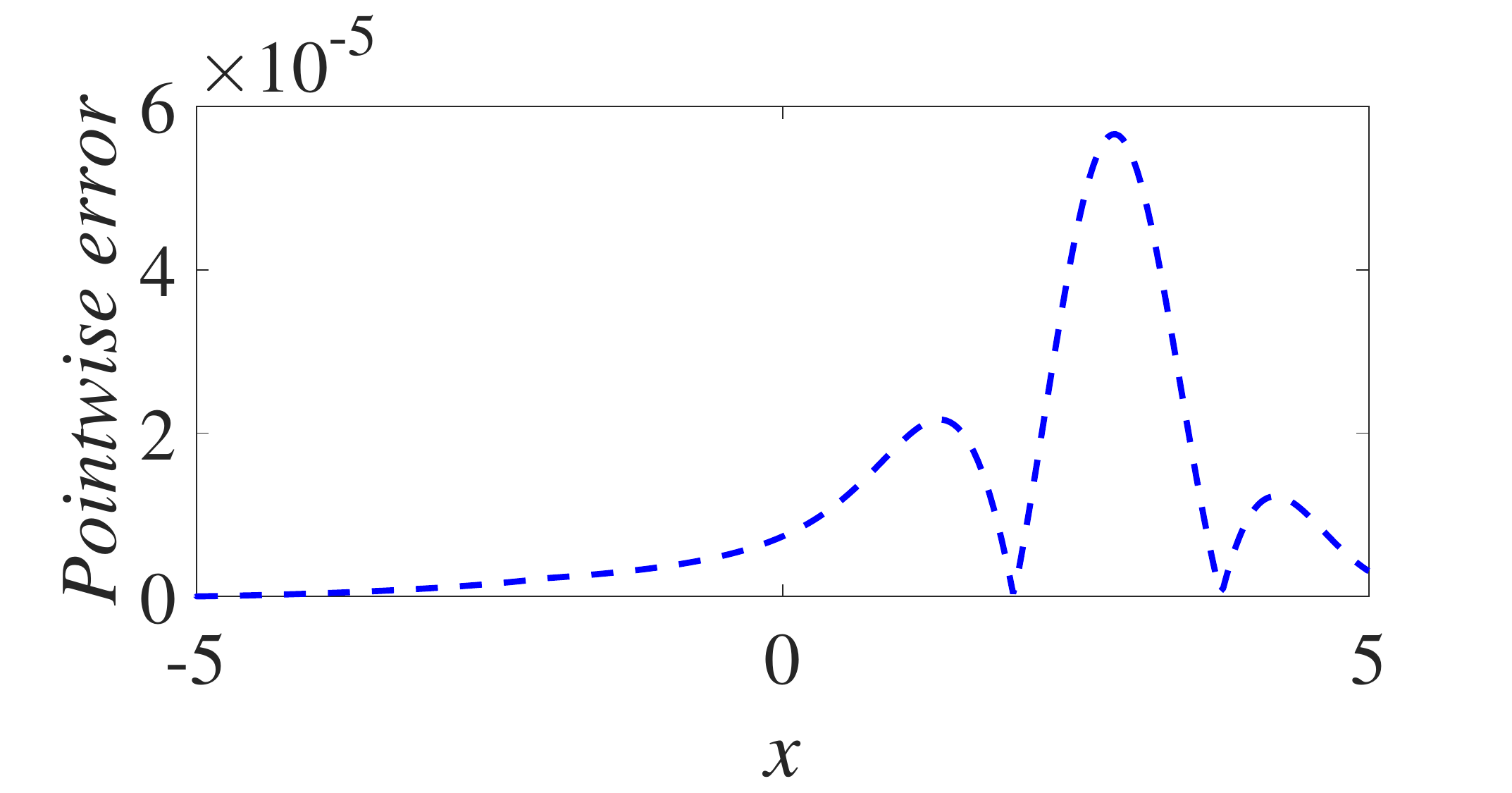}
\\[-9 pt]
$t=6$\\[-13 pt]
\noindent\rule{4cm}{0.7pt}
\\
\includegraphics[clip, trim=0.5cm 0cm 1.7cm 0cm, width=0.47\linewidth]{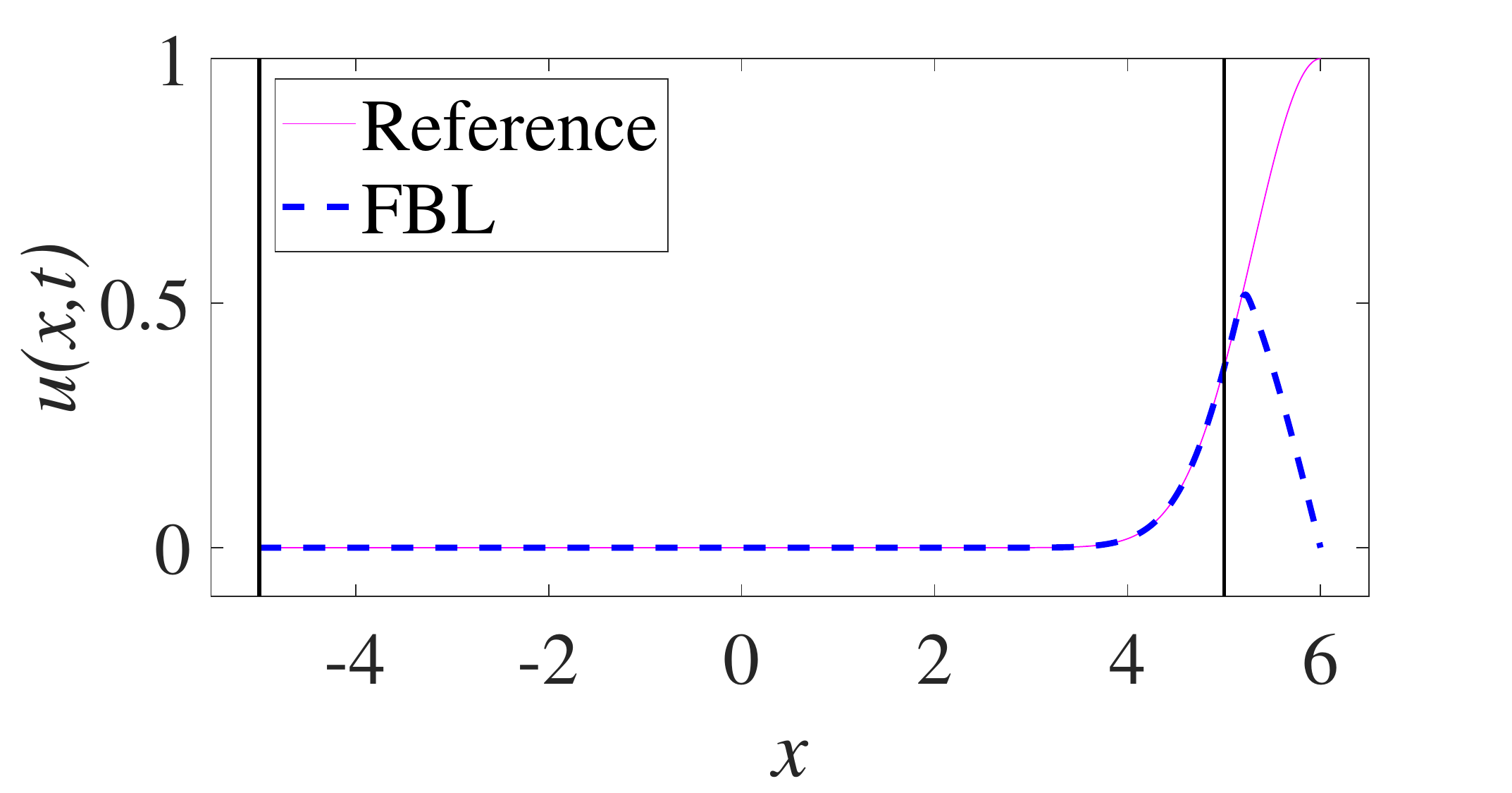}
\includegraphics[clip, trim=0.5cm 0cm 1.7cm 0cm, width=0.47\linewidth]{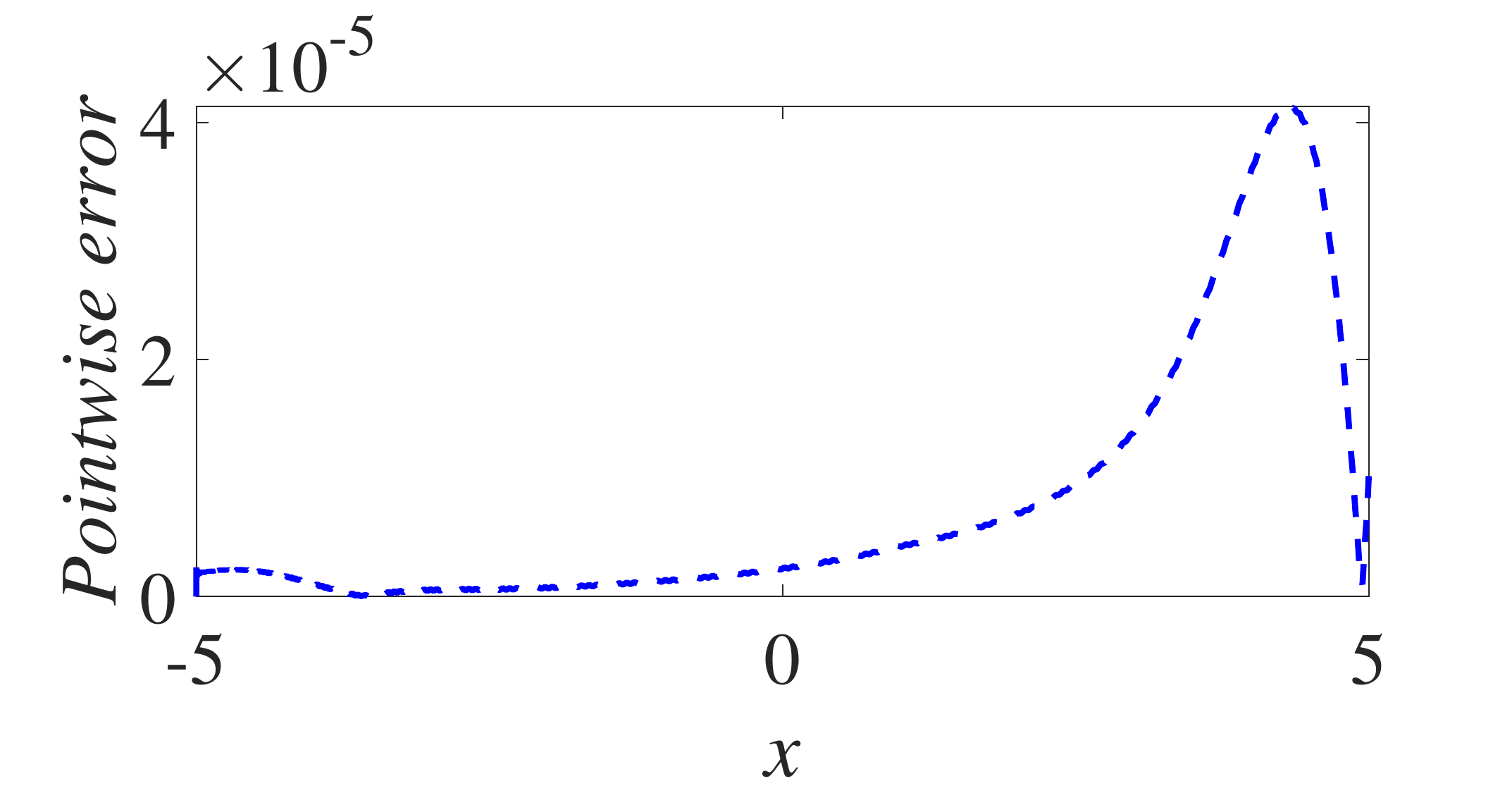}
\\[-9 pt]
$t=9$\\[-13 pt]
\noindent\rule{4cm}{0.7pt}
\\
\includegraphics[clip, trim=0.5cm 0cm 1.7cm 0cm, width=0.47\linewidth]{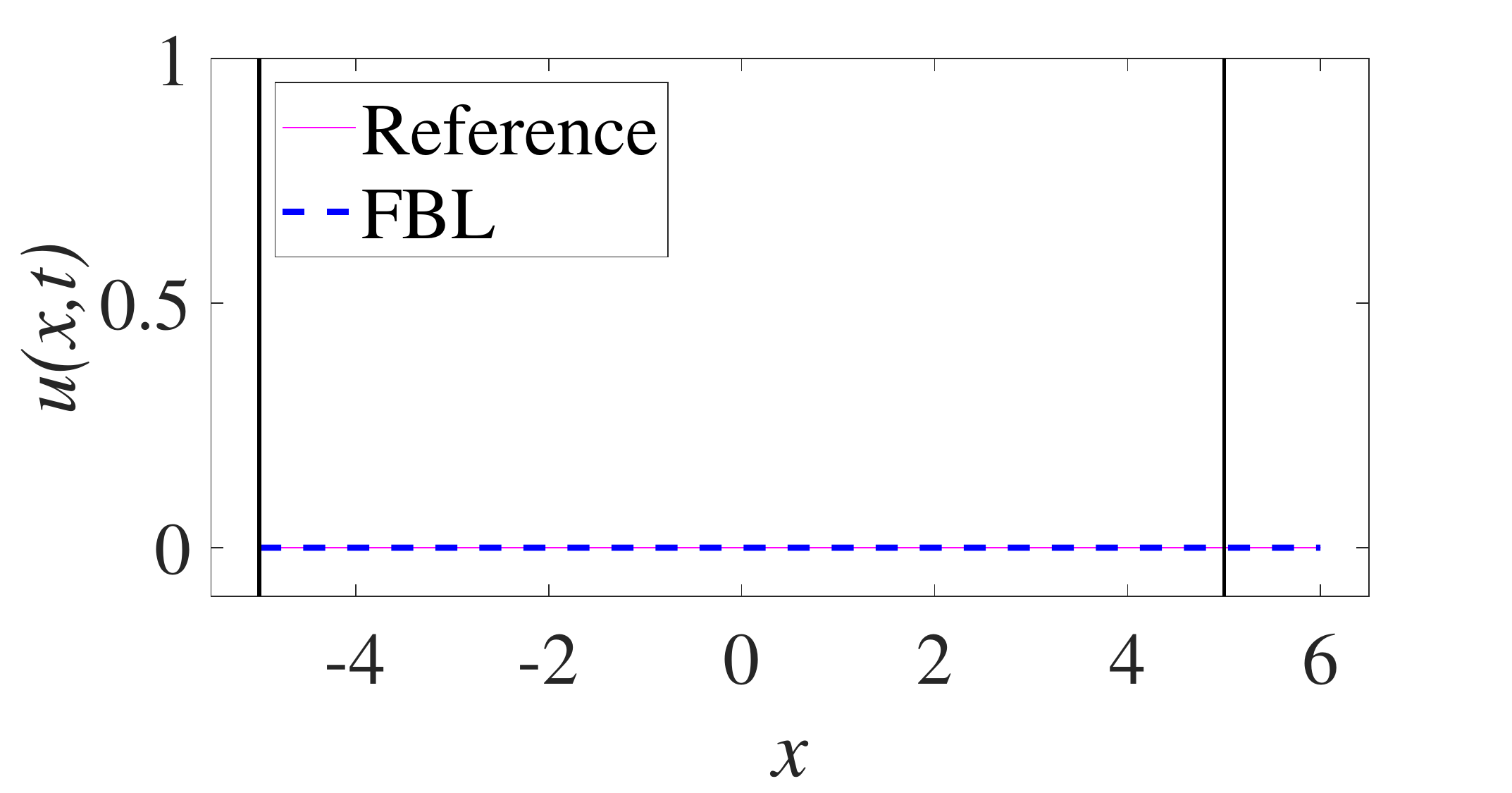}
\includegraphics[clip, trim=0.5cm 0cm 1.7cm 0cm, width=0.47\linewidth]{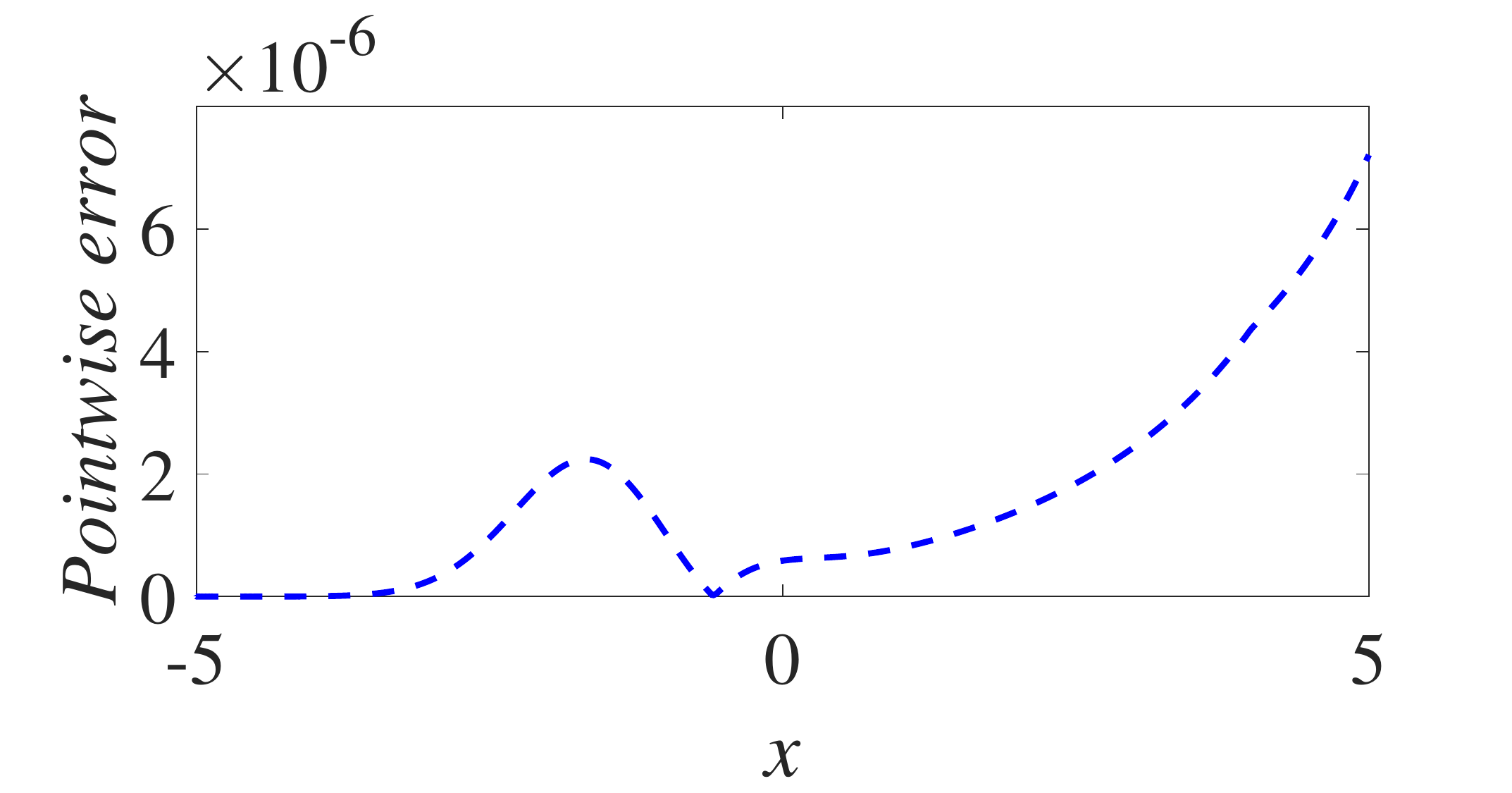}
%%%%%%%%%%
\vspace{-0.2 in}
\caption{One-dimensional one-way wave: The wave propagates to the right and then is absorbed in the FBL. Left column: FBL and reference solution. Right column: pointwise error at $t=0,3,6,9$.
The vertical black solid lines indicate the boundaries of the interior domain. The FBL is characterized as $(\delta,\bar{\delta},\omega)=(1,\frac{\delta}{2},20)$. Other details are given in Example \ref{Ex: FBL one way right}.
}
\label{Fig:NumSol_RefSol_PMLSol_C41Delt1}
\end{figure}

Figure \ref{Fig:NumSol_RefSol_PMLSol_C41Delt1} shows a successful implementation of FBL for Example \ref{Ex: FBL one way right} with parameters $(\delta,\bar{\delta},\omega)=(1,0.5,20)$. We recall that the parameters $\delta$, $\bar{\delta}$, and $\omega$ are the length of buffer layer, length of penetration region, and slope of variable function, respectively. Here the number of collocation points and time step are chosen as $(P,\tau)=(500,10^{-3})$. In Fig. \ref{Fig:NumSol_RefSol_PMLSol_C41Delt1}, the left column shows different snapshots of the wave  reaching the right boundary, penetrating into the buffer layer, and being fully absorbed. The right column shows the point-wise error in the interior domain.

% , which consists of the numerical error and a slight diffusion as mentioned latter in Remark \ref{Rem:Diff}. 

%%%%%%%%%%%%%%%
\subsection{Different sources of error in FBLs}
%%%%%%%%%%%%%%%

% \begin{remark}[Diffusive behaviour in the interior domain]\label{Rem:Diff}

\begin{figure}[h!]
\centering 
\quad 
{\scriptsize $t=1$} \qquad\qquad\qquad\qquad 
{\scriptsize $t=3$} \qquad\qquad\qquad\qquad 
{\scriptsize $t=5$}\\
\includegraphics[clip, trim=0.5cm 0cm 1.7cm 0cm, width=0.3\linewidth]{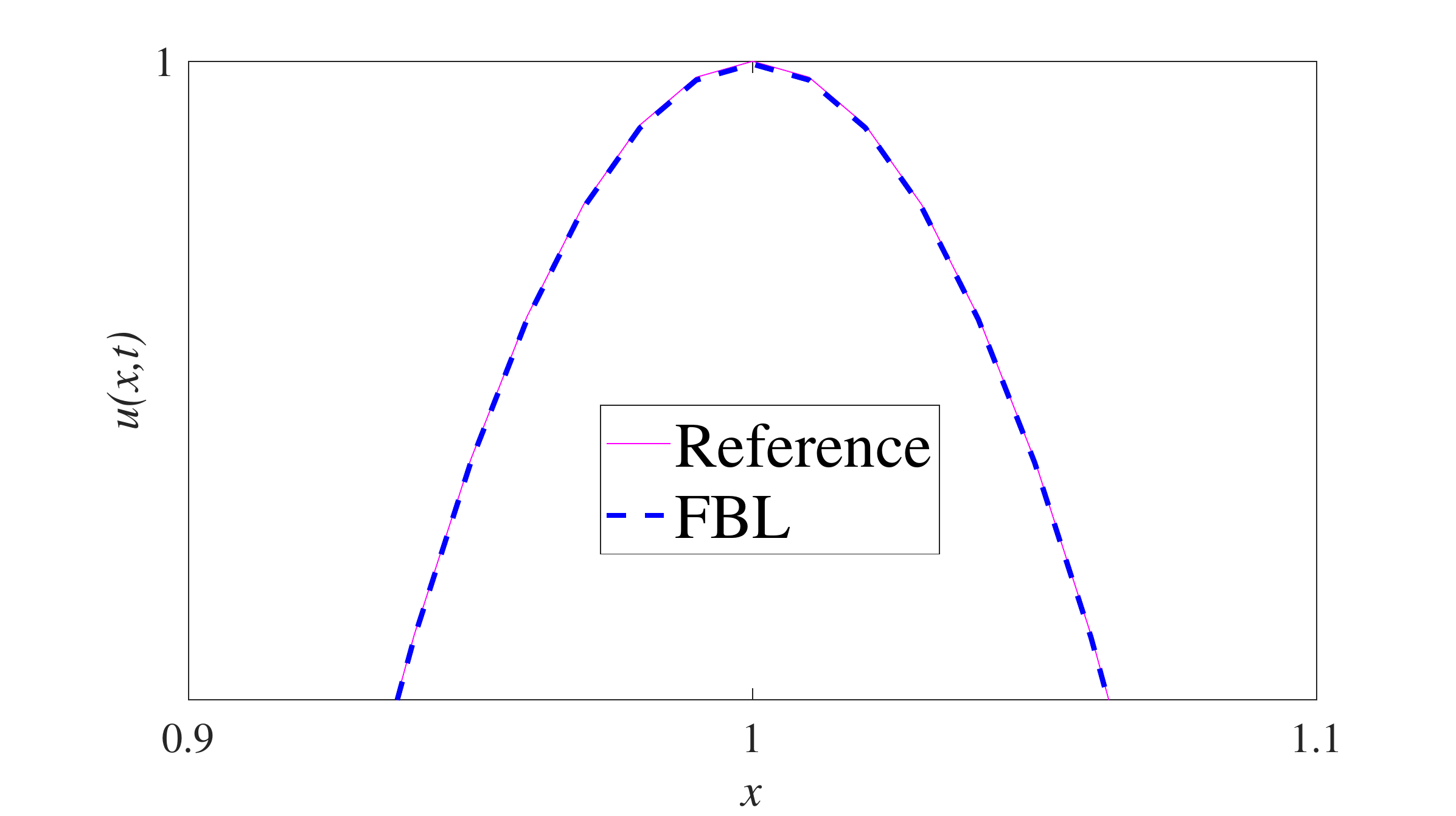} 
\includegraphics[clip, trim=0.5cm 0cm 1.7cm 0cm, width=0.3\linewidth]{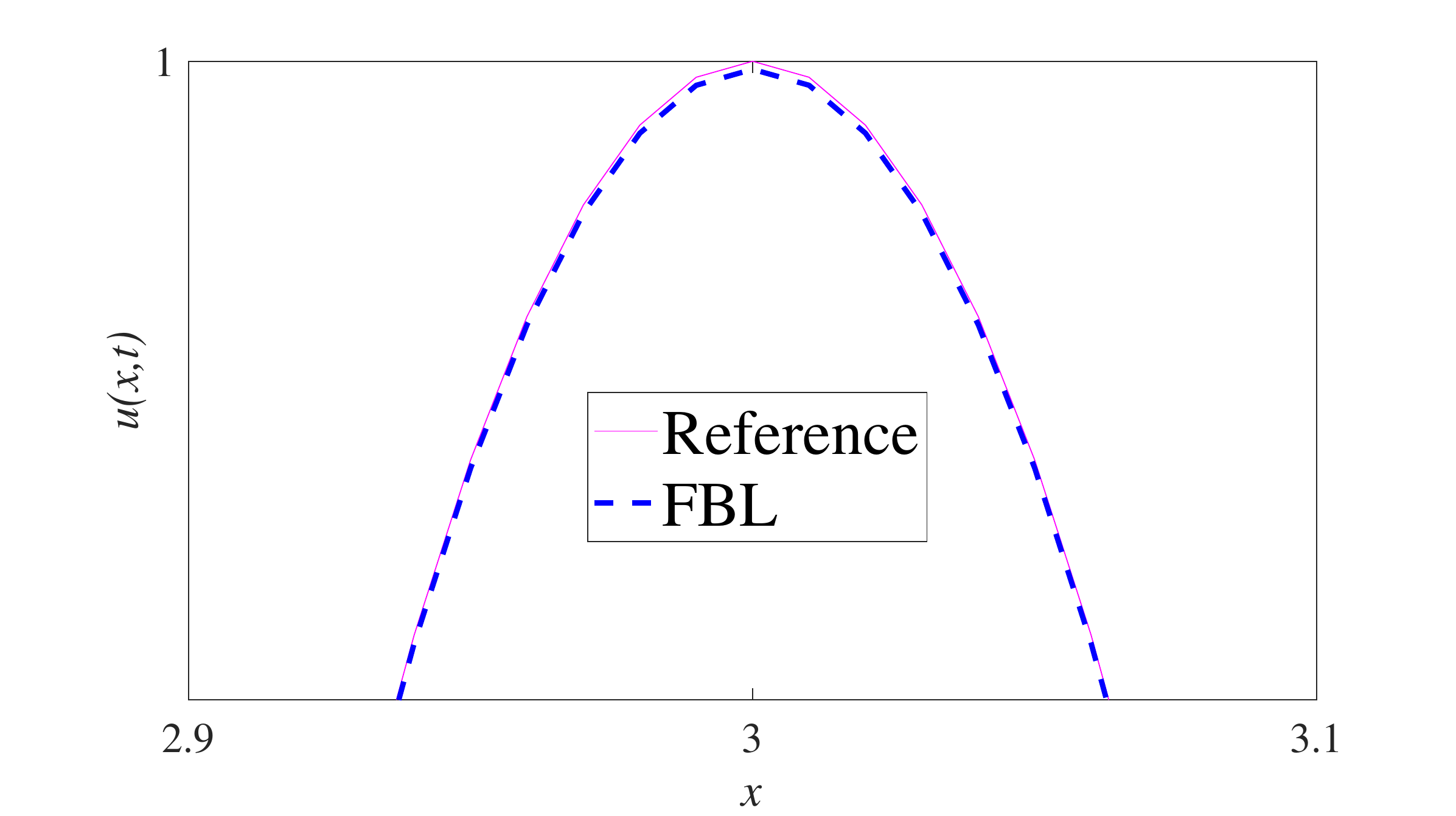} 
\includegraphics[clip, trim=0.5cm 0cm 1.7cm 0cm, width=0.3\linewidth]{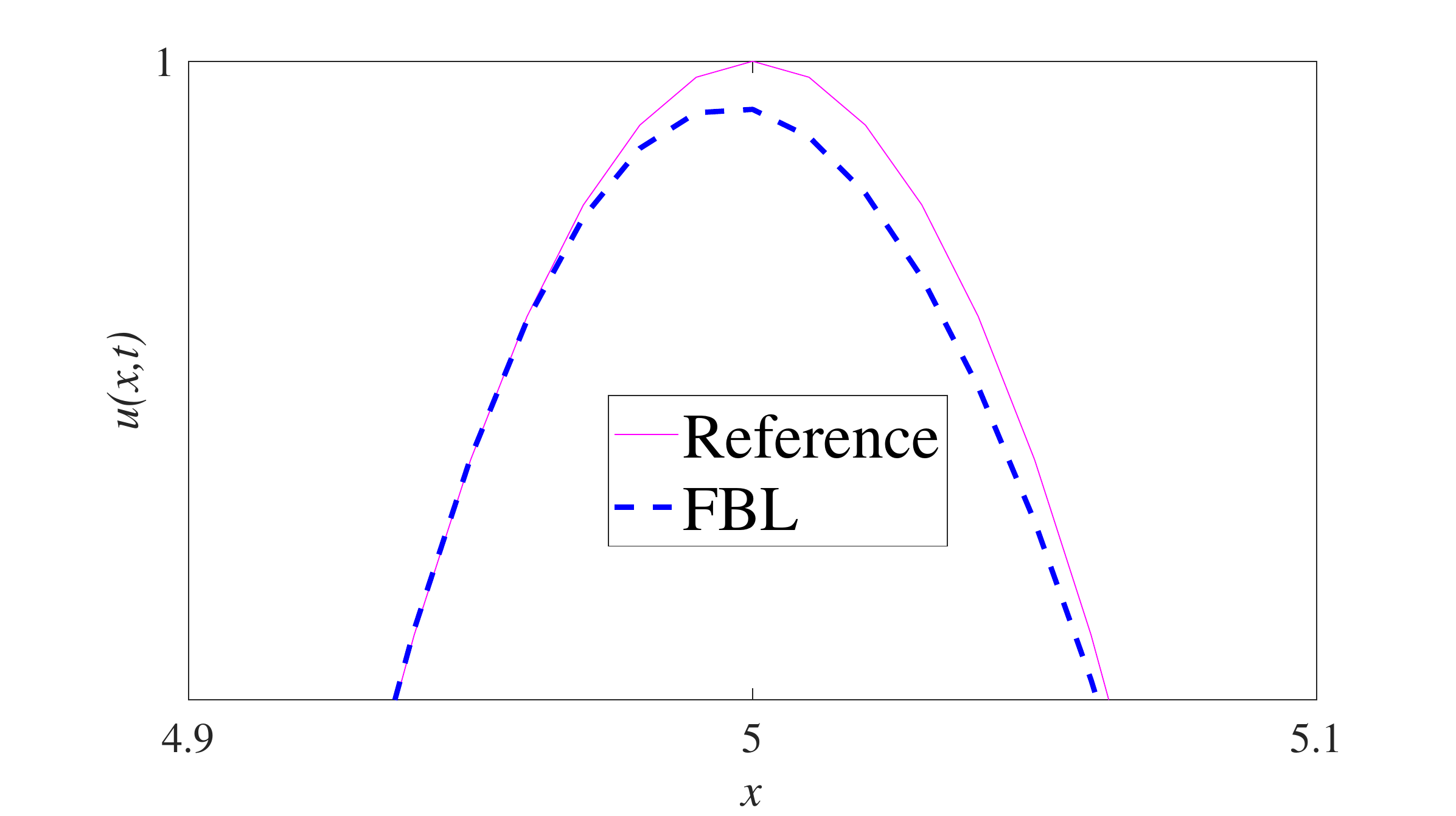}
\\
\includegraphics[width=1 \linewidth]{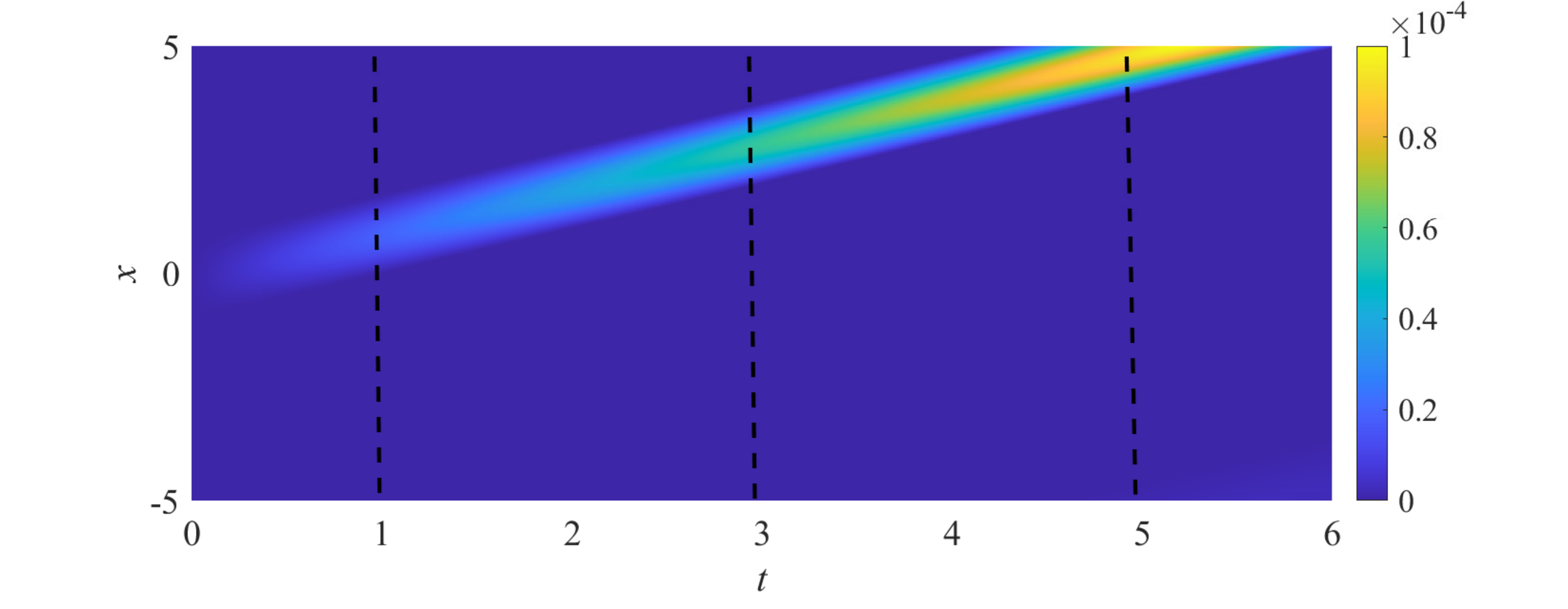}
\vspace{-0.4 in}
\caption{Model error $E_M$ for one-dimensional one-way wave (Example \ref{Ex: FBL one way right}). Top row: diffusive behaviour of FBL at the peak of the wave at $t=1$, $t=3$, and $t=5$. Bottom row: space-time point-wise error in the interior domain.}
\label{Fig:Decrease_Delta1}
\end{figure}

The total FBL error in the interior domain is comprised of three components: discretization, reflection, and model error, given as
$$E = E_N(M,\tau) + E_R(\delta,\bar\delta,\omega) + E_M(\epsilon).$$
The discretization error $E_N$ is due to the numerical scheme and can be bounded by choosing a sufficient number of collocation points $P$ and a proper time step $\tau$. We demonstrate a \textit{p}-refinement strategy for different formulations of FBLs in subsection \ref{sec: 1D Comparison with other formulations}. The reflection error $E_R$ is present because of wave reflections from the buffer layer into the interior domain when the wave reaches the boundary and depends on the FBL parameters $\{\delta,\bar\delta,\omega\}$. Different characterizations of FBL have different reflection errors. In subsection \ref{sec: Length of buffer layer}, we show the mitigation of $E_R$ by tuning the buffer layer parameters. The model error $E_M$ emerges because the FBL formulation is based on a variable-order fractional diffusion equation with the order being $1 + \epsilon$ in the interior domain. Although $\alpha(x)$ is very close to 1 in the interior domain, the equation is still a diffusion equation and thus we can observe a slight diffusive behaviour. The model error can be bounded by tuning the parameter $\epsilon$ in the definition of $\alpha(x)$. Yet, it becomes negligible when the propagation velocity is high and/or the interior domain is relatively small. Figure \ref{Fig:Decrease_Delta1} shows the model error for Example \ref{Ex: FBL one way right} by depicting the decrease of the peak of the numerical solution compared to the reference solution at $t=1,3,5$. We observe that the peak of the FBL solution is slightly below the reference solution and as time progresses this error becomes larger up to $O(10^{-4})$; see the bottom panel in Fig. \ref{Fig:Decrease_Delta1}.

% Figure \ref{Fig:Decrease_Delta1} shows the difference between the exact solution and the FBL solution for Example \ref{Ex: FBL one way right} in the settings of $(\delta,\bar{\delta},\omega)=(1,0.5,20)$ and $(M,\tau)=(500,10^{-3})$. The top row shows the the reference solution and FBL solution near the peak of the one-way wave at $t=1,3,5$. It can be observed that the peak of the FBL solution is lower than that of the reference solution, and the difference is getting bigger as $t$ increases, which indicates that the FBL solution has diffusive behaviour. The bottom row shows pointwise difference between the reference solution and the FBL solution. Here the difference is defined as $\widetilde{u}(x,t)=u^{exact}(x,t)-u(x,t)$. The slight diffusion can also be observed as we always have $\widetilde{u}(x,t)\geq0$.

% \end{remark}

%%%%%%%%%%%%%%%
\subsection{Length of buffer layer}
\label{sec: Length of buffer layer}
\begin{figure}[h!]
\centering 
\includegraphics[width=0.4 \linewidth]{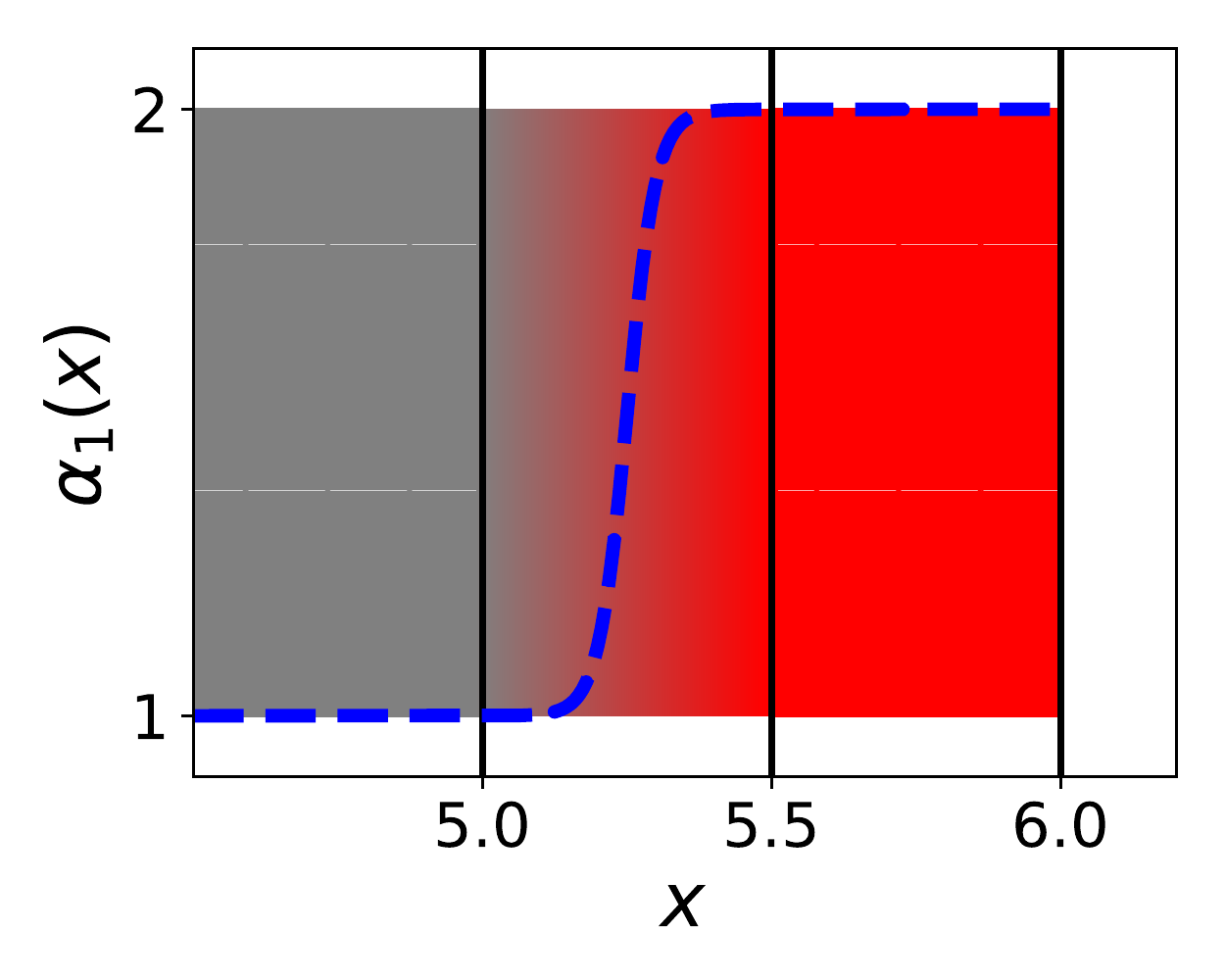}
\includegraphics[width=0.4 \linewidth]{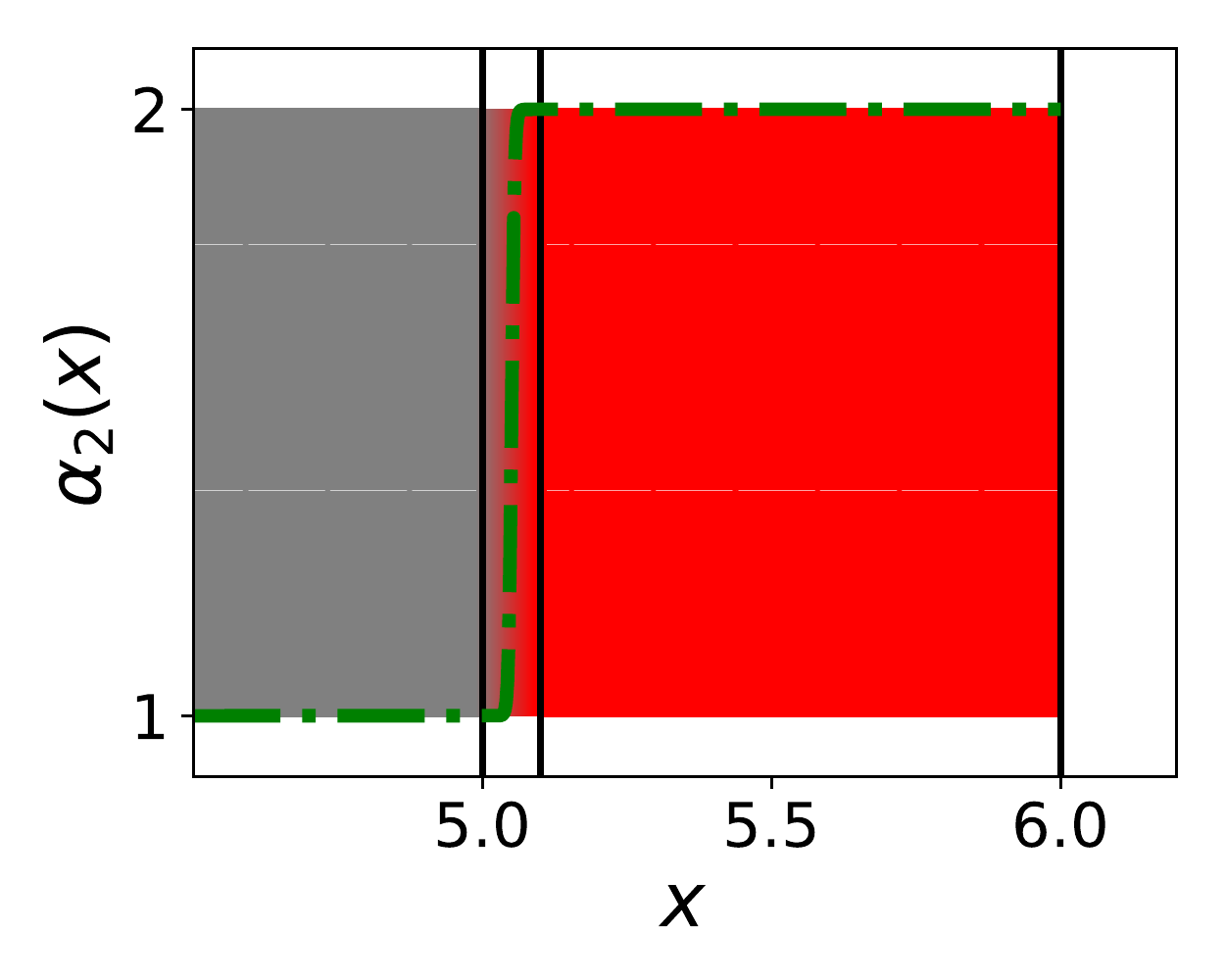} 
\vspace{-0.2 in}
\caption{Different characterizations of FBL. Left: long penetration region $\bar{\delta}=0.5$ with moderate variable-order $\omega=20$. Right: short penetration region $\bar{\delta}=0.1$ with steep variable-order $\omega=200$. The gray and red shaded areas are the advection and diffusion regions, respectively. The variable-order functions $\alpha_1(x)$ and $\alpha_2(x)$ are given in \eqref{eq:RightSmoothFunction} with buffer layer of size $\delta = 1$.}
\label{Fig:alf}
\end{figure}

\begin{figure}[h!]
\centering
{\scriptsize long penetration region $\bar{\delta}=0.5$ with moderate variable-order $\omega=20$}\\[-9 pt]
\noindent\rule{10cm}{0.7pt}
\\
\includegraphics[clip, trim=0.5cm 0cm 1.7cm 0cm, width=0.47\linewidth]{NumSol_RefSol_PMLSol_Delta1_t6.pdf}
\includegraphics[clip, trim=0.5cm 0cm 1.7cm 0cm, width=0.47\linewidth]{PointwiseErr_Delta1_t6.pdf}
\\
%%%%%%%%%%
{\scriptsize short penetration region $\bar{\delta}=0.1$ with steep variable-order $\omega=200$}\\[-9 pt]
\noindent\rule{10cm}{0.7pt}
\\
\includegraphics[clip, trim=0.5cm 0cm 1.7cm 0cm, width=0.47\linewidth]{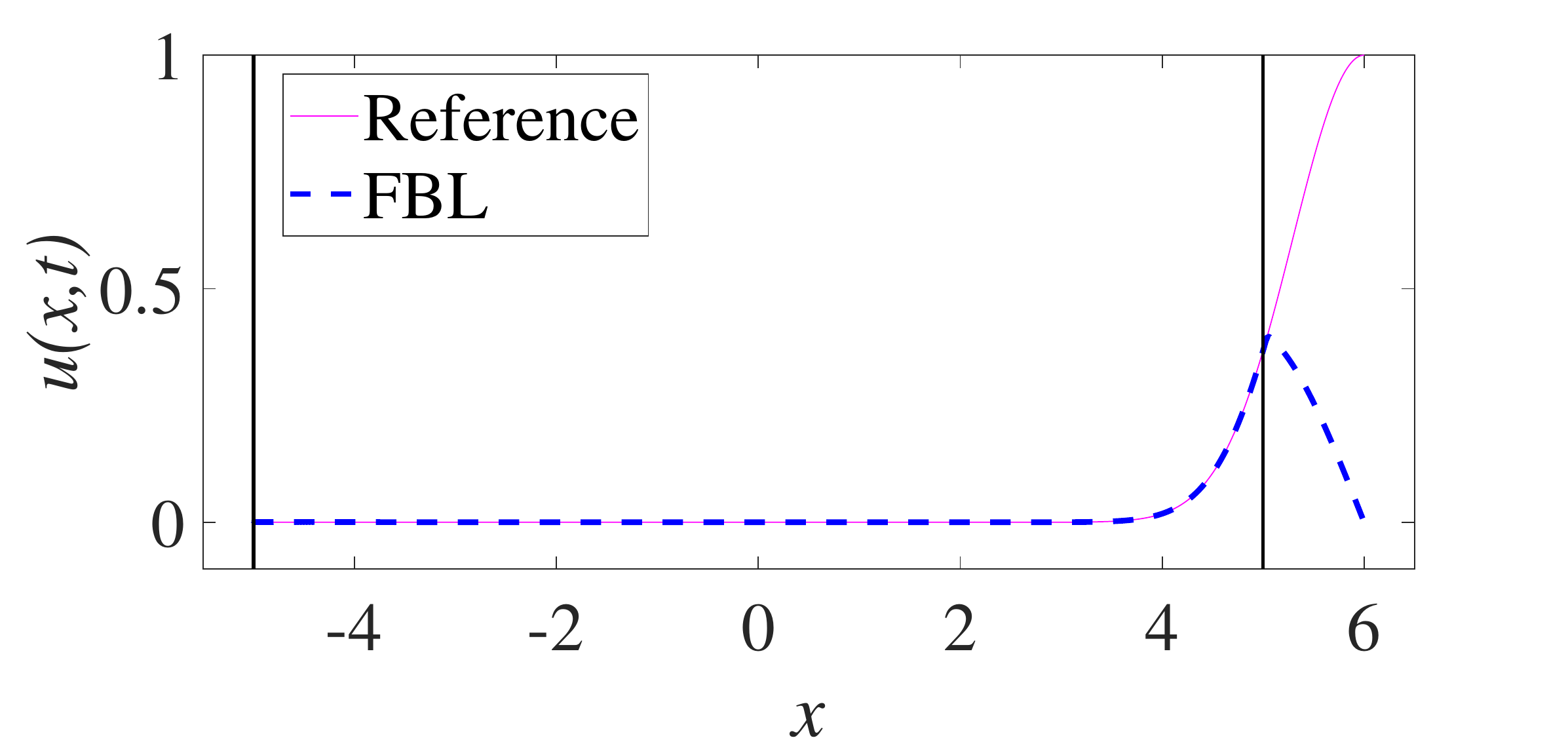}
\includegraphics[clip, trim=0.5cm 0cm 1.7cm 0cm, width=0.47\linewidth]{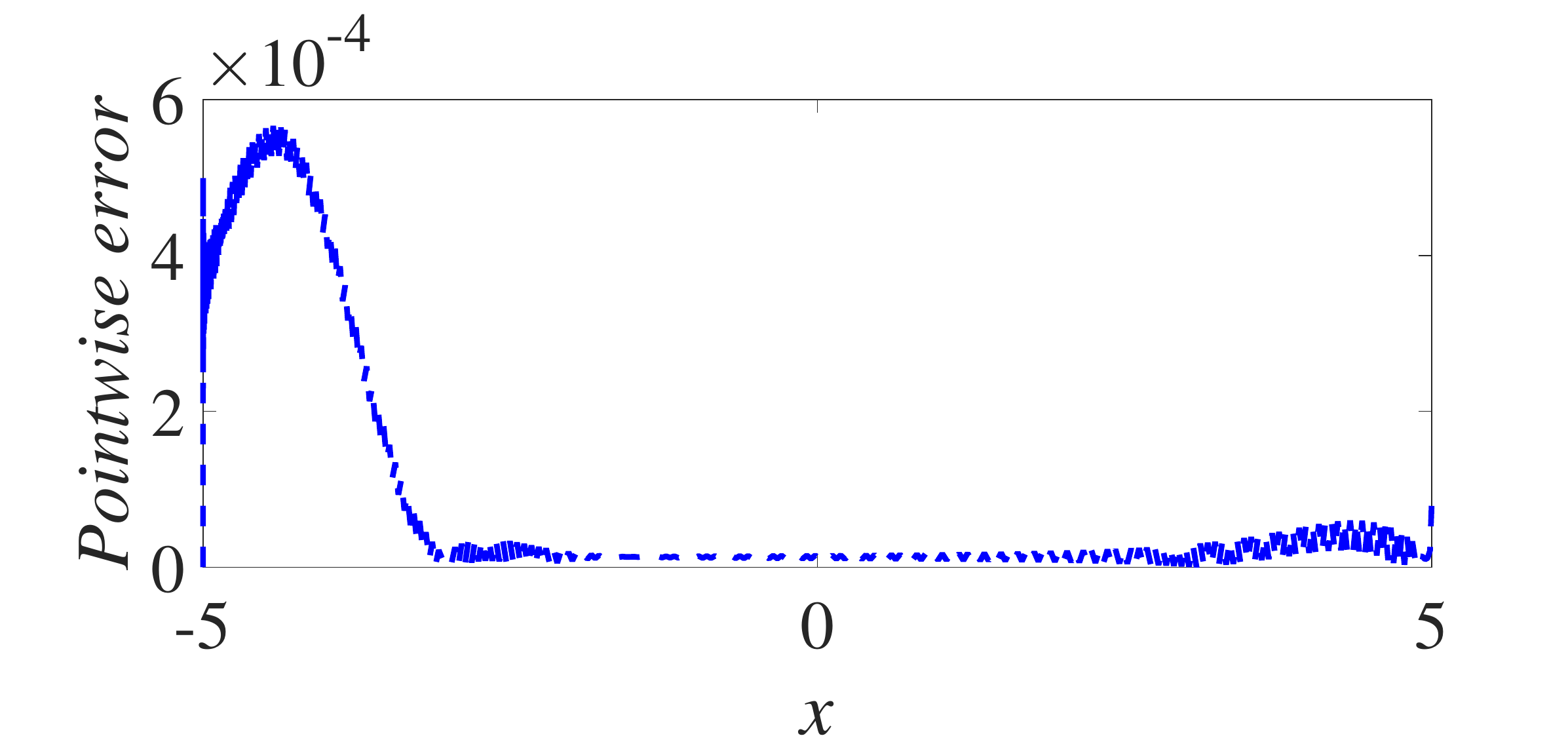}
%%%%%%%%%%
\vspace{-0.2 in}
\caption{One-dimensional one-way wave (Example \ref{Ex: FBL one way right}): FBL performance with different characterizations given in Fig. \ref{Fig:alf}. Top row: long penetration region $\bar{\delta}=0.5$ with moderate variable-order $\omega=20$. Bottom row: short penetration region $\bar{\delta}=0.1$ with steep variable-order $\omega=200$. Left column shows the FBL and reference solution at $t=6$. Right column shows the pointwise error at $t=6$.
The vertical black solid lines indicate the boundaries of the interior domain.
}
\label{Fig:NumSol_RefSol_PMLSol_compare}
\end{figure}

%%%%%%%%%%%%%%%
We study the performance of FBL by considering different descriptions of FBL. Figure \ref{Fig:alf} shows two different characterizations of the FBL in Example \ref{Ex: FBL one way right}: \textit{i}) long penetration region $\bar{\delta}=0.5$ with moderate variable-order $\alpha_1(x)$ in which $\omega=20$, and \textit{ii}) short penetration region $\bar{\delta}=0.1$ with steep variable-order $\alpha_2(x)$ in which $\omega=200$. The shorter penetration region requires a steeper slope in the variable-order function such that \eqref{eq:VO fun Advec1} always holds. The steep slope leads to a stronger reflection and thus a larger $E_R$ in the interior domain. This is shown in Fig. \ref{Fig:NumSol_RefSol_PMLSol_compare}, where we clearly observe the oscillatory behavior in the error and a much larger error forming at the left boundary. 

The key parameter in the performance of FBL is the length of penetration region $\bar\delta$, which dictates the slope of the variable-order function. The length of the buffer layer $\delta$ is less significant. 

%%%%%%%%%%%%%%%%%%%%%%%%%%%%%%%%%%%%%%%%
%%%%%%%%%%%%%%%%%%%%%%%%%%%%%%%%%%%%%%%

% \begin{figure}[h!]
% \centering
% {\scriptsize wider layer $\delta=1$ with moderate variable-order $\omega=20$}\\[-11 pt]
% \noindent\rule{10cm}{0.7pt}
% \\
% \includegraphics[clip, trim=0.5cm 0cm 1.7cm 0cm, width=0.47\linewidth]{NumSol_RefSol_PMLSol_Delta1_t6.pdf}
% \includegraphics[clip, trim=0.5cm 0cm 1.7cm 0cm, width=0.47\linewidth]{PointwiseErr_Delta1_t6.pdf}
% \\[-7 pt]
% %%%%%%%%%%
% {\scriptsize shorter layer $\delta=0.5$ with moderate variable-order $\omega=20$}\\[-11 pt]
% \noindent\rule{10cm}{0.7pt}
% \\
% \includegraphics[clip, trim=0.5cm 0cm 1.7cm 0cm, width=0.47\linewidth]{NumSol_RefSol_PMLSol_Delta05_t6.pdf}
% \includegraphics[clip, trim=0.5cm 0cm 1.7cm 0cm, width=0.47\linewidth]{PointwiseErr_Delta05_t6.pdf}
% %%%%%%%%%%
% \vspace{-0.2 in}
% \caption{One-dimensional one-way wave in Example \ref{Ex: FBL one way right} with penetration region $\bar{\delta}=\frac{\delta}{2}$. Top row: wider layer $\delta=1$ with moderate variable-order $\omega=20$. Bottom row: shorter layer $\delta=0.5$ with moderate variable-order $\omega=20$. Left column shows the FBL and reference solution at $t=6$. Right column shows the pointwise error at $t=6$.
% The vertical black solid lines indicate the boundaries of interior domain.
% }
% \label{Fig:NumSol_RefSol_PMLSol_compare_Delta}
% \end{figure}

%%%%%%%%%%%%%%%%%%%%%%%%%%%%%%%%%%%%%%%%%%%%
%%%%%%%%%%%%%%%%%%%%%%%%%%%%%%%%%%%%%%%%%%%%%%%%%%%%%%%%%%%

%%%%%%%%%%%%%%%
\subsection{Comparison with other formulations}
\label{sec: 1D Comparison with other formulations}
%%%%%%%%%%%%%%%
In order to gain a deeper insight into the performance of FBL and better characterize the corresponding errors, we compare different formulations of absorbing layer for integer- and fractional order one-dimensional one-way wave equation. The integer-order equation is given by \eqref{eq: Advec}, and has two different fractional counterparts with orders greater than one $(1+\epsilon)$ and less than one $(1-\epsilon)$. Both cases converge to the integer-order equation in the limiting case $\epsilon \rightarrow 0$ but describe different phenomena. If the order is greater than one, we have the fractional diffusion equation (FracDiff) with two boundary conditions and if the order is less than one, we have the fractional advection equation (FracAdv) with only one boundary condition.
\begin{align}
\label{eq:FracDiff} 
\text{FracDiff: \qquad\qquad\qquad}
\frac{\partial u}{\partial t}
    =|V|\prescript{RL}{x}{\mathcal{D}}_{x_R+\delta}^{1+\epsilon}u(x,t), \qquad\qquad\qquad\qquad
\\
\label{eq:FracAdv} 
\text{FracAdv: \qquad\qquad\qquad}
\frac{\partial u}{\partial t}
    =V\prescript{RL}{x_L-\delta}{\mathcal{D}}_{x}^{1-\epsilon}u(x,t),\qquad\qquad\qquad\qquad
\end{align}
Following the classical PML setting \cite{Johnson2007}, we add a damping term into the equations in the form of $-\sigma_{x}u$, where $\sigma_{x}$ is the damping function. The FBL, in contrast, transitions the advection to the diffusion via variable-order fractional differential operators. The FBL is only applicable to the integer-order \eqref{eq: Advec} and FracDiff \eqref{eq:FracDiff} and not to the FracAdv \eqref{eq:FracAdv}. Table \ref{Tab:Comparison} shows the different cases of equations describing advection dominating phenomena and different approaches in absorbing the wave. The contributing source of error is also mentioned in each case.
\begin{table}[h!]
{\scriptsize
\centering
\begin{tabular}{|l|c|c|c|}
\hline
&Integer-order Advection & Fractional Advection & Fractional Diffusion
\\ \hline &&&\\[-1em] 
& $\mathcal{L}_x = V \frac{\partial^{\alpha}}{\partial x^{\alpha}}, \,\, {\alpha} = 1$
& $\mathcal{L}_x=V\prescript{RL}{x_L-\delta}{\mathcal{D}}_{x}^{\alpha}, \,\, \alpha = 1-\epsilon$ 
& $\mathcal{L}_x=|V|\prescript{RL}{x}{\mathcal{D}}_{x_R+\delta}^{\alpha}, \,\, \alpha = 1+\epsilon$
%%%%
\\ &&&\\[-1em]  \hline\hline  &&&\\[-0.9em] 
%%%%
FBL  
& $\alpha \rightarrow \alpha(x)\in (1,2]$ 
& - 
& $\alpha \rightarrow \alpha(x)\in (1,2]$
%%%%
\\  &&&\\[-1em] 
%%%%
B.C.  
& $u(x_L) = u(x_R+\delta) = 0$ 
& -  
& $u(x_L) = u(x_R+\delta) = 0$
%%%%
\\  &&&\\[-1em] 
%%%%
error  
& $E_N + E_R + E_M$ 
& -  
& $E_N + E_R$
%%%%
\\ &&&\\[-1em]  \hline\hline &&&\\[-0.9em] 
%%%%
PML 
& $\mathcal{L}_x(u)-\sigma_x u $ 
& $\mathcal{L}_x(u)-\sigma_x u$ 
& $\mathcal{L}_x(u)-\sigma_x u$ 
%%%%
\\  &&&\\[-1em] 
%%%%
B.C.  
& $u(x_L) = 0$ 
& $u(x_L) = 0$ 
& $u(x_L) = u(x_R+\delta) = 0$
%%%%
\\  &&&\\[-1em] 
%%%%
error  
& $E_N + E_R $ 
& $E_N + E_R $ 
& $E_N + E_R $
%%%%
\\ \hline
%%%%
\end{tabular}
\caption{Different formulations of absorbing layer for one-way wave of the form $\frac{\partial}{\partial t}u
= V \mathcal{L}_x (u)$.}
\label{Tab:Comparison}
}
\end{table}

We compare these different formulations by obtaining the $L_{\infty}$ error in the interior domain. Similar to Example \ref{Ex: FBL one way right}, we only consider the right propagation advection with the velocity $V=-1$ and a smooth initial condition of the form $u_{0}(x)=e^{-x^2}$. We note that similar results can also be obtained for the case of left propagation advection. The interior domain is $(x_L,x_R)=(-5,5)$ and we append a buffer layer of length $\delta=1$ next to the right boundary. Thus, $x\in [-5,6]$. The adopted numerical scheme for all of these formulations is based on the spectral collocation method with $P$ collocation points and the Crank-Nicolson formula with time step $\tau=10^{-3}$. The error is obtained against the reference solution $u(x,t)=e^{-(x-t)^2}$. In the FBL formulations, we consider the smooth variable-order function $\alpha(x)$ given by \eqref{eq:RightSmoothFunction} with parameters $(\epsilon,\bar{\delta},\omega)=(10^{-5},0.5,20)$. In the PML formulations, the damping function $\sigma_x$ is defined similarly as $\alpha(x)$ such that $\sigma_x = \alpha(x) - 1$ with parameters $(\epsilon,\bar{\delta},\omega)=(0,0.5,20)$.

\begin{figure}[h!]
\centering 
\includegraphics[width=0.9 \linewidth]{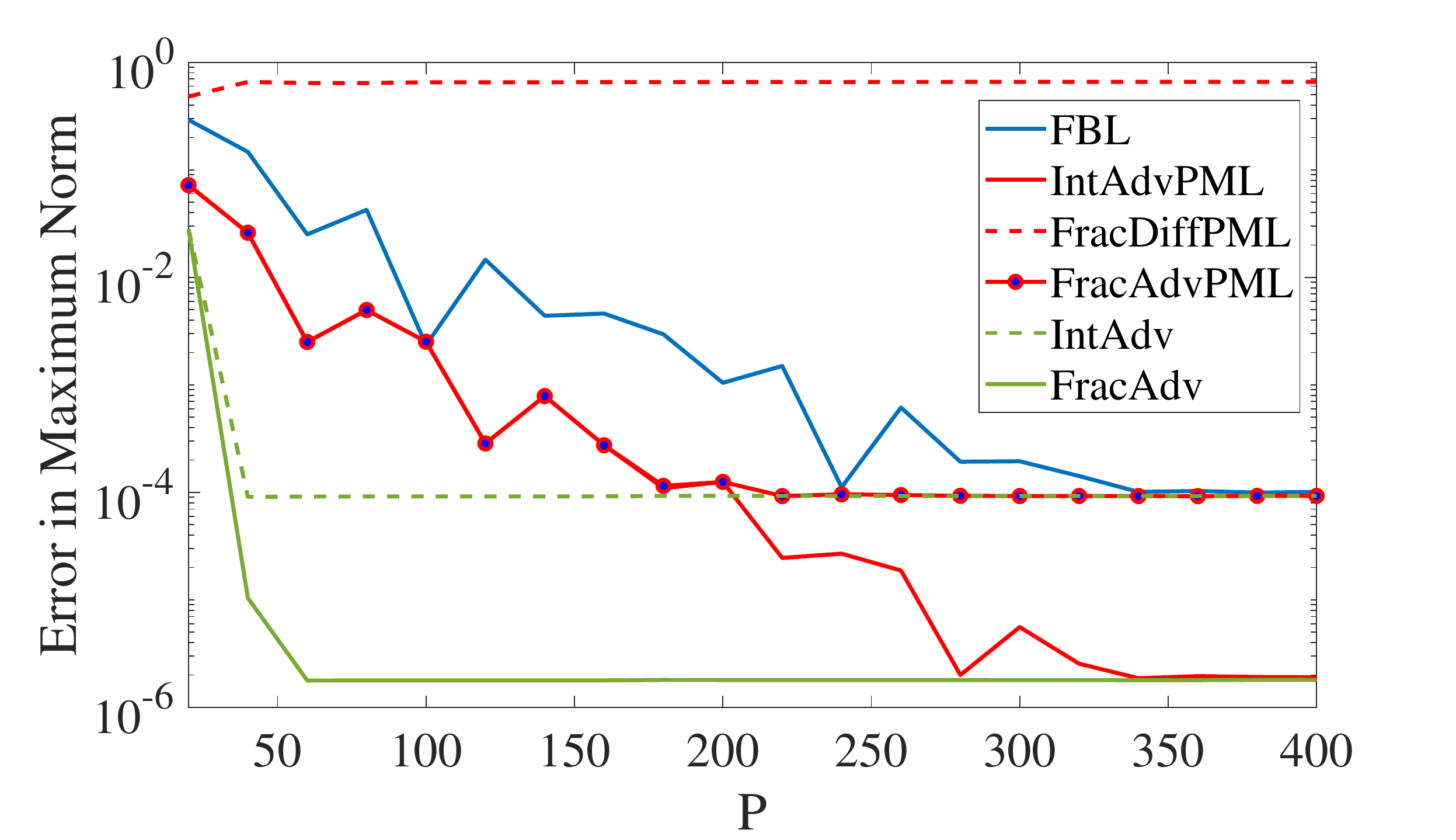} 
\vspace{-0.2 in}
\caption{P-refinement results for comparison between different formulations of PML and FBL, given in Table \ref{Tab:Comparison}.}
\label{Fig:PRefinement}
\end{figure}

Figure \ref{Fig:PRefinement} shows the \textit{p}-refinement results. For the integer-order case, the FBL (solid blue line) and PML (IntAdv-PML, solid red line) have almost similar rate of convergence. However, FBL has slightly larger error that is eventually bounded by the model error $E_M = O(10^{-4})$, which is in agreement with Fig. \ref{Fig:Decrease_Delta1}. Therefore, for integer-order advection equation, the PML formulation outperforms the FBL formulation. This is not true for the case of the fractional diffusion equation, where PML totally fails (FracDiff-PML, dashed red line). This failure is due to the additional homogeneous Dirichlet boundary condition at the right boundary for the fractional diffusion equation. The FBL formulation, however, is effective in this case (solid blue line). Thus, FBL will be a better choice if we consider an advection dominating fractional diffusion equation.  

% The FBL is less accurate than the classical PML for the integer-order advection equation as there is an additional model error in FBL, which can be verified by the similar behavior between FBL and the considered fractional advection equation. It can also be observed that imposing homogeneous Dirichlet boundary conditions to  advection dominant fractional diffusion equation can cause reflections, even in the case with its PML formulation. However, FBL can remove this deficiency. Therefore, FBL will be a better choice if we consider a advection dominant fractional diffusion equation. 

%%%%%%%%%%%%%%%%%%%%%%%%%%%%%%%%%%%%%%%
%%%%%%%%%%%%%%%%%%%%%%%%%%%%%%%%%%%%%%%

%%%%%%%%%%%%%%%%%%%%%%%%%%%%%
\section{FBL: One-dimensional two-way waves}
\label{Sec:Wave_1D}
%%%%%%%%%%%%%%%%%%%%%%%%%%%%%

We now consider the one-dimensional two-way wave on the bounded domain $\Omega = (x_L,x_R)$. The governing equation is given as  
\begin{equation}
\label{eq:WaveEq}
\frac{\partial^{2}}{\partial t^{2}}u(x,t)
=c^2\,\frac{\partial^{2}}{\partial x^{2}}u(x,t), 
\end{equation} 
with known initial conditions $u(x,0)=u_0(x)$ and $\frac{\partial}{\partial t}u(x,0)=\varphi(x)$, where the non-negative real coefficient $c$ is the wave propagation speed.  

%%%%%%%%%%%%%%%
\subsection{An existing fractional approach}
%%%%%%%%%%%%%%%
The first fractional approach in setting absorbing boundary layers for the one-dimensional two-way wave equation was introduced in \cite[Sec. 3.2]{Zhao@JCP2015}. It proposes the variable-order time fractional wave equation over the extended domain $[x_L-\delta,x_R+\delta]$ as 
\begin{equation}
\label{eq:FracWaveCaputo}  
       \prescript{C}{0}{\mathcal{D}}_{t}^{\gamma(x,t)}u(x,t)
       =c^{2}\,\frac{\partial^{2}}{\partial x^{2}}u(x,t), 
\end{equation} 
subject to initial conditions $u(x,0)=u_{0}(x)$ and $\frac{\partial}{\partial t}u(x,0)=\varphi(x)$, and homogeneous Dirichlet boundary conditions. The time fractional derivative is in the Caputo sense and the variable-order is $\gamma(x,t)\in(1,2]$. When $\gamma(x,t)=2$, \eqref{eq:FracWaveCaputo} coincides with the wave equation and when $\gamma(x,t)\rightarrow1^{+}$, it converges to the diffusion equation. The variable-order function is chosen such that it decreases from 2 to $1+\epsilon$ ($\epsilon>0$ is sufficiently small) when the wave reaches the boundary, and thus it is necessary for the variable-order function to be time-dependent. We shall demonstrate this approach via the following example. 

\begin{example}\label{Ex: FBL two way 1D Xuan}
To verify the time-dependence of the variable-order functions in  \eqref{eq:FracWaveCaputo}, we consider the absorbing boundary layer based on \eqref{eq:FracWaveCaputo} with the propagation speed $c=1$ and the smooth initial conditions $u_{0}(x)=2e^{-(x-10)^2}$, $\varphi(x)=0$. We let $(x_L,x_R)=(5,15)$ be the interior domain with buffer layers of the length $\delta=5$ next to the boundaries.
\end{example}

In Example \ref{Ex: FBL two way 1D Xuan}, in addition to adopting the time-dependent variable-order function given by \cite{Zhao@JCP2015}, we consider a time-independent variable-order function, which ensures the presence of diffusion in the buffer layers. Here, we apply the central difference method and the second-order formula proposed by \cite{Zhao@JCP2015} to numerically solve \eqref{eq:FracWaveCaputo}.
Note that the exact solution to the two-way wave equation \eqref{eq:WaveEq} in this case is  $u^{exact}(x,t)=\frac{u_{0}(x+t)+u_{0}(x-t)}{2}$. We regard it as the reference solution to compare the results in the interior domain.

\begin{figure}
\centering
\includegraphics[width=0.8\linewidth]{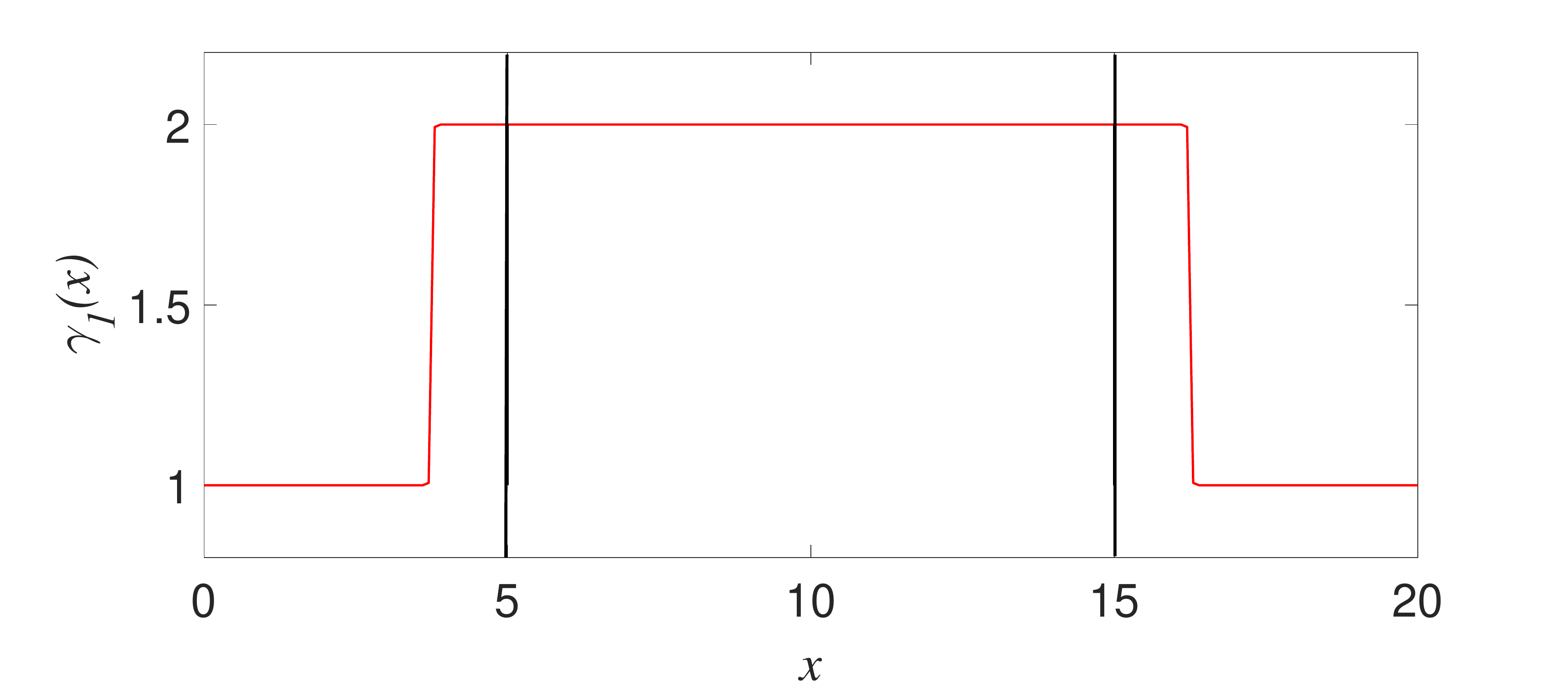}
\vspace{-0.2 in}
%%%%%%%%%%
\caption{
% Time fractional one-dimensional two-way wave (Example \ref{Ex: FBL two way 1D Xuan}). 
The time-independent variable-order function $\gamma_{1}(x)$, given by Eq. \eqref{eq:TimeIndependentOrd} with $\delta=5$. The vertical black solid lines indicate the boundaries of interior domain. }\label{Fig:TimeIndependentOrd}
\end{figure}

\begin{figure}
\centering
\qquad $t=0$\\[-0.8 pt]
% \noindent\rule{4cm}{0.7pt}
% \\
\includegraphics[clip, trim=0.5cm 0cm 1.7cm 0.46cm, width=0.75\linewidth]{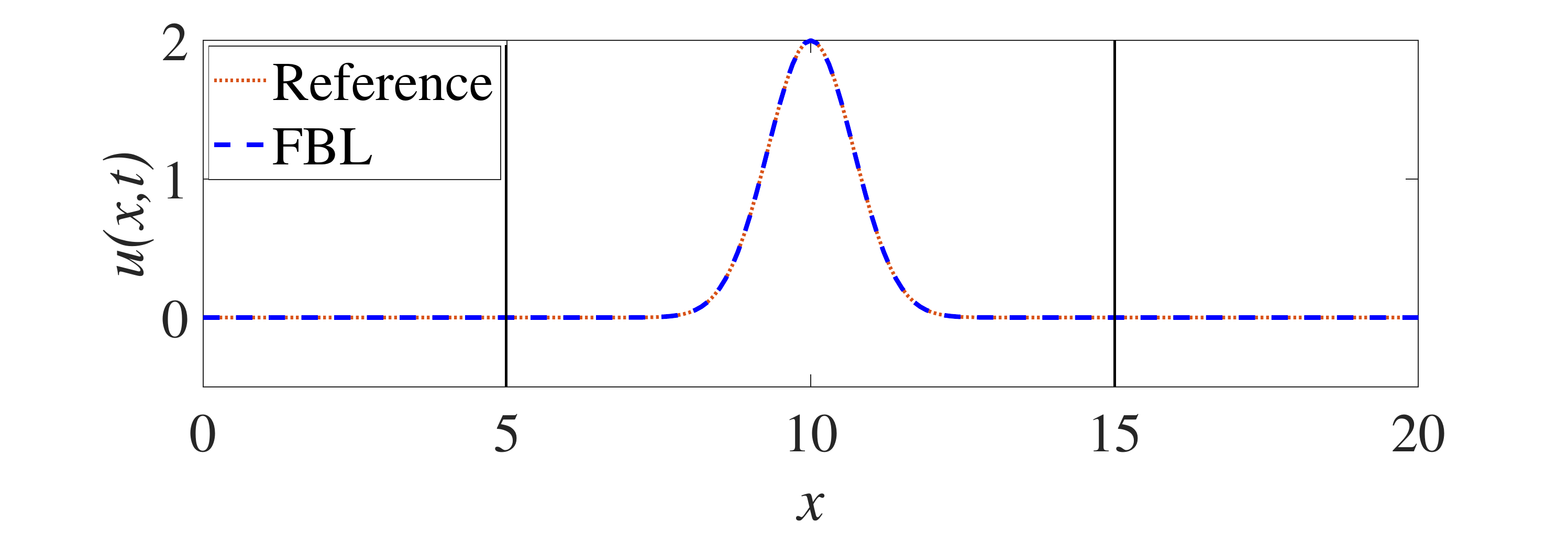}
\\
\qquad $t=3$\\[-0.8 pt]
% \noindent\rule{4cm}{0.7pt}
% \\
\includegraphics[clip, trim=0.5cm 0cm 1.7cm 0.46cm, width=0.75\linewidth]{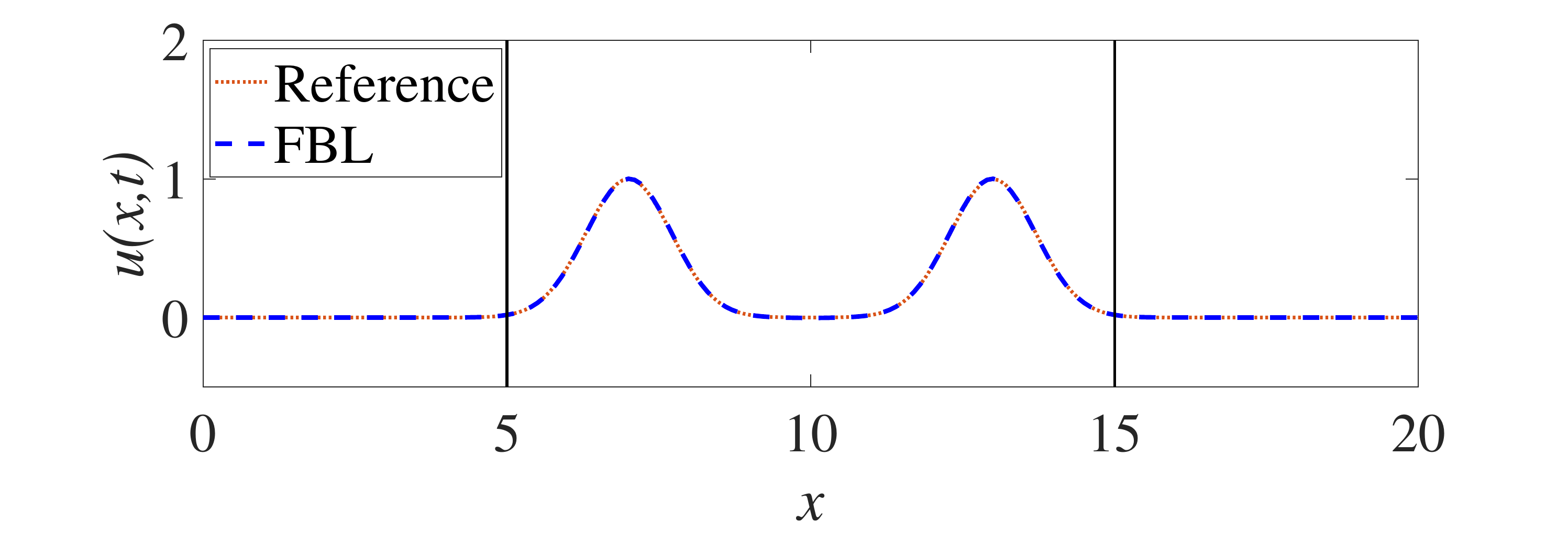}
\\
\qquad $t=6$\\[-0.8 pt]
% \noindent\rule{4cm}{0.7pt}
% \\
\includegraphics[clip, trim=0.5cm 0cm 1.7cm 0.46cm, width=0.75\linewidth]{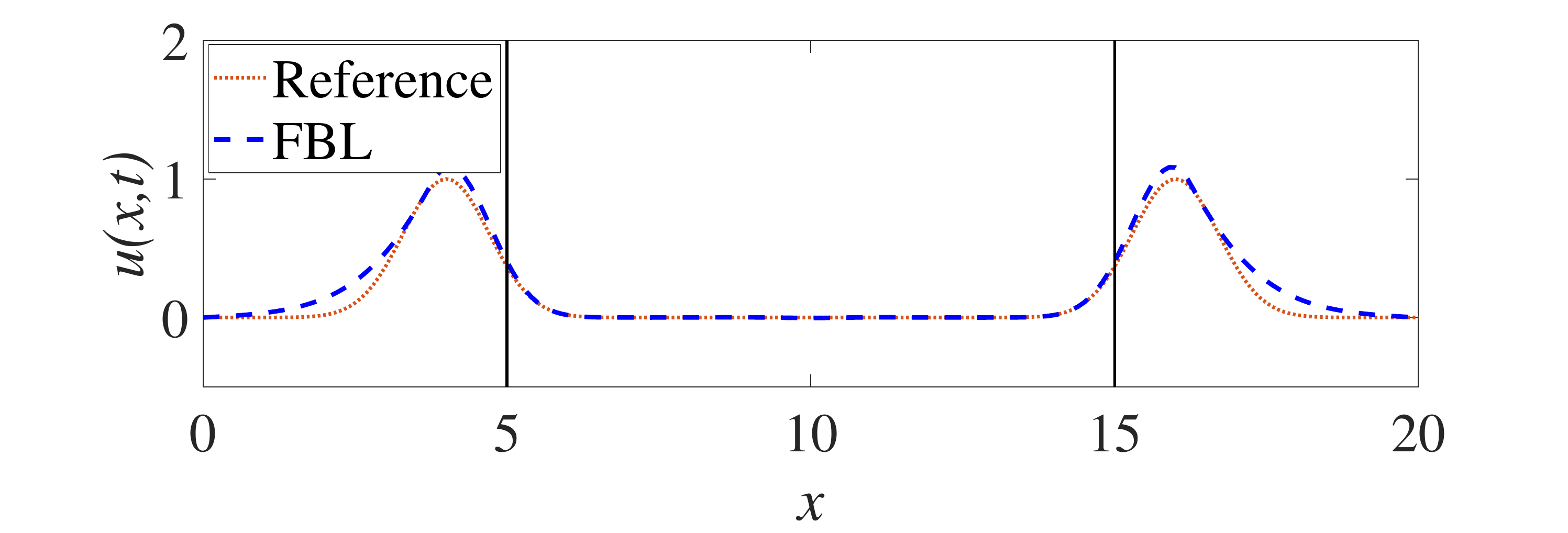}
\\
\qquad $t=9$\\[-0.8 pt]
% \noindent\rule{4cm}{0.7pt}
% \\
\includegraphics[clip, trim=0.5cm 0cm 1.7cm 0.46cm, width=0.75\linewidth]{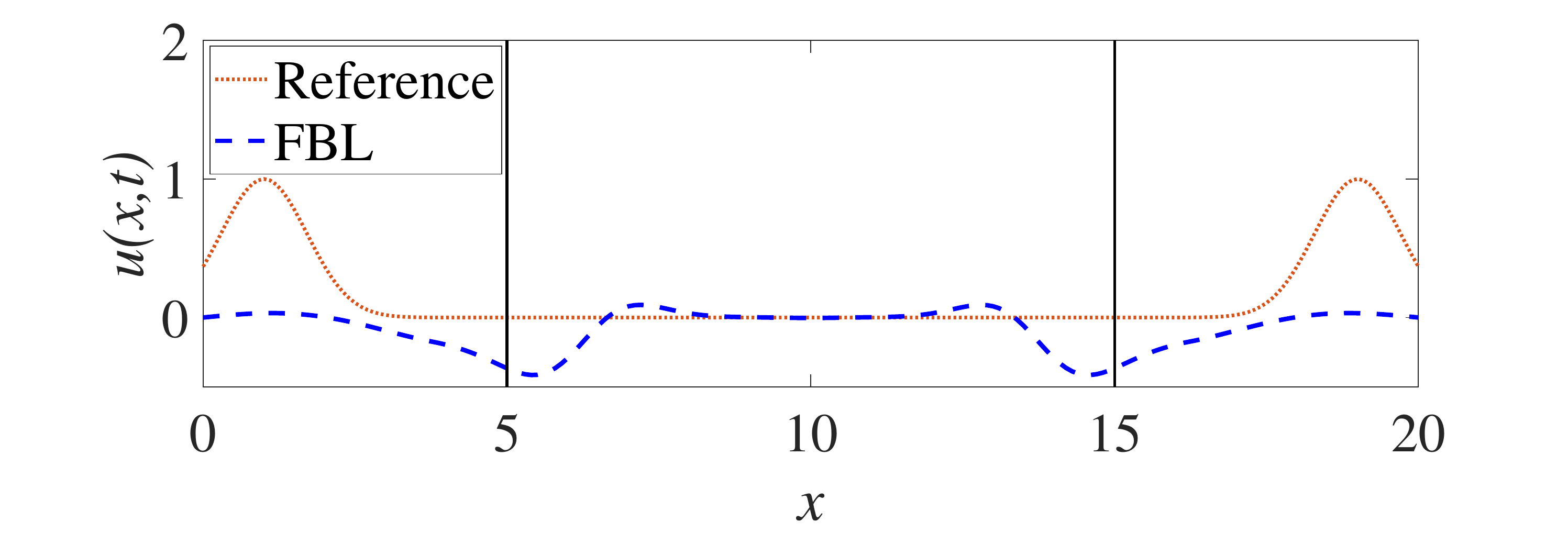}
\vspace{-0.2 in}
%%%%%%%%%%
\caption{Time fractional one-dimensional two-way wave (Example \ref{Ex: FBL two way 1D Xuan}): time-independent variable-order function $\gamma_{1}(x)$, given by \eqref{eq:TimeIndependentOrd}. The vertical black solid lines indicate the boundaries of interior domain. Reflections can be observed at $t=9$.
}\label{Fig:FracLayer4CaputoWave_SGL_1D}
\end{figure}

\begin{figure}
\centering
$t=0$\\[-13 pt]
\noindent\rule{4cm}{0.7pt}
\\
\includegraphics[clip, trim=0.5cm 0cm 1.7cm 0cm, width=0.47\linewidth]{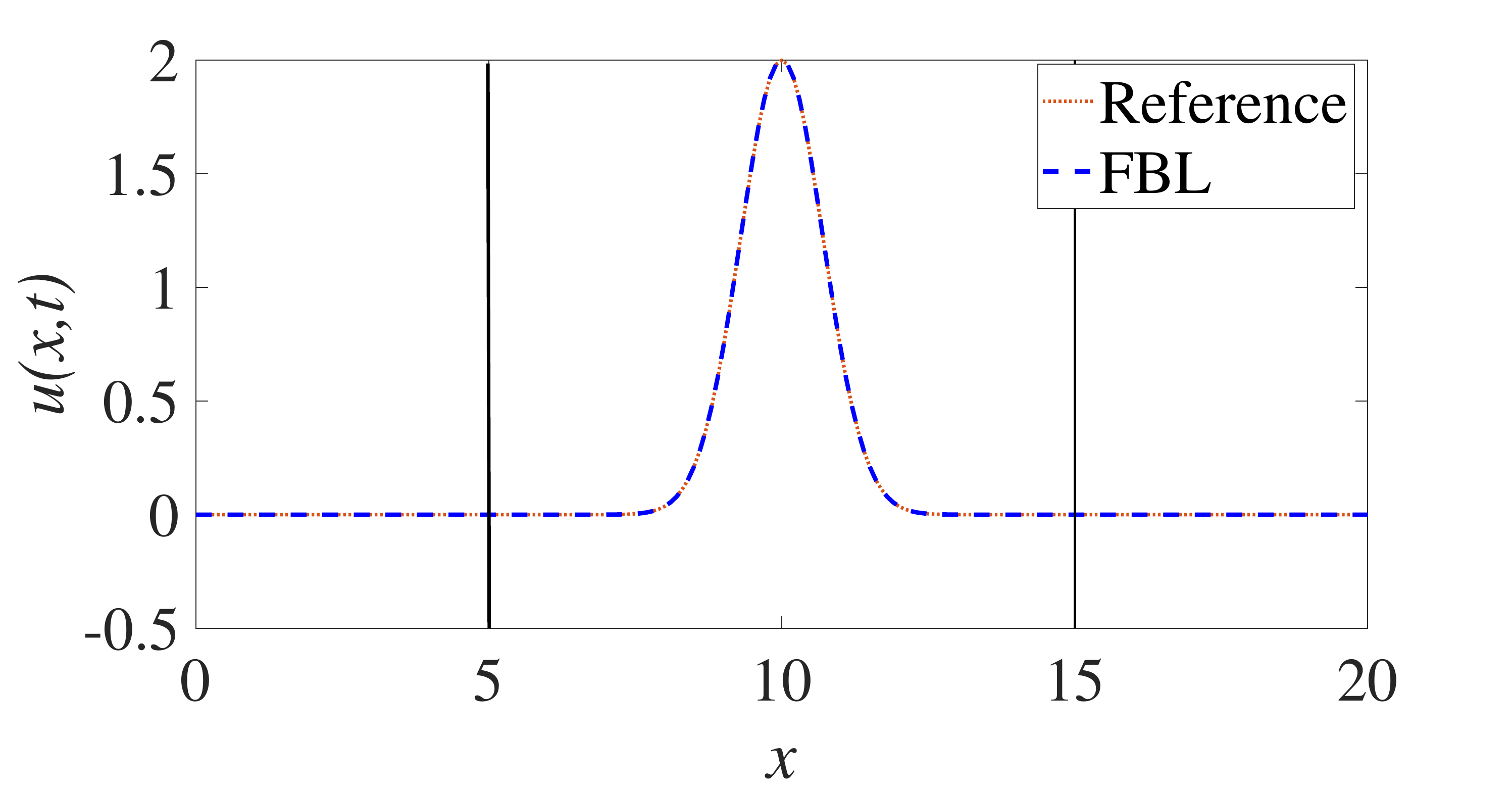}
\includegraphics[clip, trim=0.5cm 0cm 1.7cm 0cm, width=0.47\linewidth]{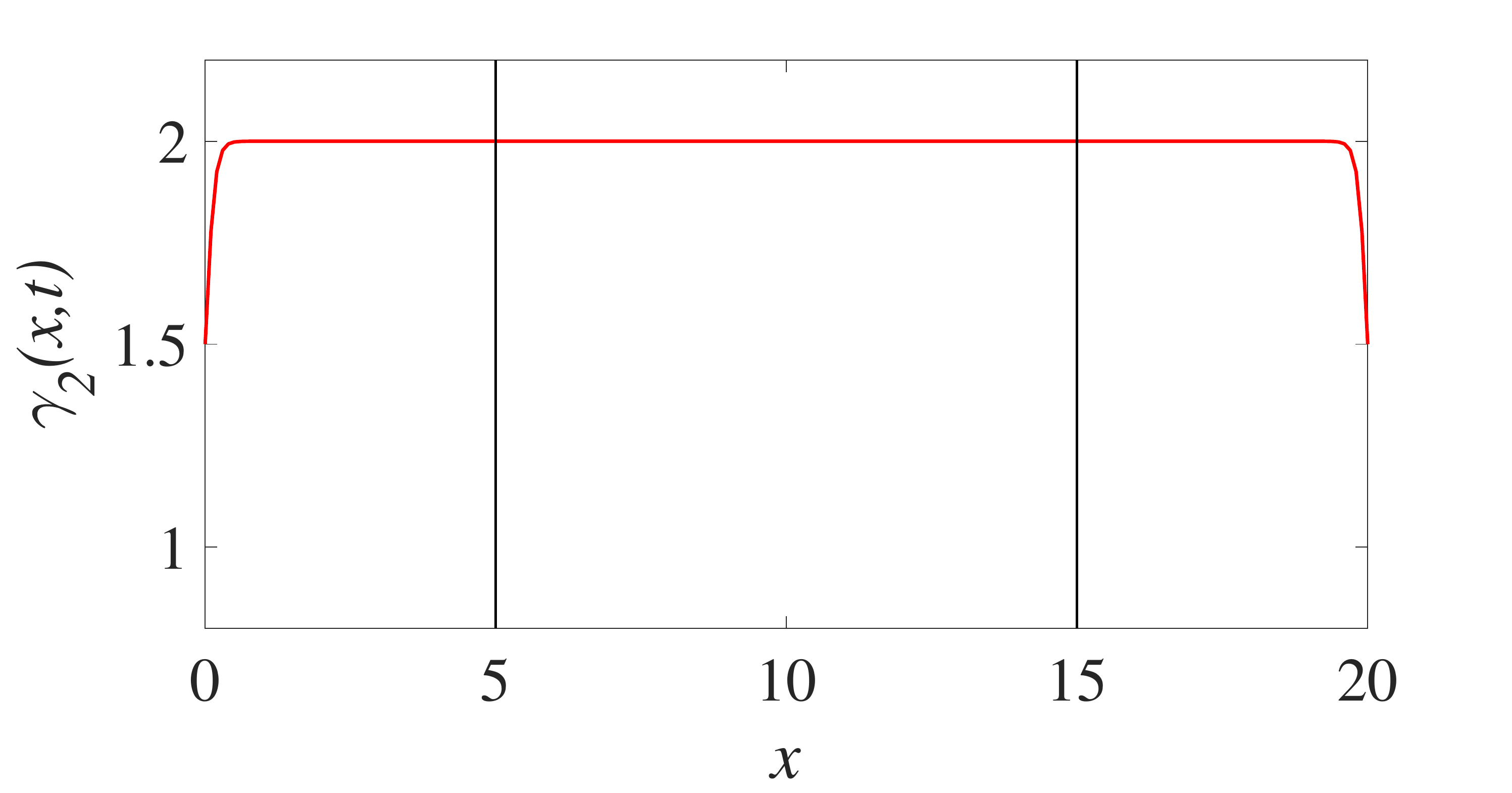}
\\[-8 pt]
$t=3$\\[-13 pt]
\noindent\rule{4cm}{0.7pt}
\\
\includegraphics[clip, trim=0.5cm 0cm 1.7cm 0cm, width=0.47\linewidth]{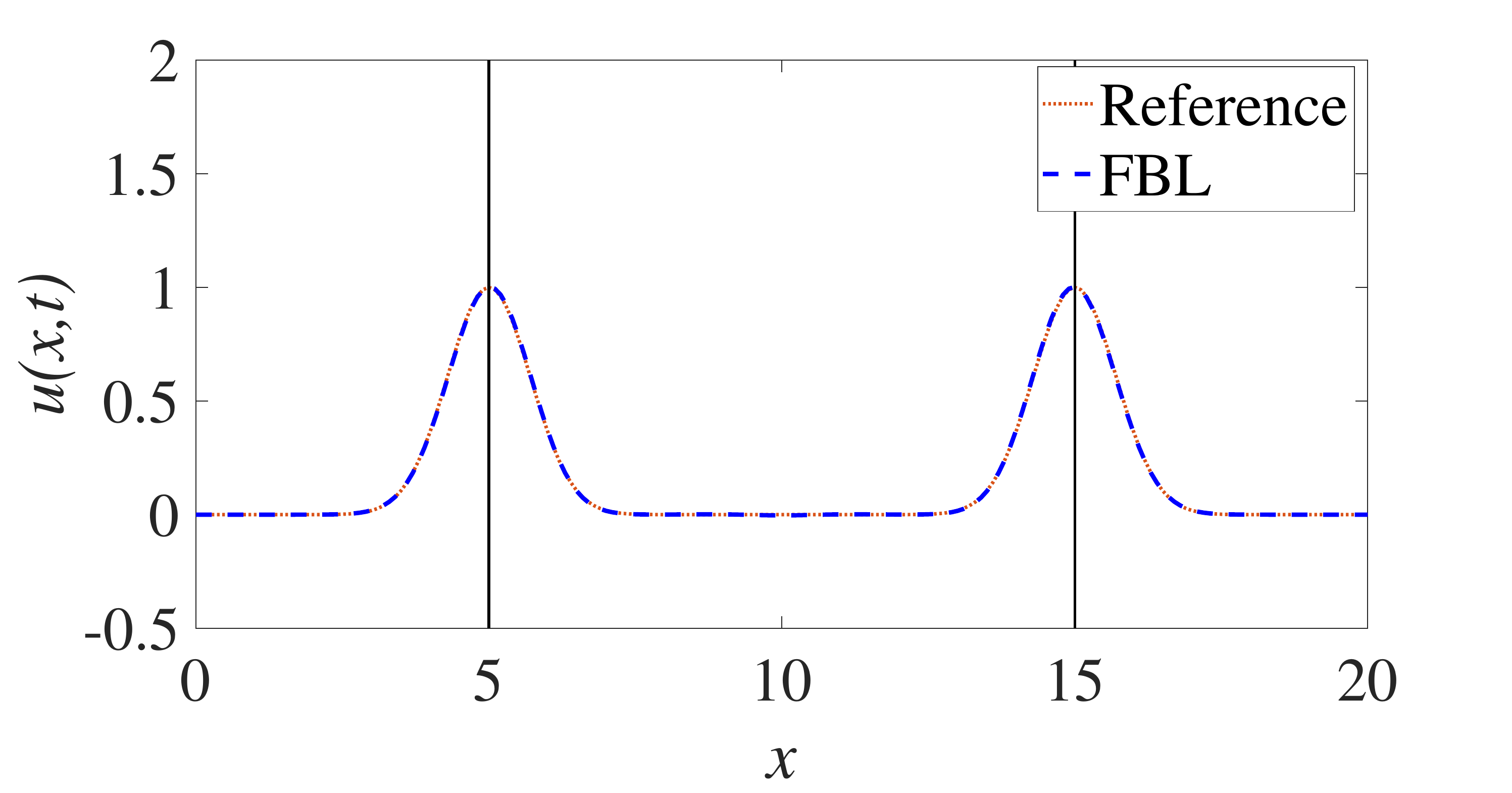}
\includegraphics[clip, trim=0.5cm 0cm 1.7cm 0cm, width=0.47\linewidth]{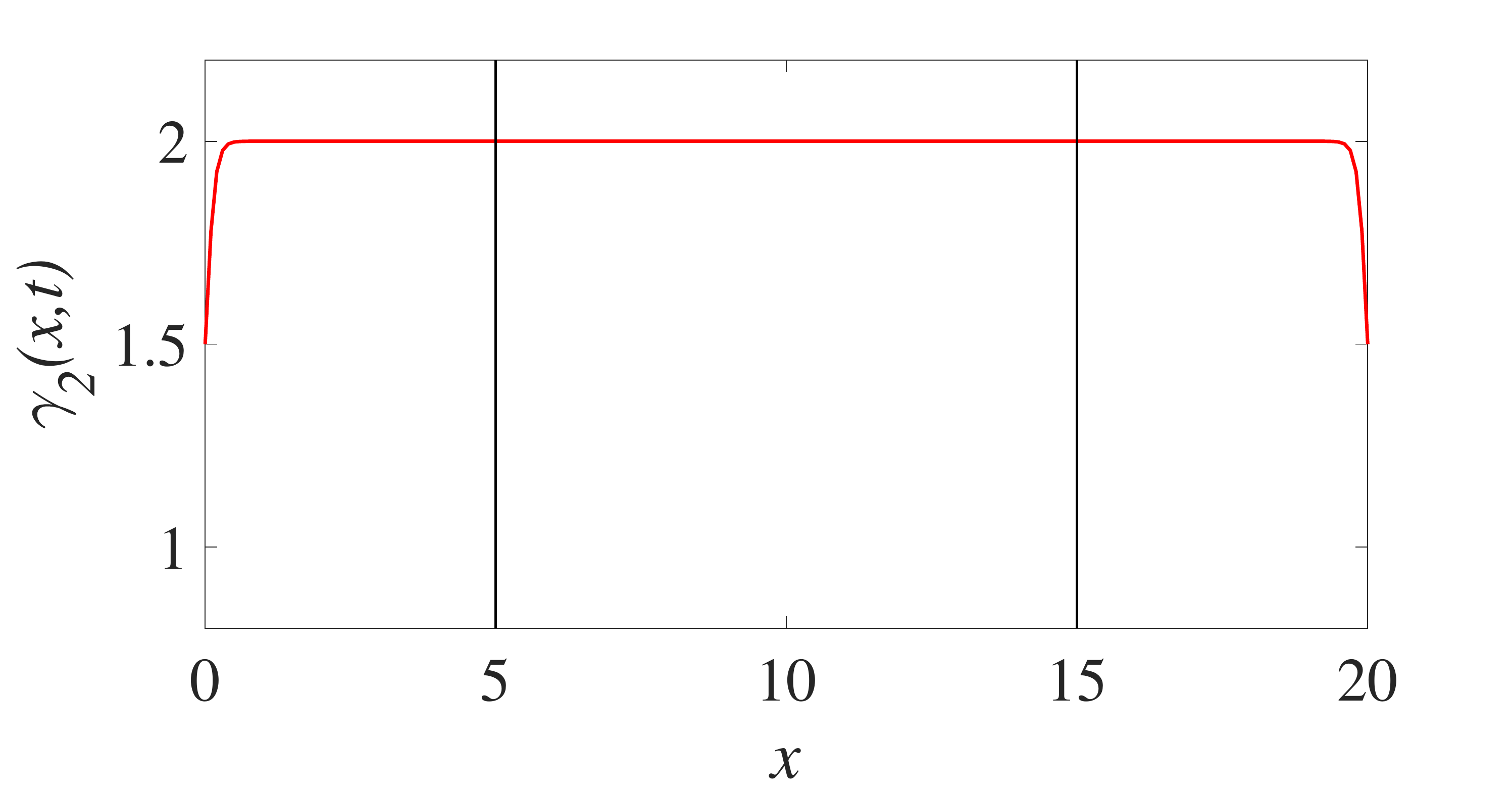}
\\[-8 pt]
$t=6$\\[-13 pt]
\noindent\rule{4cm}{0.7pt}
\\
\includegraphics[clip, trim=0.5cm 0cm 1.7cm 0cm, width=0.47\linewidth]{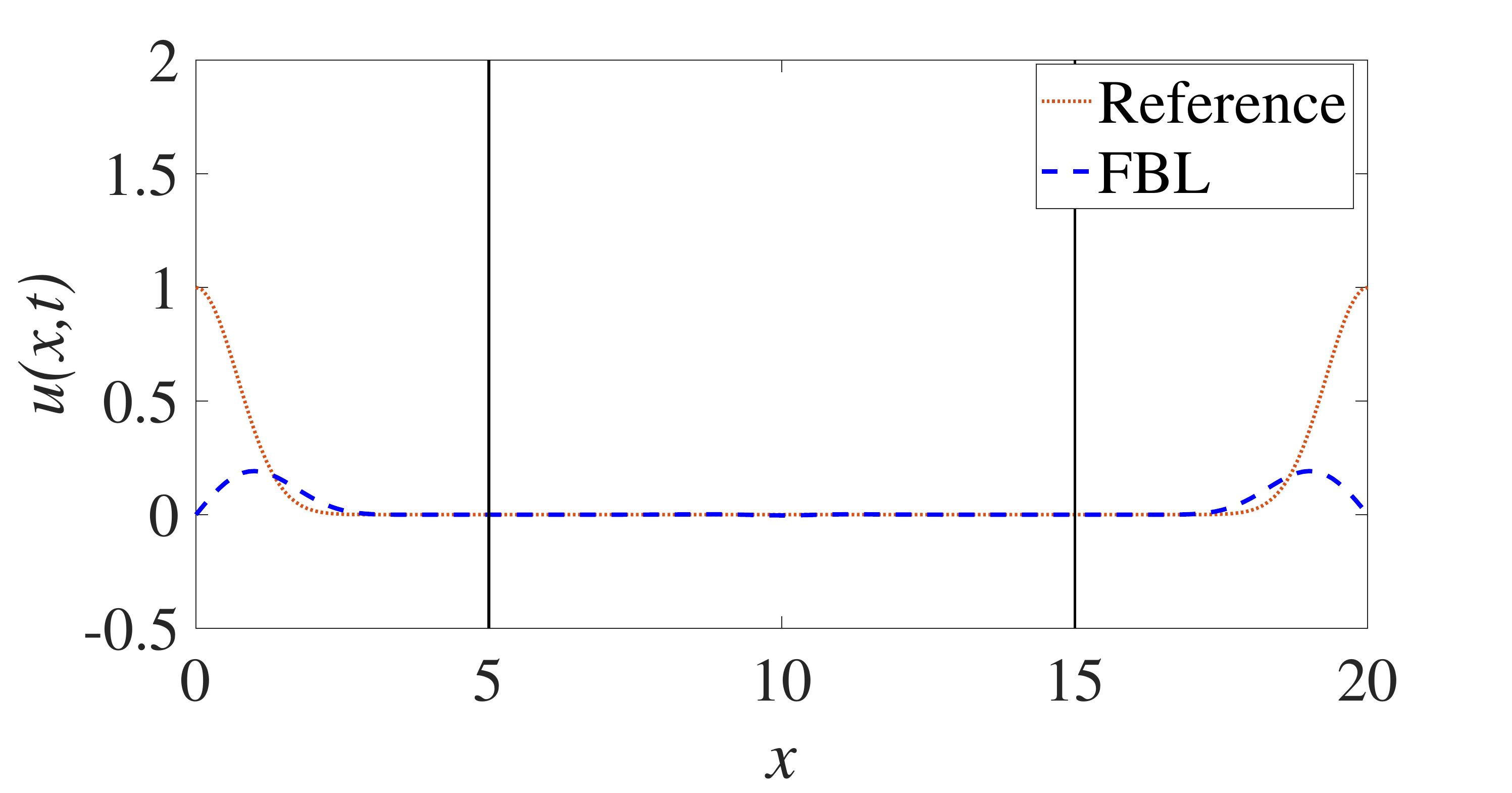}
\includegraphics[clip, trim=0.5cm 0cm 1.7cm 0cm, width=0.47\linewidth]{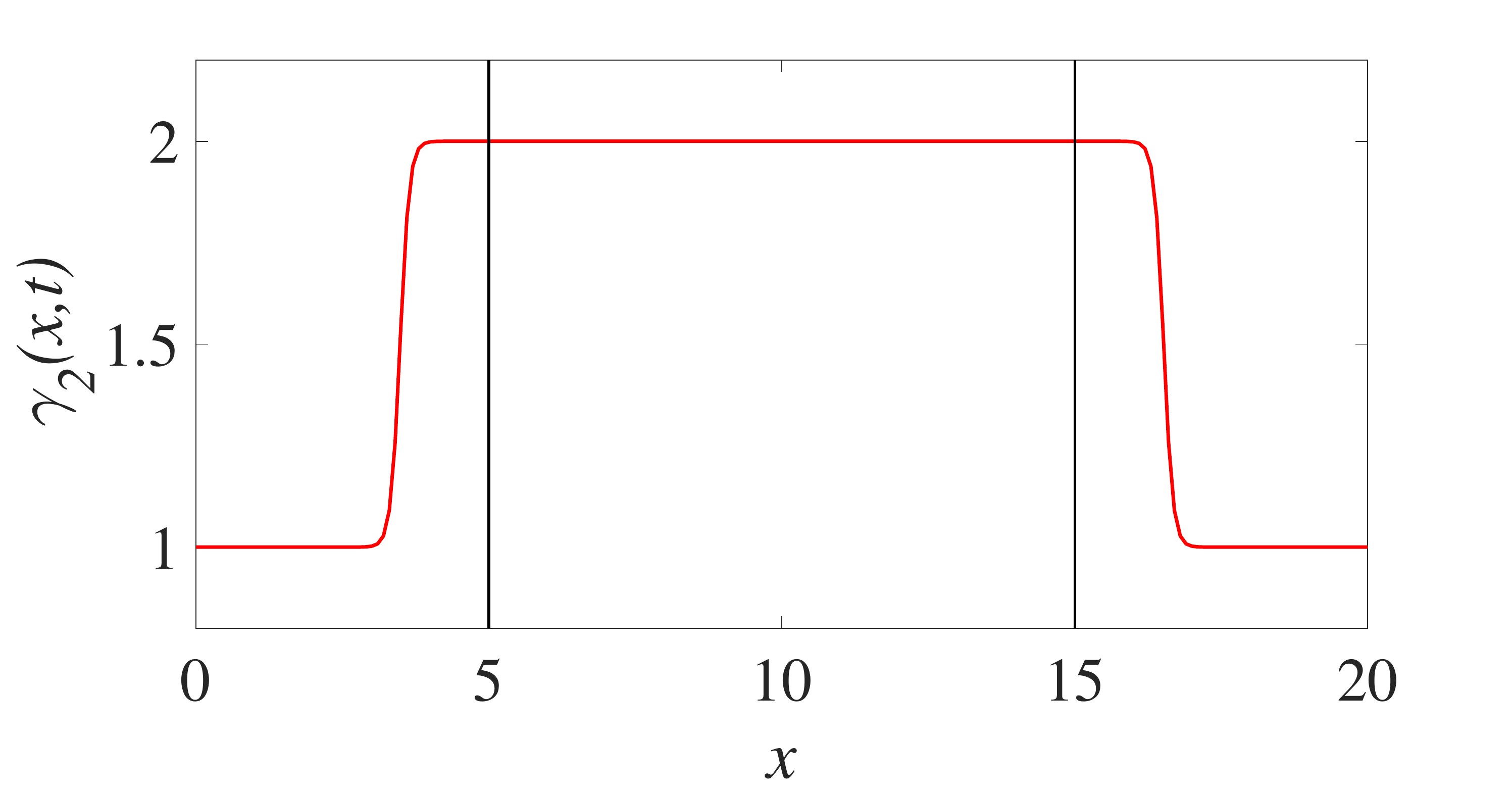}
\\[-8 pt]
$t=9$\\[-13 pt]
\noindent\rule{4cm}{0.7pt}
\\
\includegraphics[clip, trim=0.5cm 0cm 1.7cm 0cm, width=0.47\linewidth]{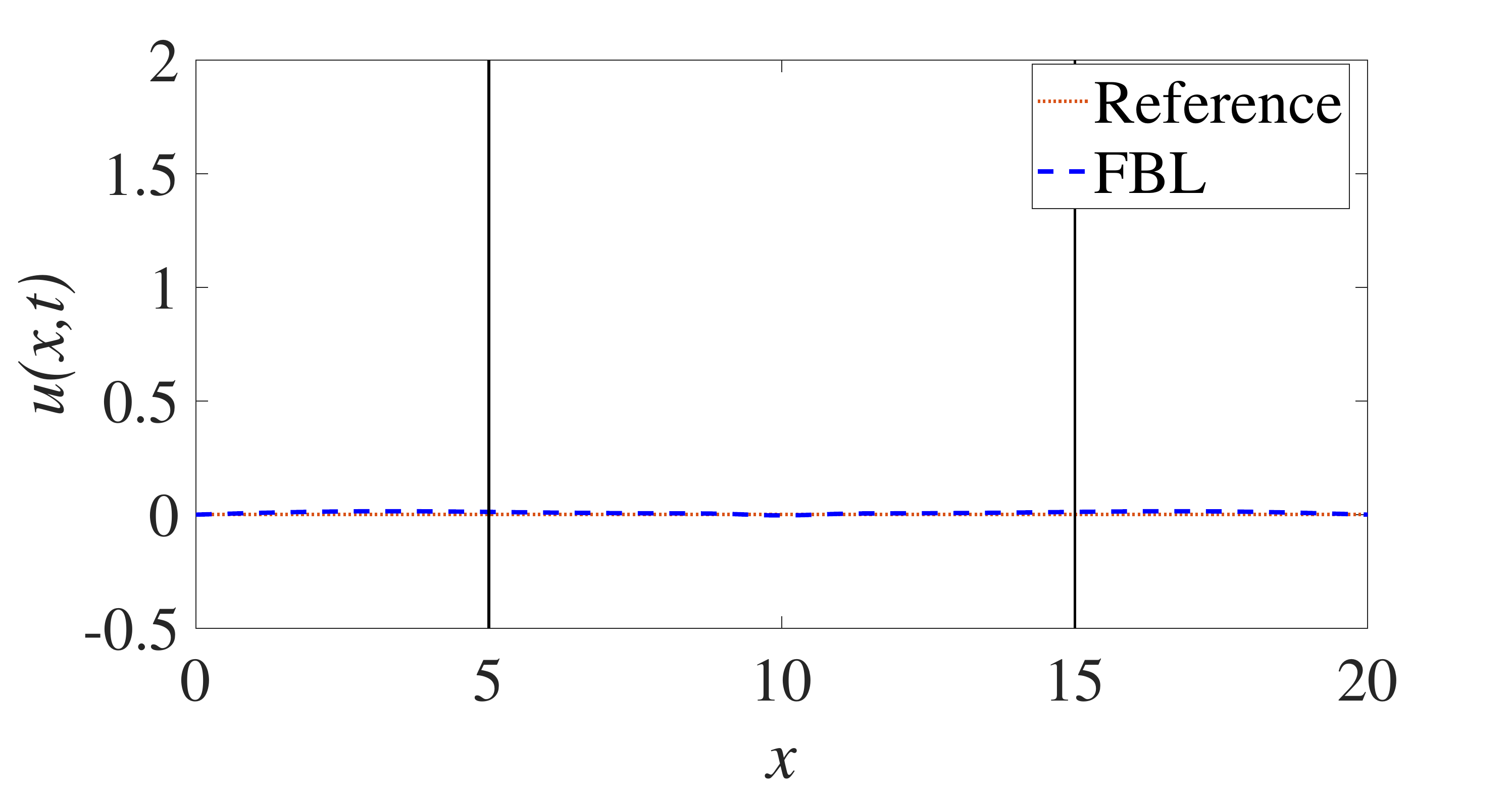}
\includegraphics[clip, trim=0.5cm 0cm 1.7cm 0cm, width=0.47\linewidth]{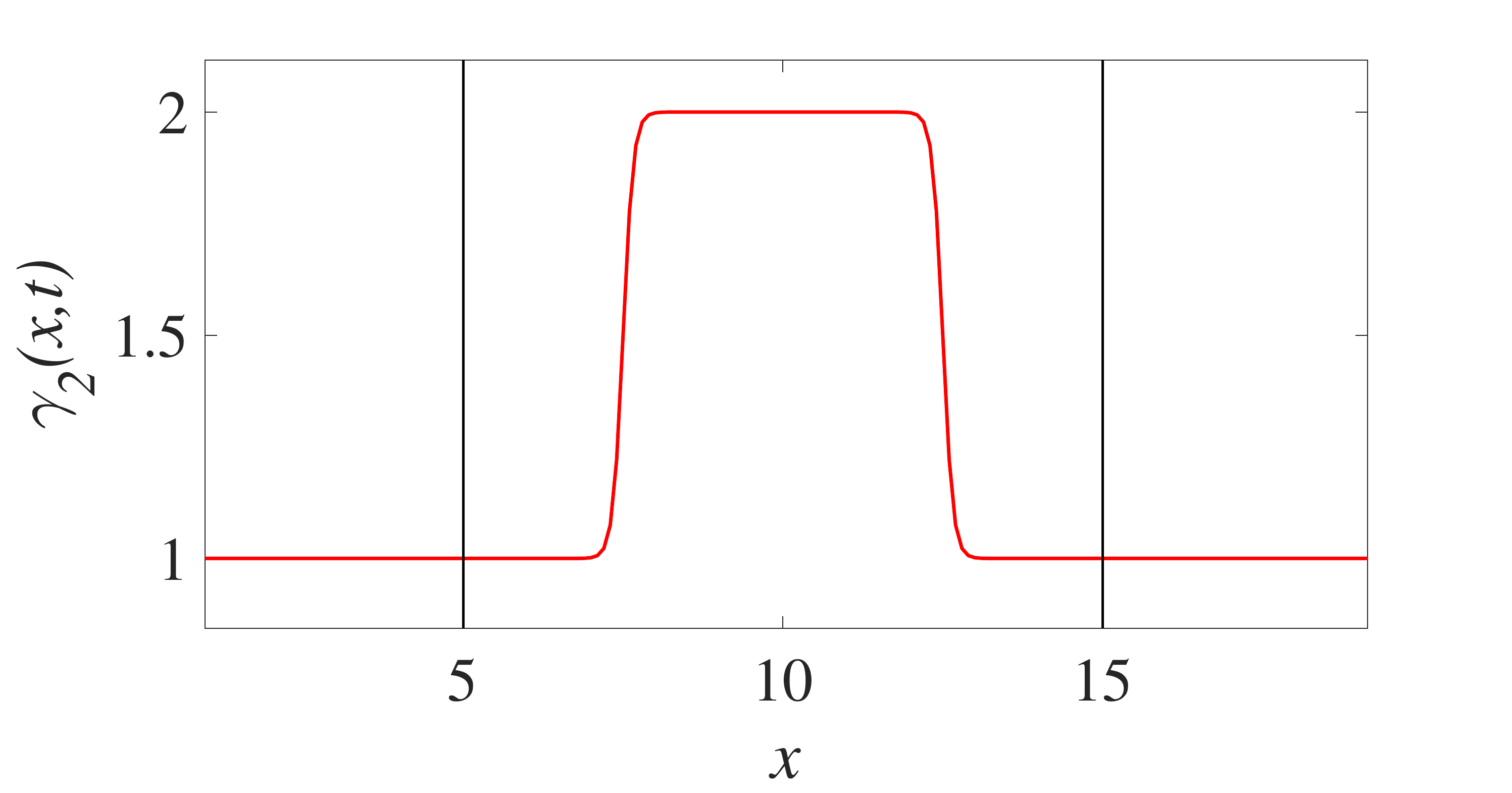}
\vspace{-0.2 in}
%%%%%%%%%%
\caption{Time fractional one-dimensional two-way wave (Example \ref{Ex: FBL two way 1D Xuan}):  time-dependent variable-order function $\gamma_2(x,t)$, given by \eqref{eq:TimeDependentOrd}. Left column: FBL and reference solution at $t=0, 5, 10, 15$. Right column: variable-order at the corresponding times.
The vertical black solid lines indicate the boundaries of interior domain. The wave propagates to both sides and then is absorbed in the buffer layers next to the boundaries.
}\label{Fig:FracLayer4CaputoWave_SGL_Xuan}
\end{figure}

We first consider the case of time-independent variable-order function, given by the following equation,
\begin{equation}\label{eq:TimeIndependentOrd}
  \gamma_{1}(x)=
  \left\{
  \begin{array}{ll}
  1, &x\in[x_L-\delta,x_L-\frac{\delta}{2})
  \\[3pt]
  1.5+(0.5-\epsilon)\tanh(\omega(x-x_L+\frac{\delta}{2})), &x\in[x_L-\frac{\delta}{2},x_L),
  \\[3pt]
  2, &x\in[x_L,x_R],
  \\[3pt]
  1.5+(0.5-\epsilon)\tanh(\omega(x_R+\frac{\delta}{2}-x)), &x\in(x_R,x_R+\frac{\delta}{2}),
  \\[3pt]
  1, &x\in(x_R+\frac{\delta}{2},x_R+\delta], 
  \end{array}
  \right.
\end{equation}
with the parameters $(\delta,\omega)=(5, 20)$; The plot of $\gamma_{1}(x)$ is shown in Fig. \ref{Fig:TimeIndependentOrd}. It is evident that in this setting, we recover the wave equation \eqref{eq:FracWaveCaputo} in the interior domain and then gradually switch to the diffusion equation in the buffer layers. The numerical results, however, show obvious reflections of the wave to the interior domain when the wave hits the boundary; see Fig. \ref{Fig:FracLayer4CaputoWave_SGL_1D}. This is because \eqref{eq:FracWaveCaputo} is local in space and remains as a wave equation in the interior domain. Therefore, even a small reflection from the boundary would propagate in the interior domain. A remedy to this issue is to consider a suitable time-dependent variable-order function \cite{Zhao@JCP2015}, such that when the wave passes through the boundaries, the equation switches from wave to diffusion equation in the interior domain to absorb the reflections. The following time-dependent variable function was proposed in \cite{Zhao@JCP2015} as 
\begin{equation}\label{eq:TimeDependentOrd}
  \begin{aligned}
      \gamma_{2}(x,t){}
     =&\frac{1}{2}\left\{\tanh\left[2\pi\left(x-\frac{t}{4}\tanh(\pi t-31)-\frac{t}{4}\right)\right]\right.
     \\[3pt]
     &\left.+\tanh\left[2\pi\left(20-\frac{t}{4}\tanh(\pi t-31)-\frac{t}{4}-x\right)\right]\right\}+1.
  \end{aligned}
\end{equation}
We see from the numerical results shown in Fig. \ref{Fig:FracLayer4CaputoWave_SGL_Xuan} (right column) how $\gamma_{2}(x,t)$ is being adjusted in the interior domain to absorb the reflections from the boundary. At $t=9$, even when the wave has completely left the interior domain, the interior sub-domains next to the boundaries are still diffusion-dominated ($\gamma_2 = 1$) to damp out any reflections. Although effective, this approach requires precise knowledge of the time that the wave reaches the boundaries to design the time-dependent variable-order function $\gamma_2(x,t)$. This may not be feasible in practical problems as such information may not be available a priori. What is often known in most of problems, however, are the boundaries of the domain of interest. This motivates us to redirect our focus on constructing an absorbing boundary layer based on time-independent but space-fractional derivatives.

%%%%%%%%%%%%%%%%%%%%%%%%%%%%%%%%%%%%%%%%%%%%%%%%%%%%%%%%%
 
%%%%%%%%%%%%%%%%%%%%%%%%%%%%%%%%%%%%%%%%%%%%%%%%%%%%%%%%%%

%%%%%%%%%%%%%%%
\subsection{Formulation of FBL}
%%%%%%%%%%%%%%%

In order to construct FBLs for the one-dimensional two-way wave, we first split the two-way wave into a system of one-way waves and then apply FBLs as was discussed in section \ref{Sec:Advc}.

The two-way wave equation \eqref{eq:WaveEq} can be split into a system of one-way wave equations by the mappings $v(x,t)=\frac{\partial}{\partial t}u(x,t)$ and $w(x,t)=c\,\frac{\partial}{\partial x}u(x,t)$. Thus, \eqref{eq:WaveEq} is equivalent to the following coupled system
\begin{equation}
  \left\{
     \begin{aligned}
     &\frac{\partial}{\partial t}v(x,t)=c\,\frac{\partial}{\partial x}w(x,t),
     \\[3pt]
     &\frac{\partial}{\partial t}w(x,t)=c\,\frac{\partial}{\partial x}v(x,t), 
     \end{aligned}
  \right.
\end{equation} 
with the initial conditions $v(x,0)=\varphi(x)$ and $w(x,0)=c\,\frac{\partial}{\partial x}u_{0}(x)$. In the matrix form, we can rewrite the equations as 
\begin{equation}
\label{eq: coupled split wave eq}
  \begin{aligned}
     \frac{\partial}{\partial t}
     \begin{bmatrix}
     v
     \\[3pt]
     w
     \end{bmatrix}
     =c\,
     \begin{bmatrix}
     0&1
     \\[3pt]
     1&0
     \end{bmatrix}
     \frac{\partial}{\partial x}
     \begin{bmatrix}
     w
     \\[3pt]
     v
     \end{bmatrix},
  \end{aligned}
\end{equation} 
which can be decoupled using the eigenvectors of the coefficient matrix and via the mappings $V=\frac{\sqrt{2}}{2}(w-v)$ and $W=\frac{\sqrt{2}}{2}(v+w)$. Therefore, we have
\begin{equation}
     \frac{\partial}{\partial t}
     \begin{bmatrix}
     V
     \\[3pt]
     W
     \end{bmatrix}
     =c\,
     \begin{bmatrix}
     -1&0
     \\[3pt]
     0&1
     \end{bmatrix}
     \frac{\partial}{\partial x}
     \begin{bmatrix}
     V
     \\[3pt]
     W
     \end{bmatrix},
\end{equation} 
which yields the following decoupled system 
\begin{equation}
\label{eq:UncoupledAdvc}
  \left\{
     \begin{aligned}
     &\frac{\partial V}{\partial t}=-c\,\frac{\partial V}{\partial x},
     \\[3pt]
     &\frac{\partial W}{\partial t}=c\,\frac{\partial W}{\partial x},
     \end{aligned}
  \right.
\end{equation}
subject to the initial conditions $V(x,0)=\frac{\sqrt{2}}{2}\left(c\,\frac{\partial}{\partial x}u_{0}(x)-\varphi(x)\right)$ and $W(x,0)=\frac{\sqrt{2}}{2}\left(\varphi(x)+c\,\frac{\partial}{\partial x}u_{0}(x)\right)$.

Instead of directly working with two-way wave equation \eqref{eq:WaveEq}, we develop the FBL for the decoupled split equation \eqref{eq:UncoupledAdvc}. This system is comprised of two uncoupled one-way wave equations with different propagation directions. Therefore, we can readily apply our method in section \ref{Sec:Advc} to construct the FBL for each one-way wave equation and then reconstruct the solution for the two-way wave equation. We first extend the domain to $[x_L-\delta,x_R+\delta]$ and let $(x,t)\in(x_L-\delta,x_R+\delta)\times(0,T]$. Then, similar to equations \eqref{eq:FracAdvcLeftRL} and  \eqref{eq:FracAdvcRightRL}, we transform \eqref{eq:UncoupledAdvc} into its variable-order fractional counterpart as 
\begin{equation}
\label{eq:FracLayer4Wave}
  \left\{
     \begin{aligned}
     &\frac{\partial}{\partial t}V(x,t)=c\,\prescript{RL}{x}{\mathcal{D}}_{x_R+\delta}^{\alpha(x)}V(x,t),
    %  &(x,t)\in(x_L-\delta,x_R+\delta)\times(0,T],
     \\[3pt]
     &\frac{\partial}{\partial t}W(x,t)=c\,\prescript{RL}{x_L-\delta}{\mathcal{D}}_{x}^{\beta(x)}W(x,t). 
     \end{aligned}
  \right.
\end{equation} 
Here, the initial conditions are  $V(x,0)=\frac{\sqrt{2}}{2}\left(c\,\frac{\partial}{\partial x}u_{0}(x)-\varphi(x)\right)$ and $W(x,0)=\frac{\sqrt{2}}{2}\left(\varphi(x)+c\,\frac{\partial}{\partial x}u_{0}(x)\right)$. We also impose the homogeneous Dirichlet boundary conditions to \eqref{eq:FracLayer4Wave}. The variable-order functions $\alpha(x)$ and $\beta(x)$ should satisfy the condition %
\begin{equation}
\label{eq:Condition1}
\left\{
\begin{array}{lll}
1+\epsilon,
&x\in[x_L,x_R],
&: \text{\textit{advection}}
\\[3pt]
(1,2],
&x\in[x_L-\bar{\delta},x_L)
\cup
(x_R,x_R+\bar{\delta}],
&: \text{\textit{penetration}}
\\[3pt]
2,
&x\in[x_L-\delta,x_L-\bar{\delta})
\cup
(x_R+\bar{\delta},x_R+\delta],
&: \text{\textit{diffusion}}
\end{array}
\right.
\end{equation} 
with the positive small parameter $\epsilon$ and $\bar{\delta}\in(0,\delta)$ being the length of the penetration region. When $\epsilon \rightarrow 0^{+}$, the FBL in \eqref{eq:FracLayer4Wave} reduces to \eqref{eq:UncoupledAdvc} in the interior domain $[x_L,x_R]$ and it becomes diffusion dominating in the buffer layers. The formulation of FBL for one-dimensional two-way wave equation \eqref{eq:WaveEq} in a bounded domain follows these steps:
\begin{enumerate}
  \item Consider a buffer layer of length $\delta$ and extend the domain to $[x_L-\delta,x_R+\delta]$.
  \item Consider proper variable-order functions $\alpha(x)$ and $\beta(x)$ and solve the uncoupled equation  \eqref{eq:FracLayer4Wave}. 
  \item Transform back the solution of \eqref{eq:FracLayer4Wave} into the coupled coordinate system by using the eigenvectors of coefficient matrix in \eqref{eq: coupled split wave eq}. 
  \item Reconstruct the solution by $u(x,t)=\int_{0}^{t}v(x,s){\rm d}s+u_{0}(x)$, or $u(x,t)=\frac{1}{c}\,\int_{x_L-\delta}^{x}w(s,t){\rm d}s$.
\end{enumerate}

\vspace{0.1 in} 

%%%%%%%%%%%%%%%
\subsection{Numerical simulations}
%%%%%%%%%%%%%%%
We investigate the performance of the FBL \eqref{eq:FracLayer4Wave} in absorbing the one-dimensional two-way wave given by \eqref{eq:WaveEq}. Without loss of generality, we consider $c=1$ and $\varphi(x,0)=0$. 

\begin{figure}[h]
\centering 
\includegraphics[width=0.8 \linewidth, clip, trim=1cm 1.8cm 0cm 1.5cm]{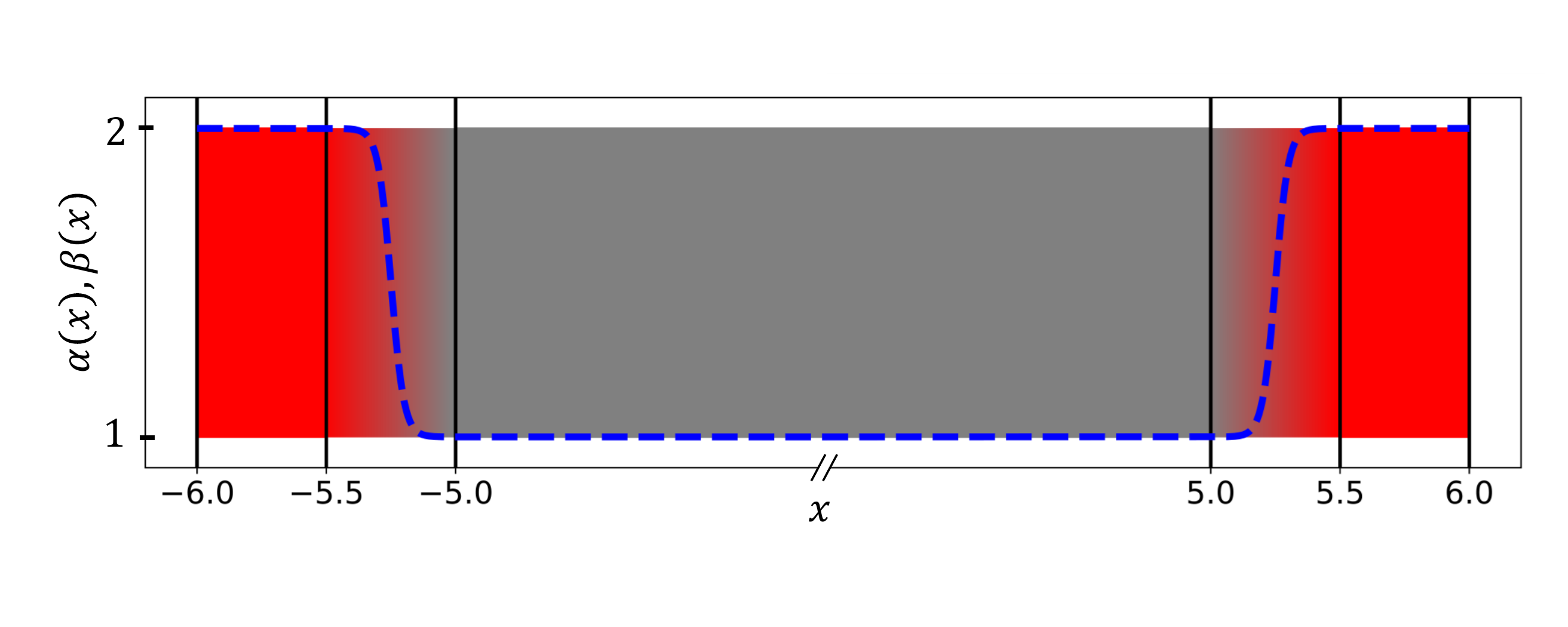}
\vspace{-0.1 in}
\caption{FBL for one-dimensional two-way wave (Example \ref{Ex: FBL two way 1D}). The variable-order functions are smooth $\tanh$ functions, given by Eq.  \eqref{eq:Vo4Wave_1D} with $\delta=1$, $\epsilon=10^{-5}$, and $\bar{\delta}=\frac{\delta}{2}$. The vertical black solid lines indicate the boundaries of interior domain and the penetration region. }\label{Fig:alf_1D_Wave}
\end{figure}

\begin{example}\label{Ex: FBL two way 1D} 
Let $(x_L,x_R)=[-5,5]$ be the interior domain. We attach two buffer layers of length $\delta=1$ next to the boundaries and extend the domain to $[-6,6]$. We assume smooth initial conditions $u_{0}(x)=e^{-x^{2}}$ and $\varphi(x)=0$. We consider the variable-order functions  $\alpha(x)$ and $\beta(x)$ both given by
\begin{equation}\label{eq:Vo4Wave_1D} 
   \left\{
   \begin{array}{ll}
   2, & x\in[x_L-\delta, x_L-\bar{\delta}),
   \\[3pt]
   1.5+(0.5-\epsilon)\tanh(20(-x+x_L-\frac{\bar{\delta}}{2})), &x\in[x_L-\bar{\delta},x_L),
   \\[3pt]
   1+\epsilon, &x\in[x_L,x_R],
   \\[3pt]
   1.5+(0.5-\epsilon)\tanh(20(x-x_R-\frac{\bar{\delta}}{2})), &x\in(x_R, x_R+\bar{\delta}],
   \\[3pt]
   2, & x\in(x_R+\bar{\delta}, x_R+\delta],
   \\[3pt]
      \end{array}
   \right.
\end{equation}
with $\epsilon=10^{-5}$ and $\bar{\delta}=\frac{\delta}{2}$.
The plots of $\alpha(x)$ and $\beta(x)$ in this case are given by Fig. \ref{Fig:alf_1D_Wave}. It is clear that $\alpha(x)$ and $\beta(x)$ satisfy condition \eqref{eq:Condition1} with the penetration regions of length $\bar{\delta}=0.5$.
\end{example}

\begin{figure}
\centering
$t=0$\\[-13 pt]
\noindent\rule{4cm}{0.7pt}
\\
\includegraphics[clip, trim=0.5cm 0cm 1.7cm 0cm, width=0.47\linewidth]{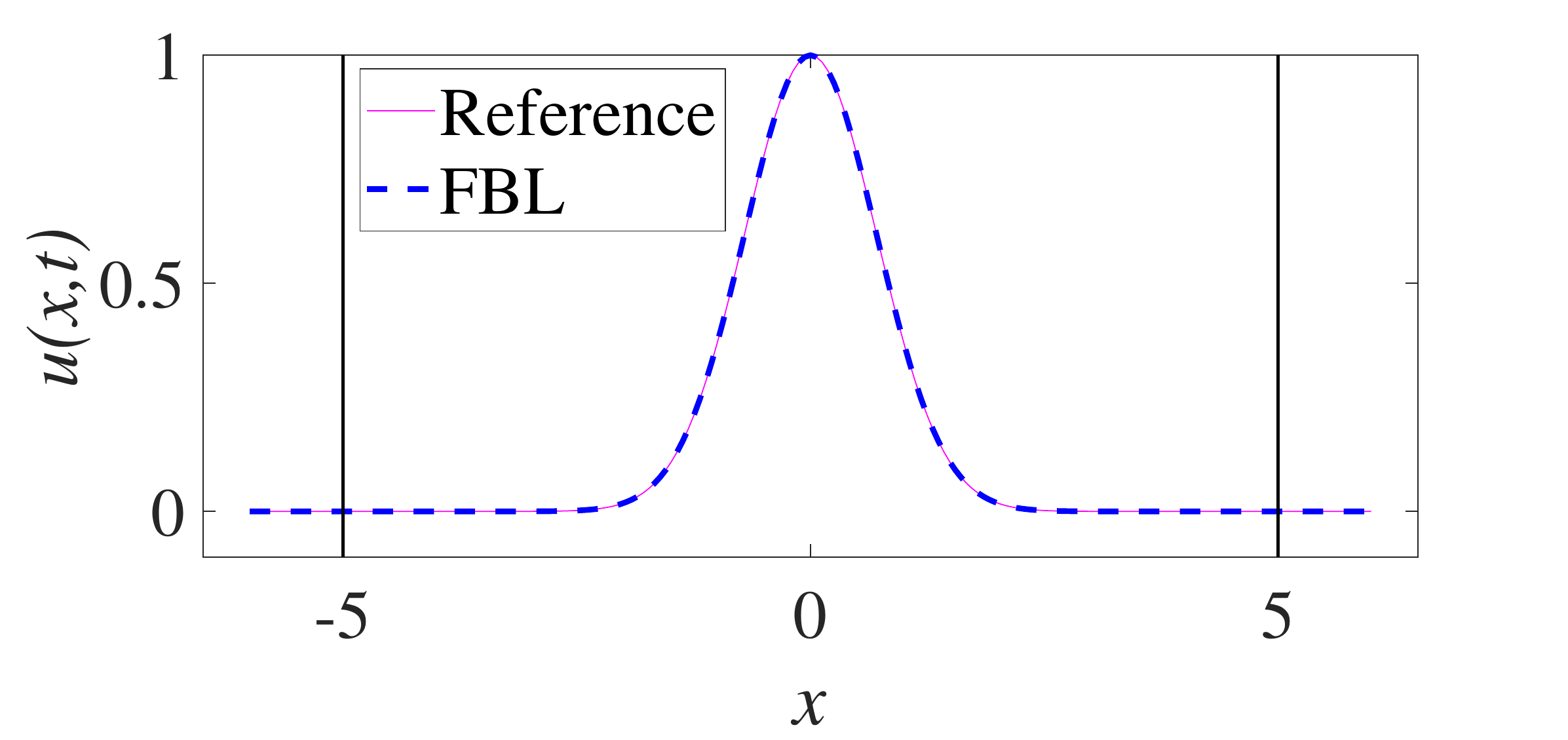} 
\\[-1 pt]
$t=3$\\[-13 pt]
\noindent\rule{4cm}{0.7pt}
\\
\includegraphics[clip, trim=0.5cm 0cm 1.7cm 0cm, width=0.47\linewidth]{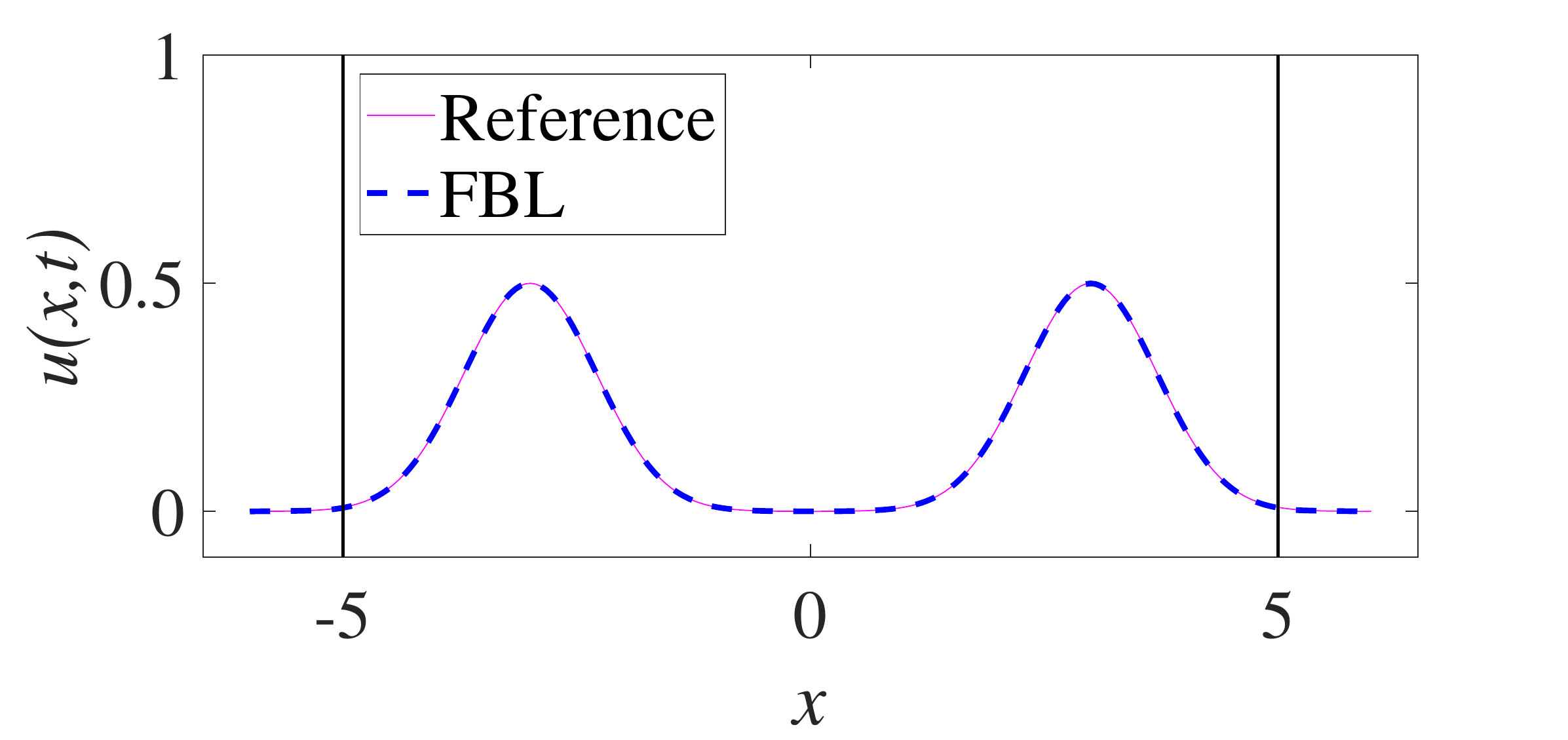}
\includegraphics[clip, trim=0.5cm 0cm 1.7cm 0cm, width=0.47\linewidth]{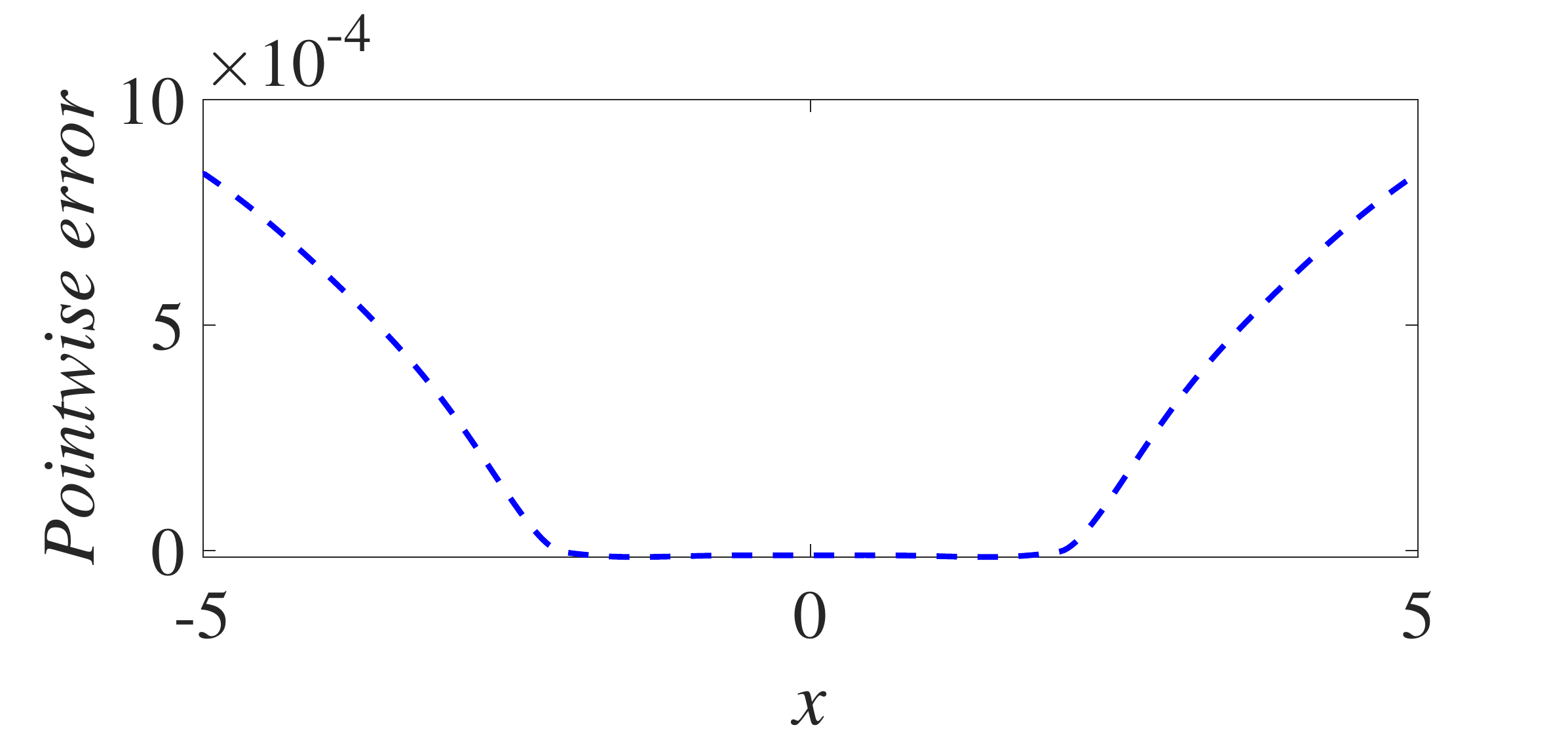}
\\[-8 pt]
$t=6$\\[-13 pt]
\noindent\rule{4cm}{0.7pt}
\\
\includegraphics[clip, trim=0.5cm 0cm 1.7cm 0cm, width=0.47\linewidth]{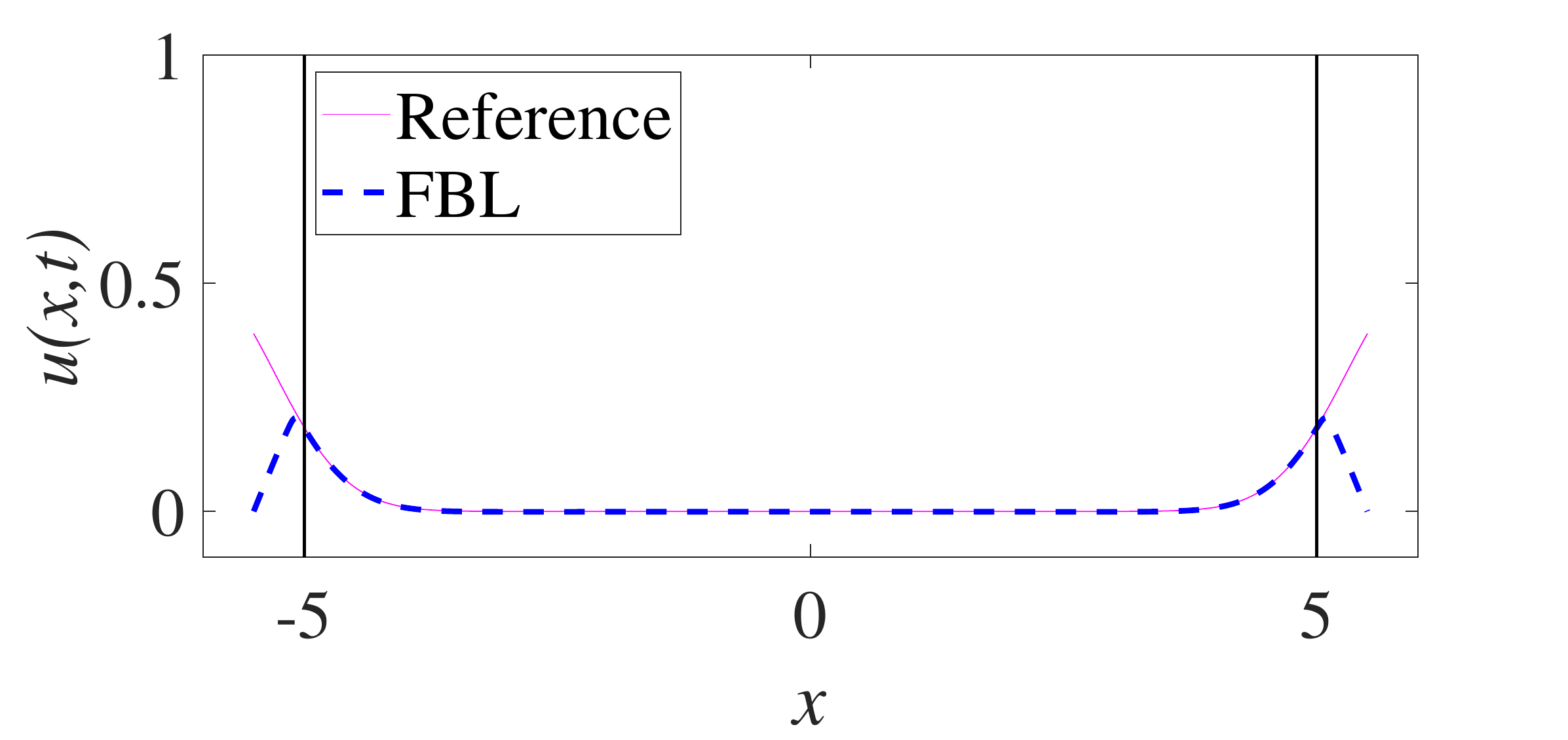}
\includegraphics[clip, trim=0.5cm 0cm 1.7cm 0cm, width=0.47\linewidth]{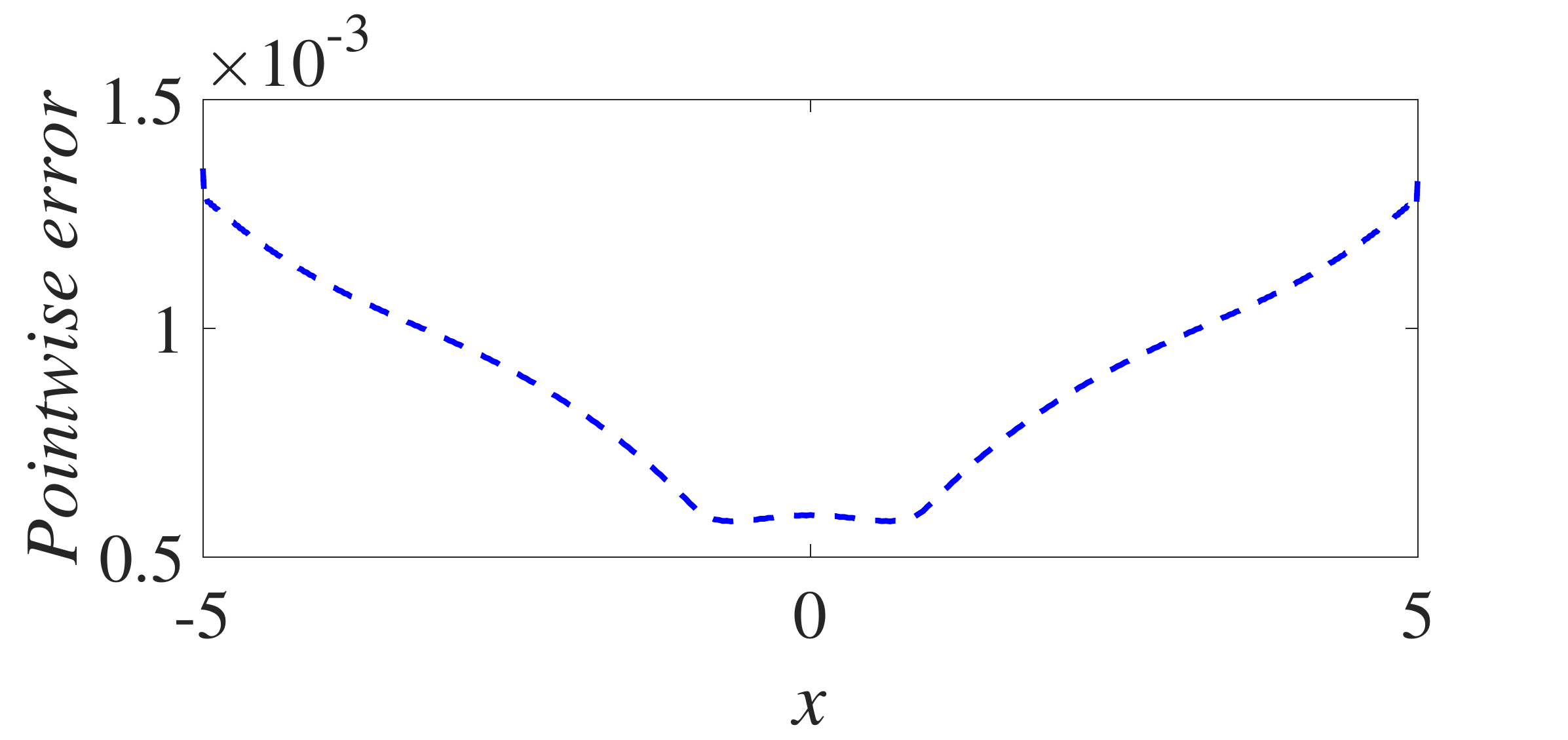}
\\[-8 pt]
$t=9$\\[-13 pt]
\noindent\rule{4cm}{0.7pt}
\\
\includegraphics[clip, trim=0.5cm 0cm 1.7cm 0cm, width=0.47\linewidth]{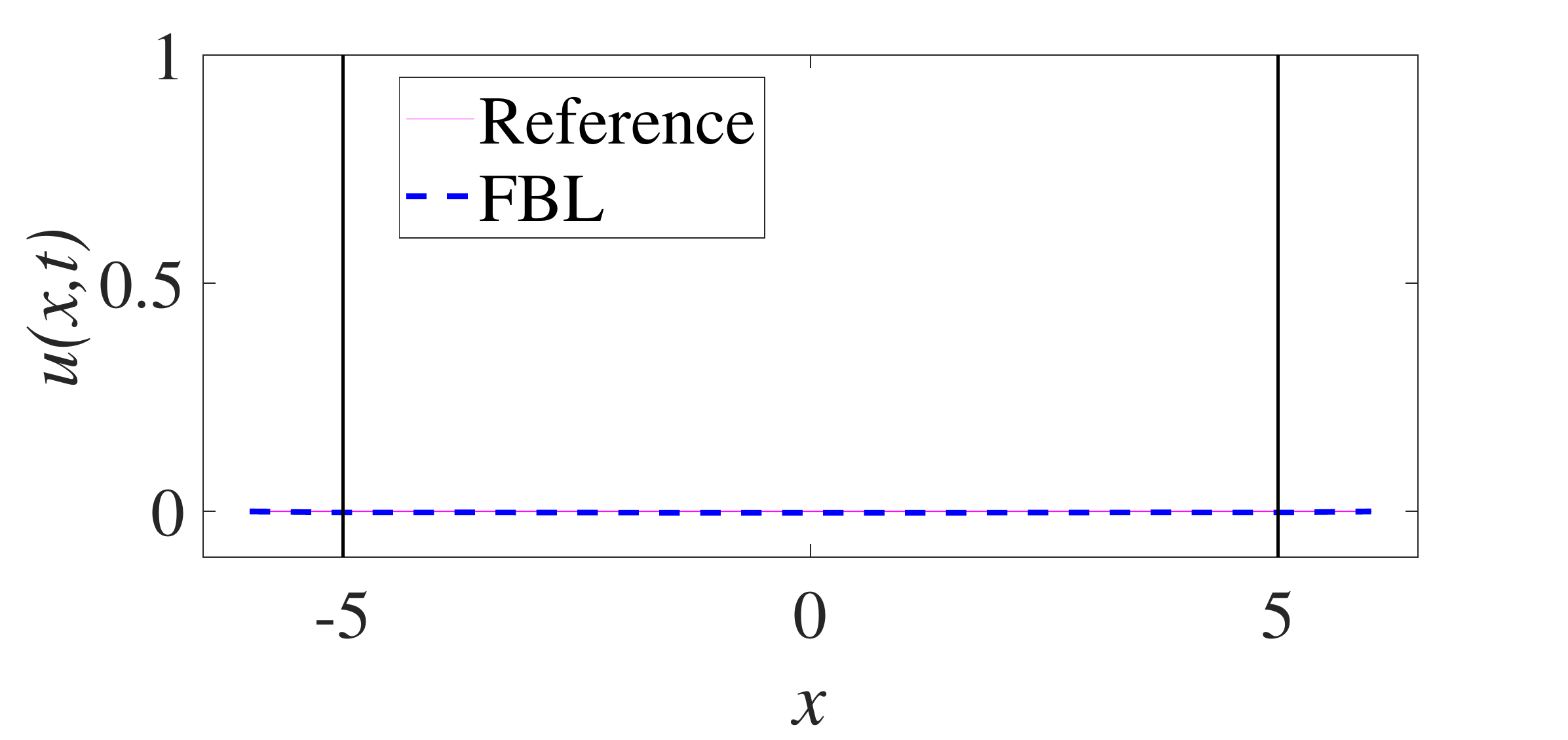}
\includegraphics[clip, trim=0.5cm 0cm 1.7cm 0cm, width=0.47\linewidth]{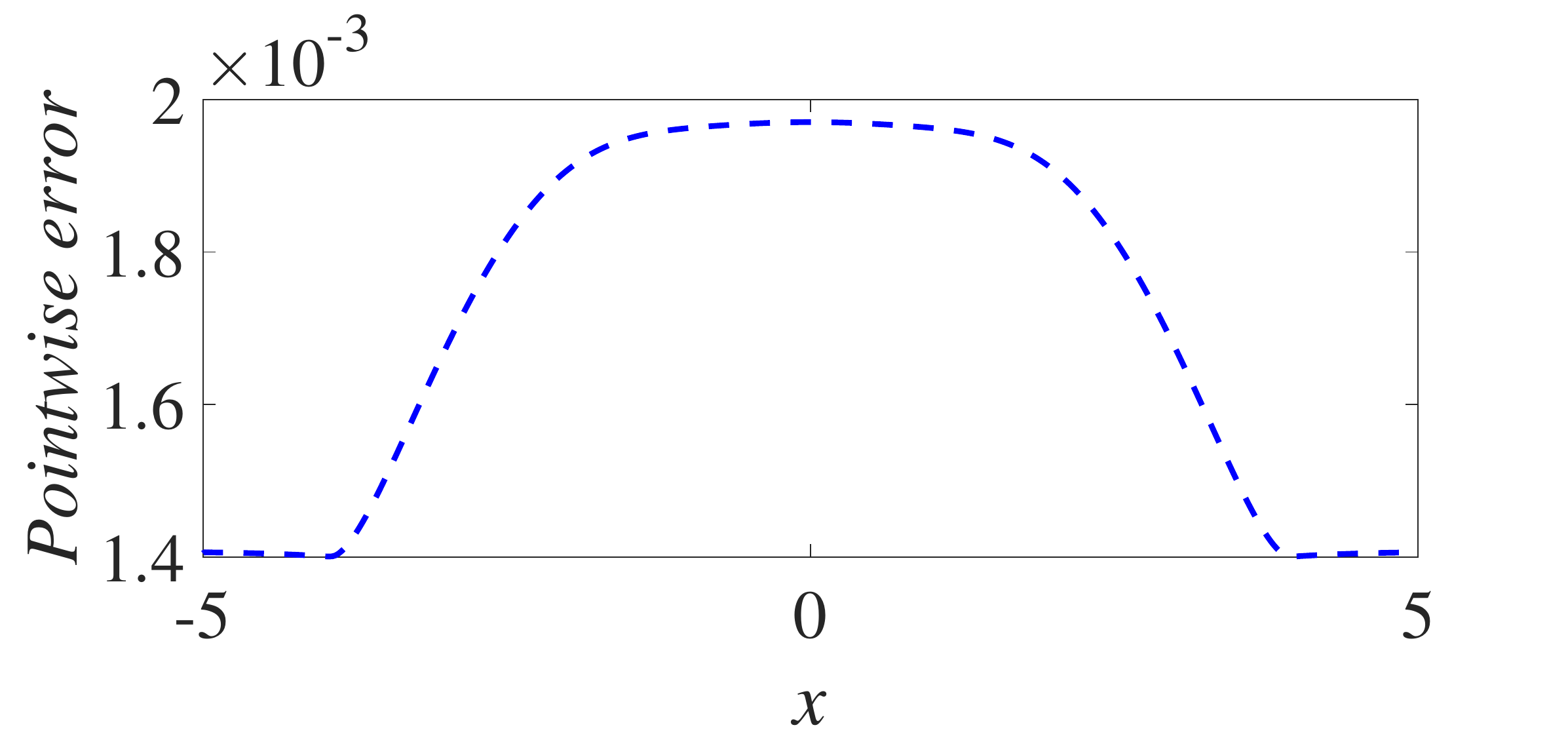}
\vspace{-0.2 in}
%%%%%%%%%%
\caption{One-dimensional two-way wave (Example \ref{Ex: FBL two way 1D}). Left column: FBL and reference solution at $t=0, 3, 6, 9$. Right column: pointwise error at the corresponding times.
The vertical black solid lines indicate the boundaries of interior domain. The wave propagates to both sides and then is absorbed in the buffer layers next to the boundaries.
}\label{Fig:NumSol_RefSol_PMLSol}
\end{figure}

In Example \ref{Ex: FBL two way 1D}, the numerical results are obtained from the fully discrete scheme \eqref{eq:FullyDiscrete_1D}, which is based on the spectral collocation method with $P=500$ collocation points in space and Crank-Nicolson method with the time step $\tau=10^{-3}$ for time integration. Note that the exact solution to the two-way wave equation \eqref{eq:WaveEq} in this case is  $u^{exact}(x,t)=\frac{u_{0}(x+t)+u_{0}(x-t)}{2}$. We regard it as the reference solution to compare our results in the interior domain. Figure \ref{Fig:NumSol_RefSol_PMLSol} shows a successful performance of the FBL \eqref{eq:FracLayer4Wave}. The left column shows different snapshots of the wave reaching to the left and right boundaries, penetrating into the buffer layers, and then being completely absorbed. The right column shows the pointwise error in the interior domain. Similar to the discussion on the one-way waves, here the pointwise error $|u^{exact}(x,t)-u(x,t)|$ consists of slight diffusion and numerical errors. 

%%%%%%%%%%%%%%%%%%%%%%%%%%%%%%%%%%%%%%%%%%%%%%%%

%%%%%%%%%%%%%%%%%%%%%%%%%%%%%%%%%%%%%%%%%%%%%%%%%%% 
 
%%%%%%%%%%%%%%%
\section{FBL: Two-dimensional two-way waves}\label{Sec:Wave_2D}
%%%%%%%%%%%%%%%
%\section{Fractional buffer layers for two-way waves in two dimensions}\label{Sec:Wave_2D}
We consider the two-dimensional two-way wave equation
\begin{equation}\label{eq:WaveEq_2D}
    \frac{\partial^{2}}{\partial t^{2}}u(x,y,t)
   =c^2\,\Delta u(x,y,t),
\end{equation}
with the initial conditions $u(x,y,0)=u_{0}(x,y)$ and $\frac{\partial}{\partial t}u(x,y,0)=\varphi(x,y)$. Here the Laplacian operator is $\Delta =\frac{\partial^{2}}{\partial x^{2}}+\frac{\partial^{2}}{\partial y^{2}}$ and the coefficient $c\geq 0$ is the propagation velocity. We split the two-way wave equation \eqref{eq:WaveEq_2D} into a coupled system of one-way equations via $v=\frac{\partial}{\partial t}u$, $w_{1}=c\,\frac{\partial}{\partial x}u$,
and $w_{2}=c\,\frac{\partial}{\partial y}u$. Thus,
\begin{equation}\label{eq:ScalarWaveEq_2D}
\left\{
    \begin{aligned}
       &\frac{\partial}{\partial t}v(x,y,t)
       =c\,\frac{\partial}{\partial x}w_{1}(x,y,t)
       +c\,\frac{\partial}{\partial y}w_{2}(x,y,t),
       \\[3pt]
       &\frac{\partial}{\partial t}w_{1}(x,y,t)
       =c\,\frac{\partial}{\partial x}v(x,y,t),
       \\[3pt]
       &\frac{\partial}{\partial t}w_{2}(x,y,t)
       =c\,\frac{\partial}{\partial y}v(x,y,t),
    \end{aligned}
\right.
\end{equation}
along with the initial conditions $v(x,y,0)=\varphi(x,y)$, $w_{1}(x,y,0)=c\,\frac{\partial}{\partial x}u_{0}(x,y)$, and $w_{2}(x,y,0)=c\,\frac{\partial}{\partial y}u_{0}(x,y)$. 

%%%%%%%%%%%%%%%
\subsection{Formulation of FBL}
%%%%%%%%%%%%%%%
We construct the FBLs that let the two-dimensional wave penetrate through the boundaries of the interior bounded domain $\Omega=[x_L,x_R]\times[y_L,y_R]$ and then get fully absorbed in the buffer layer of length $\delta>0$ attached to the boundaries. Therefore, we extend the domain to $\overline{\Omega}=[x_L-\delta,x_R+\delta]\times[y_L-\delta,y_R+\delta]$. Similar to the one-dimensional case, we develop the FBL by considering the split equation \eqref{eq:ScalarWaveEq_2D} and then reconstruct the solution for the wave equation. Since we cannot readily decouple \eqref{eq:ScalarWaveEq_2D}, the approach introduced in Sec. \ref{Sec:Advc} is not directly applicable here. Nevertheless, we introduce the fractional counterpart of the coupled system \eqref{eq:ScalarWaveEq_2D} by adding extra diffusion terms in the buffer layers.
 
Let $(x,y,t)\in\overline{\Omega}\times[0,T]$ and denote 
\begin{align*}
\prescript{RL}{-}{\mathcal{D}}_{x}^{\alpha(x)}
&=\frac{1}{2}\left[\prescript{RL}{x_L-\delta}{\mathcal{D}}_{x}^{\alpha(x)}-\prescript{RL}{x}{\mathcal{D}}_{x_R+\delta}^{\alpha(x)}\right],
\,\,\,
\prescript{RL}{+}{\mathcal{D}}_{x}^{\alpha(x)}
=\frac{1}{2}\left[\prescript{RL}{x_L-\delta}{\mathcal{D}}_{x}^{\alpha(x)}+\prescript{RL}{x}{\mathcal{D}}_{x_R+\delta}^{\alpha(x)}\right] ,
\\
\prescript{RL}{-}{\mathcal{D}}_{y}^{\beta(y)}
&=\frac{1}{2}\left[\prescript{RL}{y_L-\delta}{\mathcal{D}}_{y}^{\beta(y)}-\prescript{RL}{y}{\mathcal{D}}_{y_R+\delta}^{\beta(y)}\right] , 
\,\,\,
\prescript{RL}{+}{\mathcal{D}}_{y}^{\beta(y)}
=\frac{1}{2}\left[\prescript{RL}{y_L-\delta}{\mathcal{D}}_{y}^{\beta(y)}+\prescript{RL}{y}{\mathcal{D}}_{y_R+\delta}^{\beta(y)}\right]. 
\end{align*} 
Then, the FBL formulation for the coupled system \eqref{eq:ScalarWaveEq_2D} will take the form
\begin{equation}\label{eq:FracLayer4ScalarWaveEq_2D}
\left\{
    \begin{aligned}
        \frac{\partial}{\partial t}v(x,y,t)
       &=c\,\prescript{RL}{-}{\mathcal{D}}_{x}^{\alpha(x)}w_{1}(x,y,t) 
       +c\,\prescript{RL}{+}{\mathcal{D}}_{x}^{\alpha(x)}v(x,y,t)
       \\[3pt] 
       &+c\,\prescript{RL}{-}{\mathcal{D}}_{y}^{\beta(y)}w_{2}(x,y,t) 
       +c\,\prescript{RL}{+}{\mathcal{D}}_{y}^{\beta(y)}v(x,y,t), 
       \\[3pt]
        \frac{\partial}{\partial t}w_{1}(x,y,t)
       &=c\,\prescript{RL}{+}{\mathcal{D}}_{x}^{\alpha(x)}w_{1}(x,y,t)
       +c\,\prescript{RL}{-}{\mathcal{D}}_{x}^{\alpha(x)}v(x,y,t), 
       \\[3pt]
        \frac{\partial}{\partial t}w_{2}(x,y,t)
       &=c\,\prescript{RL}{+}{\mathcal{D}}_{y}^{\beta(y)}w_{2}(x,y,t)
       +c\,\prescript{RL}{-}{\mathcal{D}}_{y}^{\beta(y)}v(x,y,t), 
    \end{aligned}
\right.
\end{equation}
with the same initial conditions as in \eqref{eq:ScalarWaveEq_2D} and homogeneous Dirichlet boundary conditions. Here, the variable-order functions satisfy
\begin{equation}\label{eq:Condition3}
   \left\{
   \begin{array}{lll}
   \alpha(x)\rightarrow1^{+},
   &x\in[x_L,x_R],
   &: \text{\textit{advection}}
   \\[3pt]
   \alpha(x)\in(1,2],
   &x\in[x_L-\bar{\delta},x_L)\cup(x_R,x_R+\bar{\delta}],
   &: \text{\textit{penetration}}
   \\[3pt]
   \alpha(x)=2,
   &x\in[x_L-\delta,x_L-\bar{\delta})\cup(x_R+\bar{\delta},x_R+\delta],
   &: \text{\textit{diffusion}}
   \end{array}
   \right.
\end{equation}
and
\begin{equation}\label{eq:Condition4}
   \left\{
   \begin{array}{lll}
   \beta(y)\rightarrow1^{+},
   &y\in[y_L,y_R],
   &: \text{\textit{advection}}
   \\[3pt]
   \beta(y)\in(1,2],
   &y\in[y_L-\bar{\delta},y_L)\cup(y_R,y_R+\bar{\delta}],
   &: \text{\textit{penetration}}
   \\[3pt]
   \beta(y)=2,
   &y\in[y_L-\delta,y_L-\bar{\delta})\cup(y_R+\bar{\delta},y_R+\delta],
   &: \text{\textit{diffusion}}
   \end{array}
   \right.
\end{equation}
with $\bar{\delta}\in(0,\delta)$ being the length of the penetration region.

\begin{remark}[Anisotropic Propagation]
Different definitions of variable-order functions $\alpha(x)$ and $\beta(y)$ provide the flexibility to formulate FBLs for anisotropic wave propagation. In each direction $x$ and $y$, we can consider different characteristics of FBLs with different sizes $\delta$ and penetration regions $\bar{\delta}$ as shown for example in Fig. \ref{Fig:alf}.
\end{remark}

It follows from the consistency relation \eqref{Eq: consistency} that \eqref{eq:FracLayer4ScalarWaveEq_2D} recovers the split equation \eqref{eq:ScalarWaveEq_2D} in the interior domain when $\alpha(x), \beta(y) \rightarrow 1^{+} $ and becomes diffusion dominant in the buffer layers when $\alpha(x), \beta(y) = 2 $. By comparing the integer-order equation \eqref{eq:ScalarWaveEq_2D} with the fractional one in \eqref{eq:FracLayer4ScalarWaveEq_2D}, we can consider the four additional terms  $c\,\prescript{RL}{+}{\mathcal{D}}_{x}^{\alpha(x)}v$, $c\,\prescript{RL}{+}{\mathcal{D}}_{y}^{\beta(y)}v$,
$c\,\prescript{RL}{+}{\mathcal{D}}_{x}^{\alpha(x)}w_{1}$, and $c\,\prescript{RL}{+}{\mathcal{D}}_{y}^{\beta(y)}w_{2}$ as damping terms in the buffer layers. In summary, the formulation of FBL for two-dimensional two-way wave equation \eqref{eq:WaveEq_2D} in a bounded domain follows these steps:
\begin{enumerate}
  \item Consider a buffer layer of length $\delta_x$ in $x$ direction and $\delta_y$ in $y$ direction and extend the computational domain to $[x_L-\delta_x,x_R+\delta_x]\times[y_L-\delta_y,y_R+\delta_y]$.
  \item Consider proper variable-order functions $\alpha(x)$ and $\beta(y)$ and solve the coupled  equation \eqref{eq:FracLayer4ScalarWaveEq_2D}.
  \item Reconstruct the solution by  $u(x,y,t)=\int_{0}^{t}v(x,y,s){\rm d}s+u_{0}(x,y)$, or $u(x,y,t)=\frac{1}{c}\,\int_{x_{L}-\delta}^{x}w_{1}(s,y,t){\rm d}s$, or $u(x,y,t)=\frac{1}{c}\,\int_{y_{L}-\delta}^{y}w_{2}(x,s,t){\rm d}s$.
  %This can be obtained  via quadrature rules.
\end{enumerate}

%%%%%%%%%%%%%%%
\subsection{Other possible fractional approaches}
%%%%%%%%%%%%%%%

One alternative fractional approach to construct the absorbing boundary layer for the two-dimensional wave equation \eqref{eq:WaveEq_2D} is to consider the variable-order time fractional wave equation as in Example \ref{Ex: FBL two way 1D Xuan}. In this case, we let the extended domain be $\overline{\Omega}=[x_L-\delta,x_R+\delta]\times[y_L-\delta,y_R+\delta]$ and consider
\begin{equation}\label{eq:FracWaveCaputo_2D}
   {}_{C}{\rm D}_{0,t}^{\gamma(x,y,t)}u(x,y,t)
   =c\,\Delta u(x,y,t),
\end{equation}
subject to the initial conditions $u(x,y,0)=u_{0}(x,y)$ and $\frac{\partial}{\partial t}u(x,y,0)=\varphi(x,y)$, and homogeneous Dirichlet boundary conditions. As discussed before, this formulation requires a prior knowledge of the time that the wave reaches the boundaries to precisely define the time-dependent variable-order function; this information may not be available in practice. 

Another instinctive approach to design the absorbing boundary layer for the two-dimensional wave propagation is by following the one-dimensional case and replacing the Laplacian operator $(-\Delta)$ by its fractional counterpart $(-\Delta)^{\frac{\gamma}{2}}$. We refer to \cite{Laplacian2020} for an extensive review of the fractional Laplacian and 
recall the following useful properties of these operators. 
\begin{equation}\label{eq:Semigroup4FracLaplace}
\left\{
  \begin{aligned}
   &\lim\limits_{\gamma\rightarrow 2^{-}}(-\Delta)^{\frac{\gamma}{2}}=-\Delta,
   \\[3pt]
   &   (-\Delta)^{\gamma_1}(-\Delta)^{\gamma_2}=(-\Delta)^{\gamma_1+\gamma_2},
   \ \gamma_1, \gamma_2\geq0.
  \end{aligned}
\right.
\end{equation}
Here, we adopt the following spectral definition of fractional Laplacian, given as
\begin{equation}
     (-\Delta)^{\frac{\gamma(x,y,t)}{2}}f(x,y)
    =\sum\limits_{j=1}^{\infty}\sum\limits_{k=1}^{\infty}
    c_{j,k}(\lambda_{j,k})^{\frac{\gamma(x,y,t)}{2}}\phi_{j,k}(x,y),
\end{equation}
where $\lambda_{j,k}=\frac{j^{2}\pi^{2}}{(x_R-x_L)^{2}}+\frac{k^{2}\pi^{2}}{(y_R-y_L)^{2}}$ and $\phi_{j,k}=\frac{2}{\sqrt{(x_R-x_L)(y_R-y_L)}}\sin\frac{j\pi x}{x_R-x_L}\sin\frac{k\pi y}{y_R-y_L}$ are eigenvalues and eigenfunctions of the Laplacian operator $(-\Delta)$ on the rectangular domain $\Omega=[x_L,x_R]\times[y_L,y_R]$. The coefficients are given by the inner product $c_{i,j}=(f, \phi_{i,j})$.

We first split the wave equation \eqref{eq:WaveEq_2D} into a decoupled system by denoting $v(x,y,t)=\frac{\partial}{\partial t}u(x,y,t)$ and $w(x,y,t)=i\,c\,(-\Delta)^{\frac{1}{2}}u(x,y,t)$ with $i = \sqrt{-1}$ and using the mappings $V=\frac{\sqrt{2}}{2}(w-v)$ and $W=\frac{\sqrt{2}}{2}(v+w)$. It then follows from \eqref{eq:WaveEq_2D} and \eqref{eq:Semigroup4FracLaplace} that 
\begin{equation}\label{eq:UncoupledShordinger}
  \left\{
     \begin{aligned}
     &\frac{\partial V}{\partial t}=-ic\,(-\Delta)^{\frac{1}{2}}V,
     \\[3pt]
     &\frac{\partial W}{\partial t}=ic\,(-\Delta)^{\frac{1}{2}}W, 
     \end{aligned}
  \right.
\end{equation}
subject to the initial conditions $V(x,y,0)=\frac{\sqrt{2}}{2}\left(ic\,(-\Delta)^\frac{1}{2}u_{0}(x,y)-\varphi(x,y)\right)$ and $W(x,y,0)=\frac{\sqrt{2}}{2}\left(ic\,(-\Delta)^\frac{1}{2}u_{0}(x,y)+\varphi(x,y)\right)$. The instinctive approach here is to follow the formulation in \eqref{eq:FracLayer4Wave} by extending the computational domain to $\overline{\Omega}=[x_L-\delta,x_R+\delta]\times[y_L-\delta,y_R+\delta]$ and considering the variable-order fractional differential equation of the form
\begin{equation}\label{eq:FracLayer4TwoDWaveComplexValuedVariableCoef}
  \left\{
     \begin{aligned}
     &\frac{\partial}{\partial t}V(x,y,t)
     =-ic\,(-\Delta)^{\frac{\gamma_{1}(x,y,t)}{2}}V(x,y,t), 
     \\[3pt]
     &\frac{\partial}{\partial t}W(x,y,t)
     =ic\,(-\Delta)^{\frac{\gamma_{2}(x,y,t)}{2}}W(x,y,t). 
     \end{aligned}
  \right.
\end{equation}
subject to the same initial conditions as \eqref{eq:UncoupledShordinger} and homogeneous Dirichlet boundary conditions. However, the decoupled system \eqref{eq:UncoupledShordinger} in fact consists  of two fractional  Shr\"{o}dinger's equations with order $\frac{1}{2}$, rather than the advection equations as in \eqref{eq:UncoupledAdvc}. Therefore, the solution to the decoupled system \eqref{eq:FracLayer4TwoDWaveComplexValuedVariableCoef} becomes oscillatory and does not decay for any arbitrary positive-valued $\gamma_{1}(x,y,t)$ and $\gamma_{2}(x,y,t)$ because all the eigenvalues of $-ic\,(-\Delta)^{\frac{\gamma_{1}(x,y,t)}{2}}$ and $ic\,(-\Delta)^{\frac{\gamma_{2}(x,y,t)}{2}}$ are imaginary numbers.
In this scenario, we need to modify the coefficients of \eqref{eq:FracLayer4TwoDWaveComplexValuedVariableCoef} in the following form of 
\begin{equation}
\label{eq:FracLayerLaplacian}
      \left\{
     \begin{aligned}
     &\frac{\partial}{\partial t}V(x,y,t)
     =-cp(x,y,t)(-\Delta)^{\frac{\gamma_{1}(x,y,t)}{2}}V(x,y,t),
     \\[3pt]
     &\frac{\partial}{\partial t}W(x,y,t)
     =cq(x,y,t)(-\Delta)^{\frac{\gamma_{2}(x,y,t)}{2}}W(x,y,t). 
     \end{aligned}
  \right.
\end{equation}
When the variable coefficients $p(x,y,t)=q(x,y,t)=i$ and the variable orders $\gamma_{1}(x,y,t)=\gamma_{2}(x,y,t)=1$, this coupled system recovers \eqref{eq:UncoupledShordinger}. However, when $p(x,y,t)=1$, $q(x,y,t)=-1$, and $\gamma_{1}(x,y,t)=\gamma_{2}(x,y,t)=2$, it reduces to a diffusive system. Hence, we can use \eqref{eq:FracLayerLaplacian} to formulate an absorbing boundary layer. However, we see that variable coefficients and orders should be defined as time-dependent functions. Thus, to fully absorb the waves passing through the boundaries, a prior knowledge of the time that the wave reaches the boundaries is required for accurate definition of $p, q, \gamma_{1}$, and $\gamma_{2}$. As discussed before, this may not be applicable in practical problems. Therefore, we only focus on the formulation of FBL based on \eqref{eq:FracLayer4ScalarWaveEq_2D} rather than the one given by \eqref{eq:FracLayerLaplacian}.
% with the properly defined coefficient functions $p(x,y)$ and $q(x,y)$ given by
% \begin{equation}
% \left\{
% \begin{array}{ll}
% p
% = q
% = i, &(x,y)\in[x_L,x_R]\times[y_L,y_R],
% \\[3pt]
% p
% =1, q
% = -1, &{\textit{elsewhere}}.
% \end{array}
% \right.
% \end{equation}
% Moreover, defining the appropriate time-dependent variable-order functions $\gamma_{1}(x,y,t)\in(1,2]$ and $\gamma_{2}(x,y,t)\in(1,2]$ that switch from 1 to 2 in suitable ways  requires a prior knowledge of the time that the wave reaches the boundaries, which is not applicable for practical problems. Therefore, we focus on the FBL based on \eqref{eq:FracLayer4ScalarWaveEq_2D} rather than the one given by \eqref{eq:FracLayerLaplacian}.

%%%%%%%%%%%%%%%
\subsection{Numerical simulations}
%%%%%%%%%%%%%%%
We investigate the performance of the proposed FBL \eqref{eq:FracLayer4ScalarWaveEq_2D} for the two-dimensional wave equation \eqref{eq:WaveEq_2D}. We compare FBL and PML to show the 
effectiveness of FBL in removing reflections in two-dimensional problems. To numerically solve the FBL and PML formulations, we employ a spectral collocation method in space and apply the second-order Adams-Bashforth method for the time integration. Here the numbers of the collocation points in $x$ and $y$ directions are  $(P_x,P_y)=(50,50)$, and the time step is $\tau=10^{-5}$. 

\begin{figure}[h!]
\centering 
\includegraphics[width=0.8\linewidth]{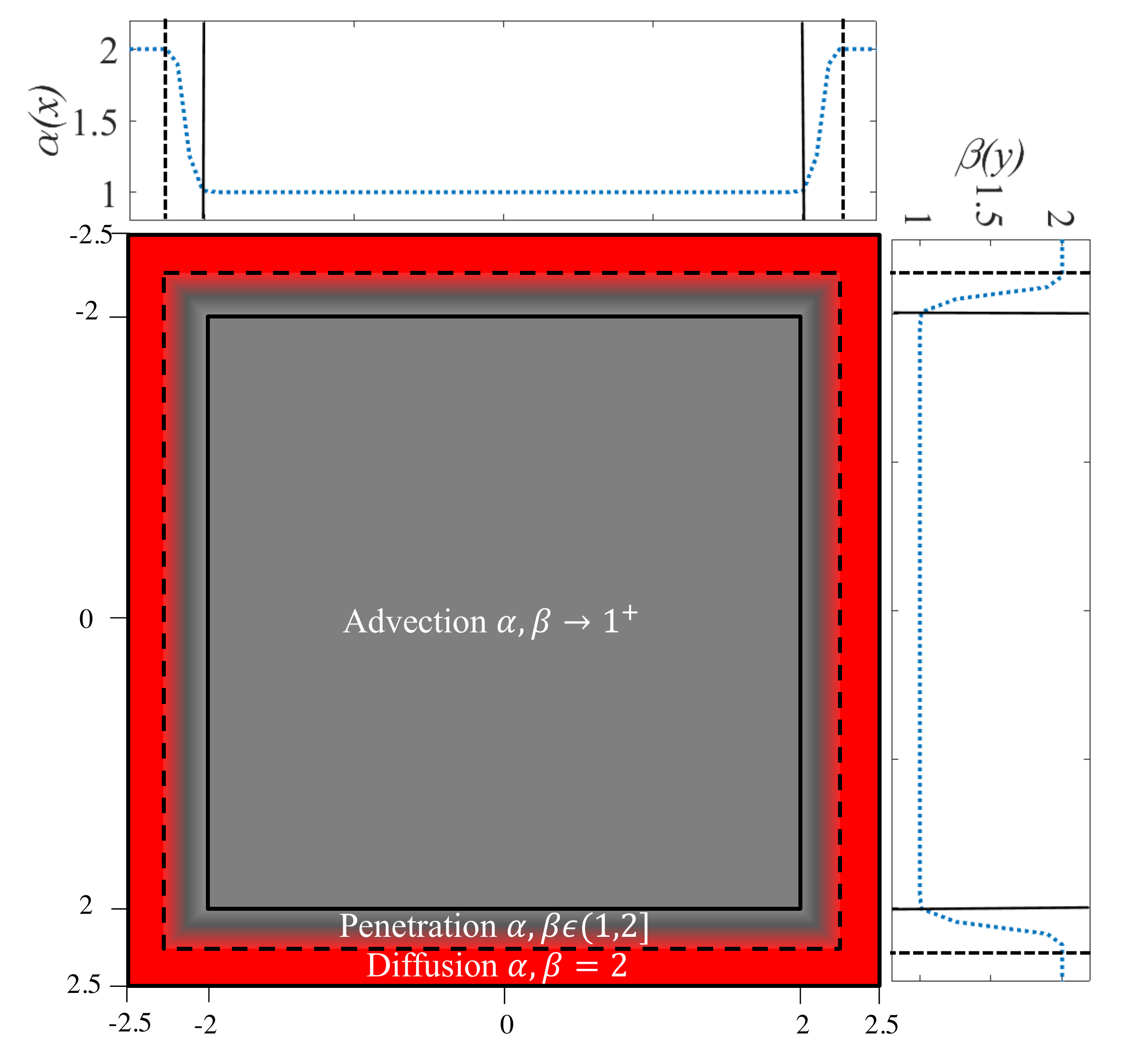}
\caption{FBL for two-dimensional two-way wave in Example. \ref{Ex: FBL two way 2D}.  The variable-order functions are given by \eqref{eq:Ord4Wave_2Dx} and \eqref{eq:Ord4Wave_2Dy} with $\delta=0.5$. The gray and red shaded areas are the advection and diffusion regions, respectively.}\label{Fig:VariableOrd}
\end{figure}

\begin{example}\label{Ex: FBL two way 2D}
We consider the two-dimensional two-way wave in  \eqref{eq:FracLayer4ScalarWaveEq_2D} with the propagation speed $c=1$. We let $\Omega=[-2,2]\times[-2,2]$ be the interior domain with a buffer layer of length $\delta=0.5$ next to its boundaries. Thus, $(x,y)\in[-2.5,2.5]\times[-2.5,2.5]$. We assume smooth initial conditions $u_{0}(x,y)=e^{-5(x^2+y^2)}$ and $\varphi(x,y)=0$. We choose the variable-order functions
\begin{equation}\label{eq:Ord4Wave_2Dx}
   \alpha(x)
   =\left\{
   \begin{array}{ll}
   2, &x\in[x_L-\delta,x_L-\bar{\delta}),
   \\[3pt]
   1.5-(0.5-\epsilon)\tanh(20(x-x_L+\bar{\delta})), &x\in[x_L-\bar{\delta},x_L),
   \\[3pt]
   1+\epsilon, &x\in[x_L,x_R],
   \\[3pt]
   1.5-(0.5-\epsilon)\tanh(20(x_R+\bar{\delta}-x)), &x\in(x_R, x_R+\bar{\delta}],
   \\[3pt]
   2, &x\in(x_R+\bar{\delta}, x_R+\delta],
   \end{array}
   \right.
\end{equation}
and
\begin{equation}\label{eq:Ord4Wave_2Dy}
   \beta(y)
   =\left\{
   \begin{array}{ll}
   2, &y\in[y_L-\delta,y_L-\bar{\delta})
   \\[3pt]
   1.5-(0.5-\epsilon)\tanh(20(y-y_L+\bar{\delta})), &y\in[y_L-\bar{\delta},y_L),
   \\[3pt]
   1+\epsilon, &y\in[y_L,y_R],
   \\[3pt]
   1.5-(0.5-\epsilon)\tanh(20(y_R+\bar{\delta}-y)), &y\in(y_R, y_R+\bar{\delta}],
   \\[3pt]
   2, &y\in(y_R+\bar{\delta},y_R+\delta],
   \end{array}
   \right.
\end{equation}
with $\epsilon=10^{-5}$ and $\bar{\delta}=\frac{\delta}{2}$. The FBL in this setting is shown in Fig. \ref{Fig:VariableOrd}. 
\end{example}

% a similar smooth function defined as 
% \begin{equation}
%     \sigma_{x}(x)
%     =\left\{
%     \begin{array}{ll}
%          1, &x\in[x_L-\delta,x_L-\bar{\delta}),
%          \\[3pt]
%          0.5-0.5\tanh(20(x-x_L+\bar{\delta})), &x\in[x_L-\bar{\delta},x_L),
%          \\[3pt]
%          0, &x\in[x_L,x_R],
%          \\[3pt]
%          0.5-0.5\tanh(20(x_R+\bar{\delta}-x)), &x\in(x_R, x_R+\bar{\delta}],
%          \\[3pt]
%          1, &x\in(x_R+\bar{\delta}, x_R+\delta],
%     \end{array}
%     \right.
% \end{equation}

\begin{figure}
\centering 
$t=1$\\[-13 pt]
\noindent\rule{6cm}{0.7pt}
\\
\includegraphics[width=0.4 \linewidth]{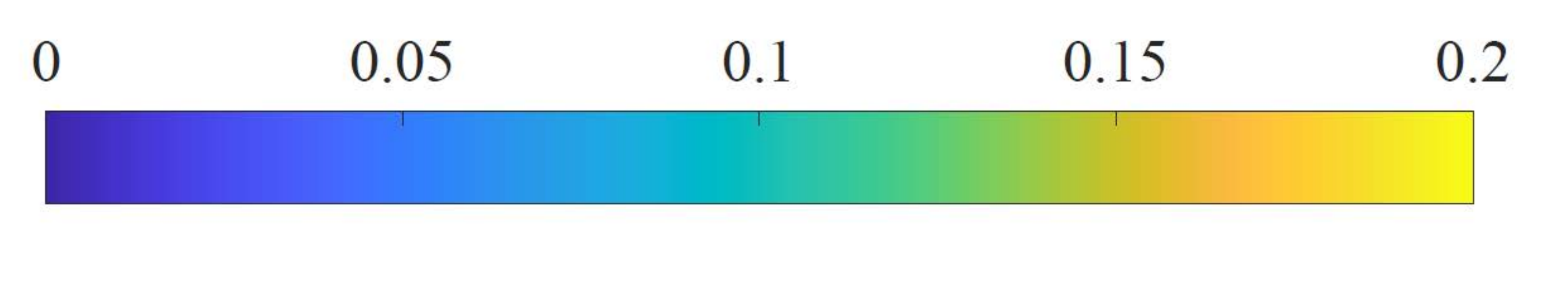}
\\
\includegraphics[clip, trim=0.5cm 0cm 1.7cm 0.8cm, width=0.25\linewidth]{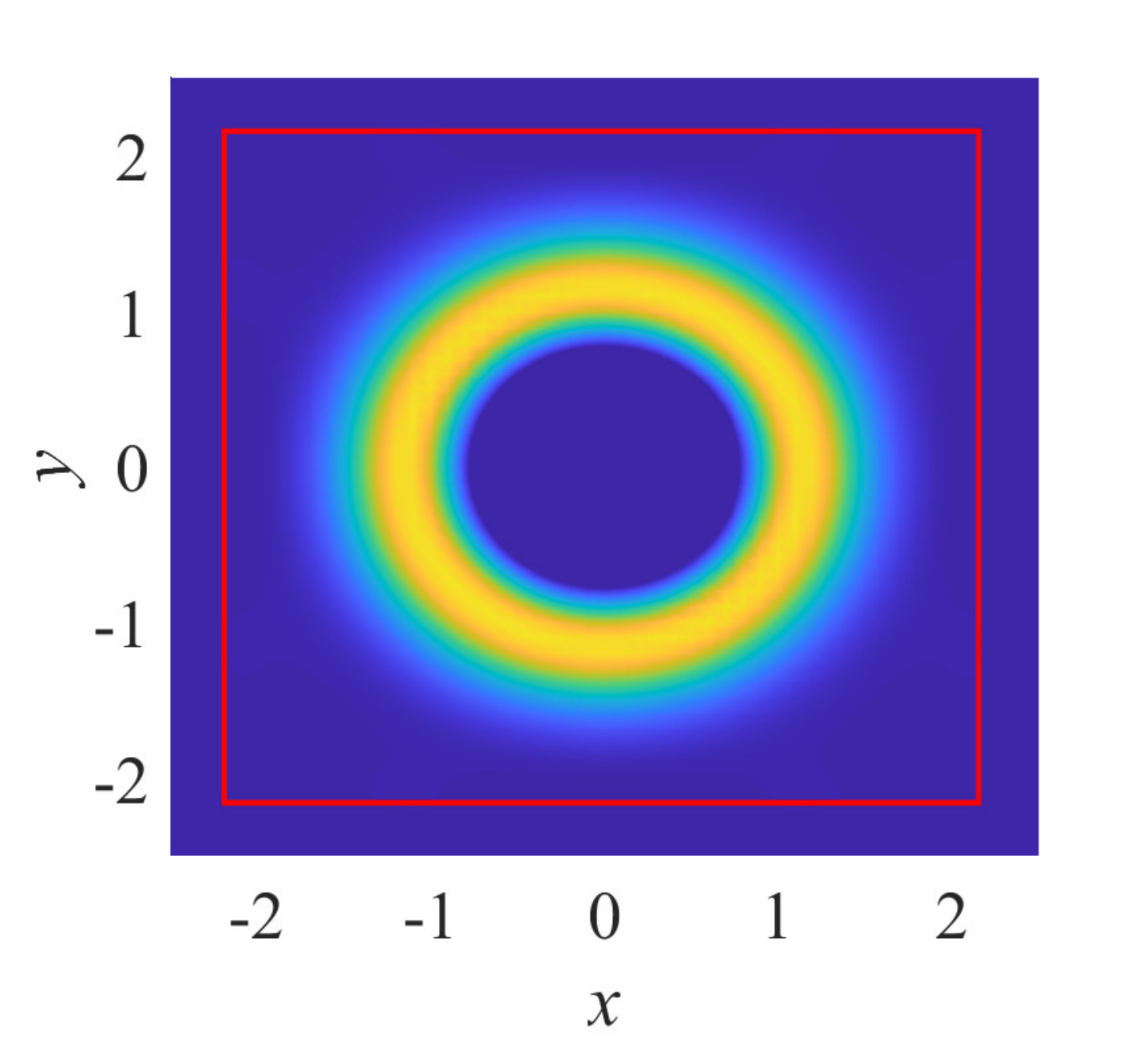}
\includegraphics[clip, trim=0.5cm 0cm 1.7cm 0.8cm, width=0.25\linewidth]{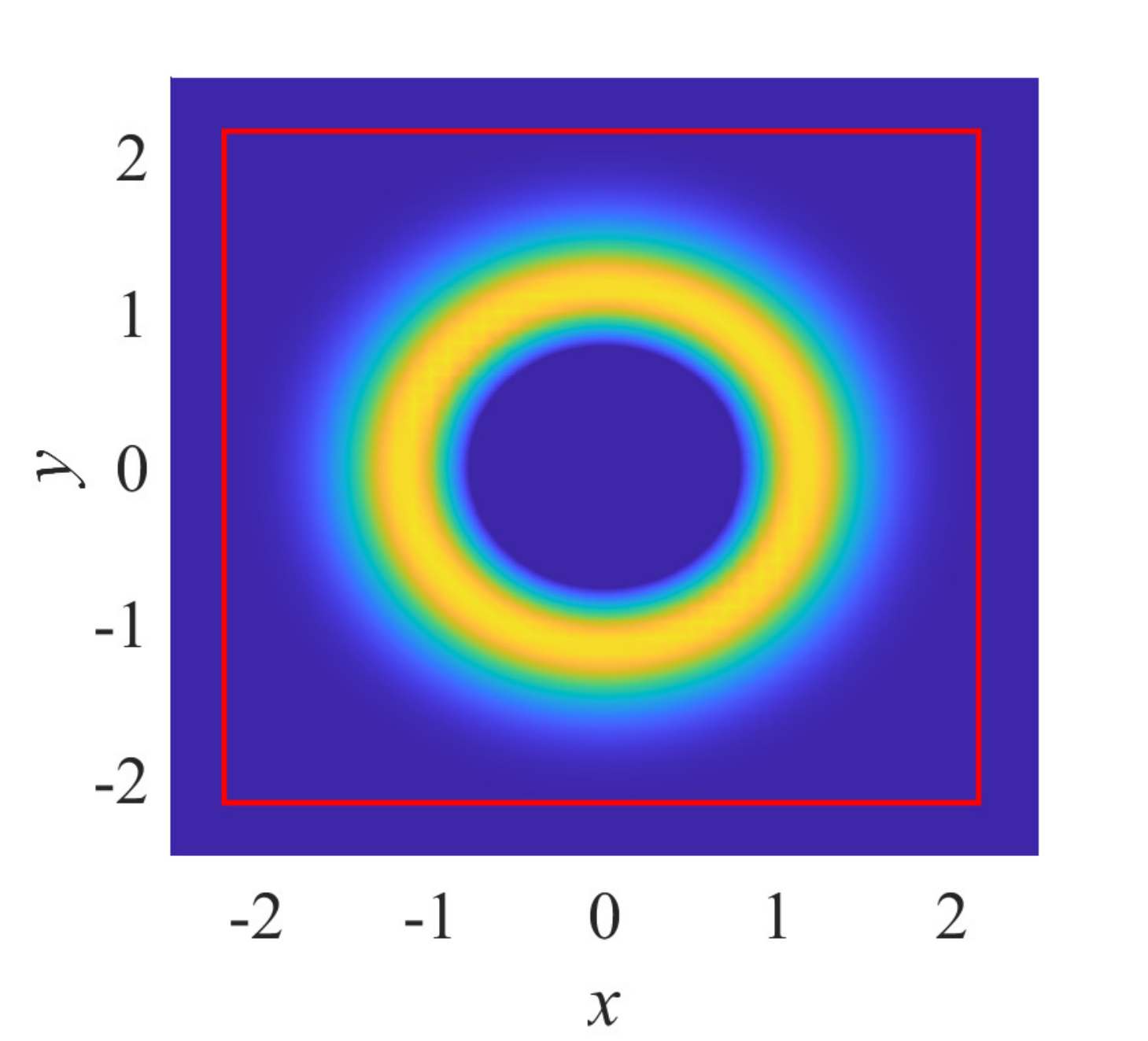}
\includegraphics[clip, trim=0.5cm 0cm 1.7cm 0.8cm, width=0.25\linewidth]{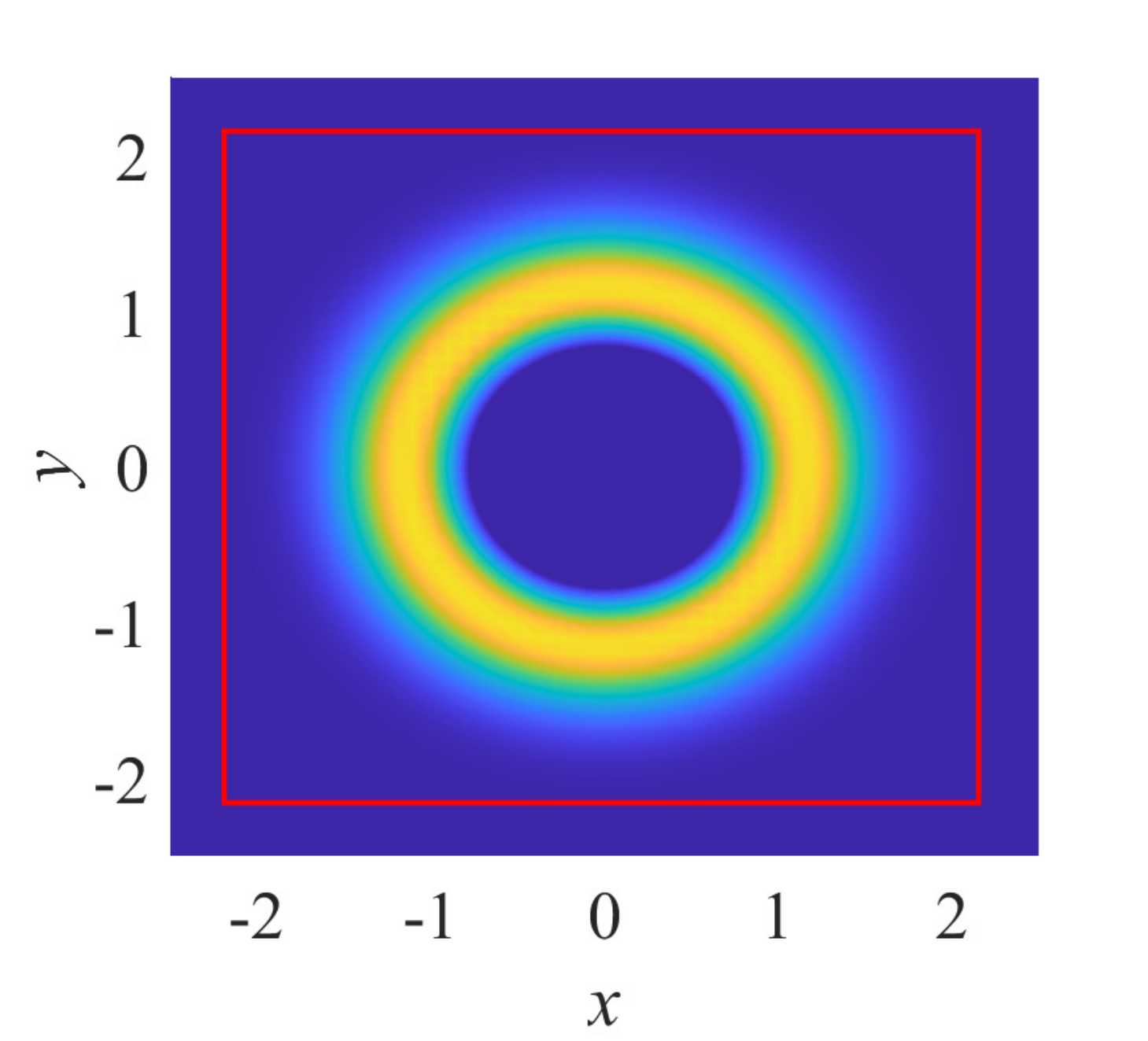}
\\ [-5 pt]
$t=3$\\[-13 pt]
\noindent\rule{6cm}{0.7pt}
\\
\includegraphics[width=0.4 \linewidth]{ColorBar.pdf}
\\
\includegraphics[clip, trim=0.5cm 0cm 1.7cm 0.8cm, width=0.25\linewidth]{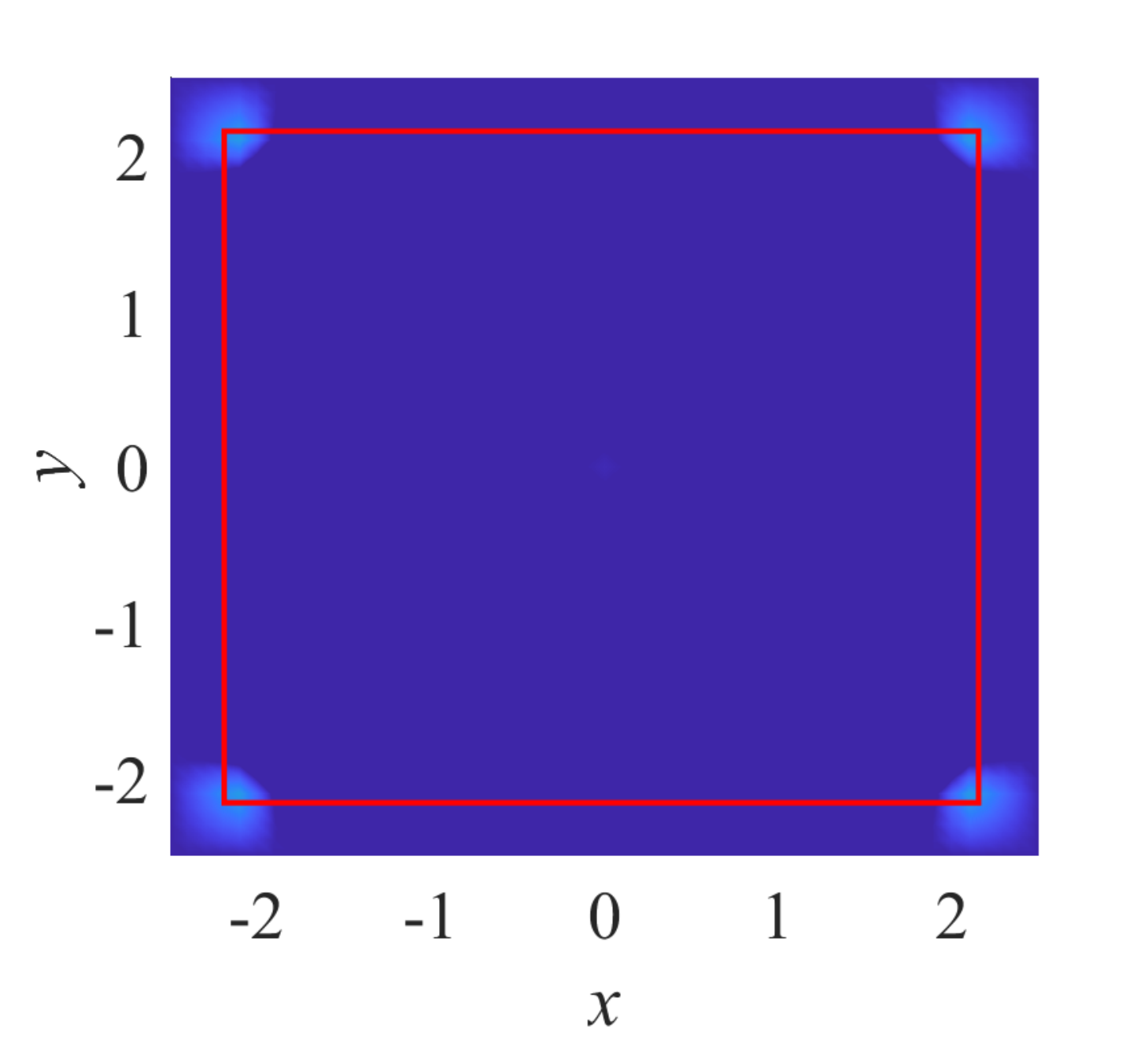}
\includegraphics[clip, trim=0.5cm 0cm 1.7cm 0.8cm, width=0.25\linewidth]{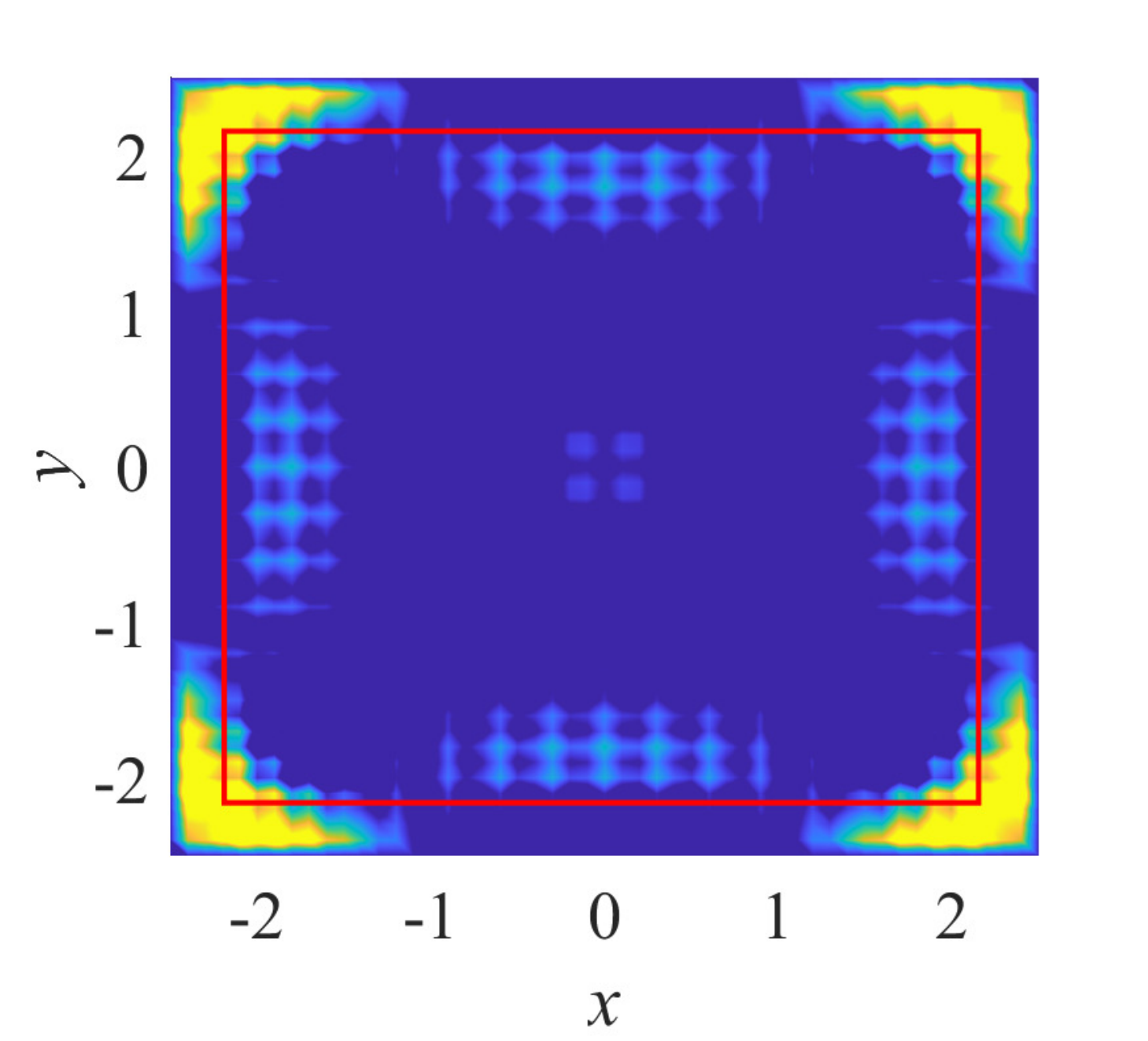}
\includegraphics[clip, trim=0.5cm 0cm 1.7cm 0.8cm, width=0.25\linewidth]{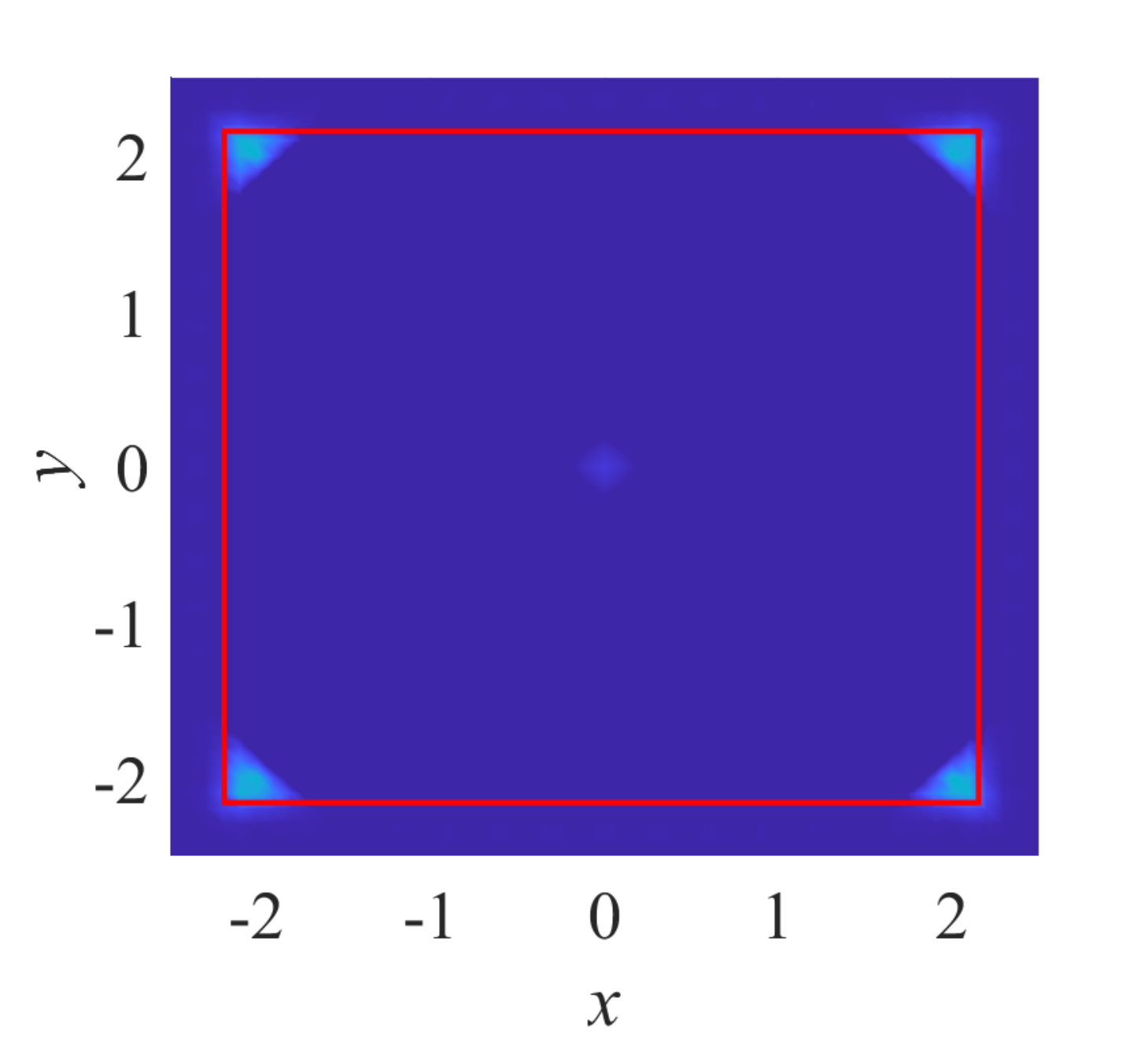}
\\[-5 pt] 
$t=5$\\[-13 pt]
\noindent\rule{6cm}{0.7pt}
\\
\includegraphics[width=0.4 \linewidth]{ColorBar.pdf}
\\
\includegraphics[clip, trim=0.5cm 0cm 1.7cm 0.8cm, width=0.25\linewidth]{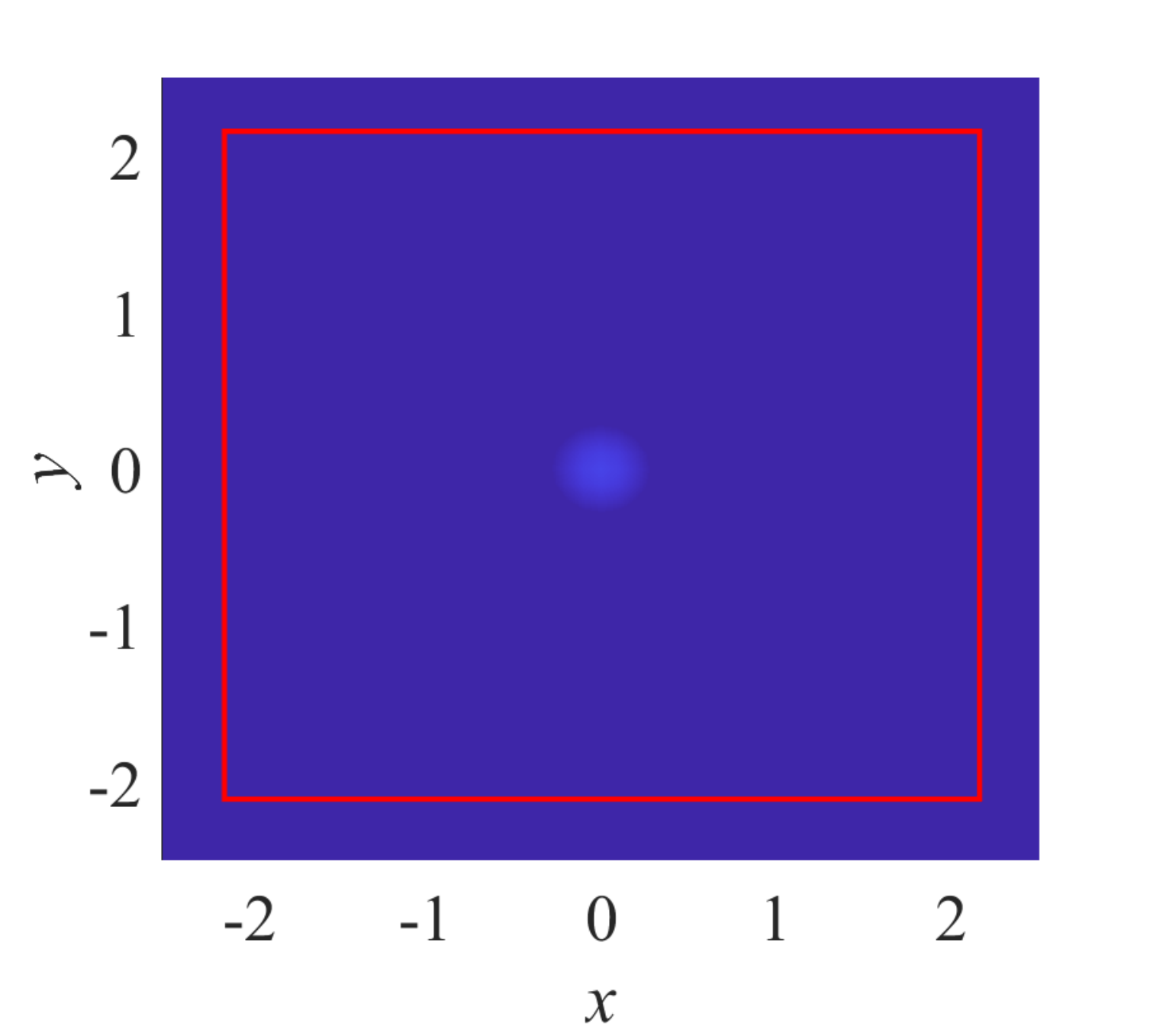}
\includegraphics[clip, trim=0.5cm 0cm 1.7cm 0.8cm, width=0.25\linewidth]{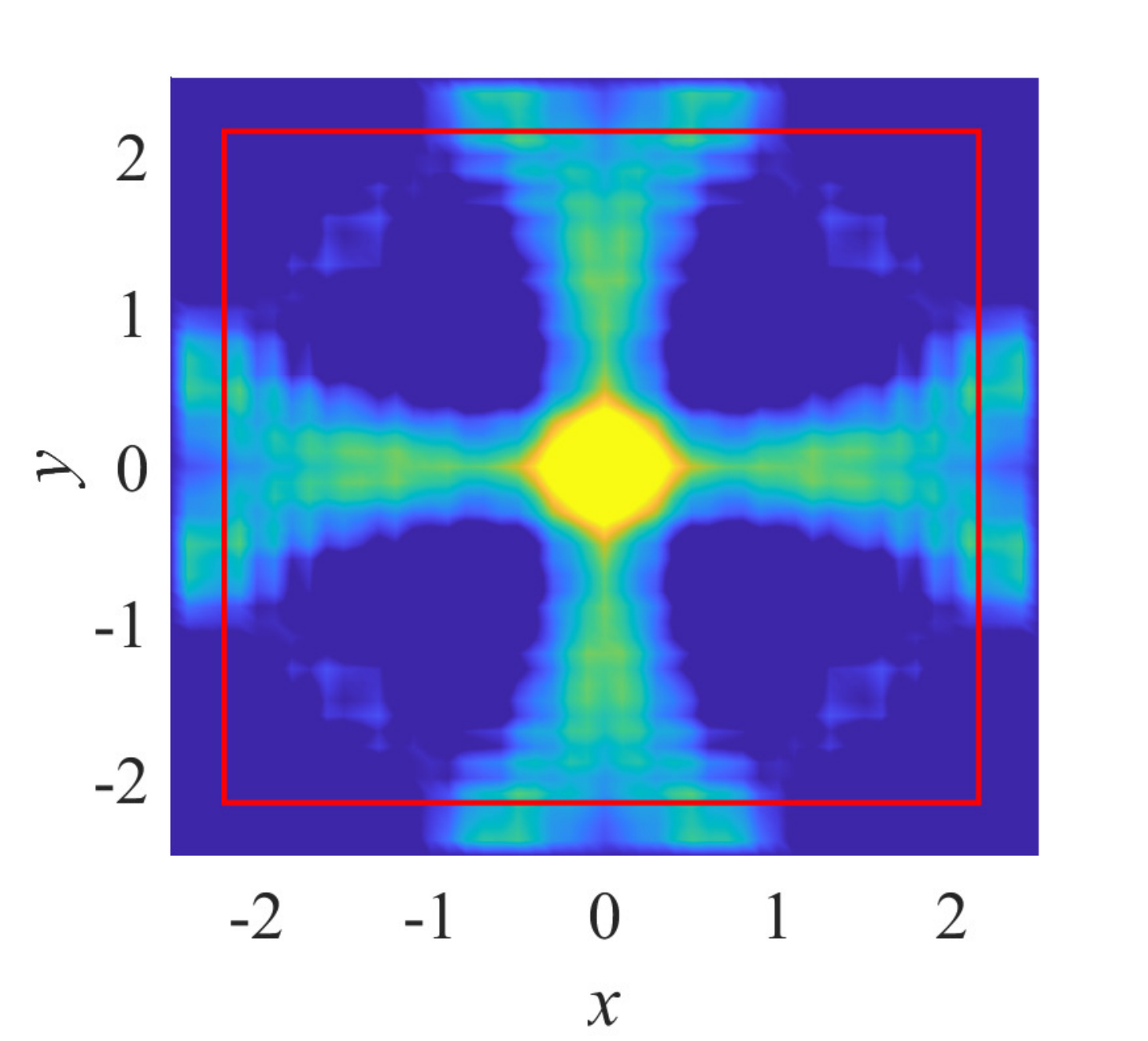}
\includegraphics[clip, trim=0.5cm 0cm 1.7cm 0.8cm, width=0.25\linewidth]{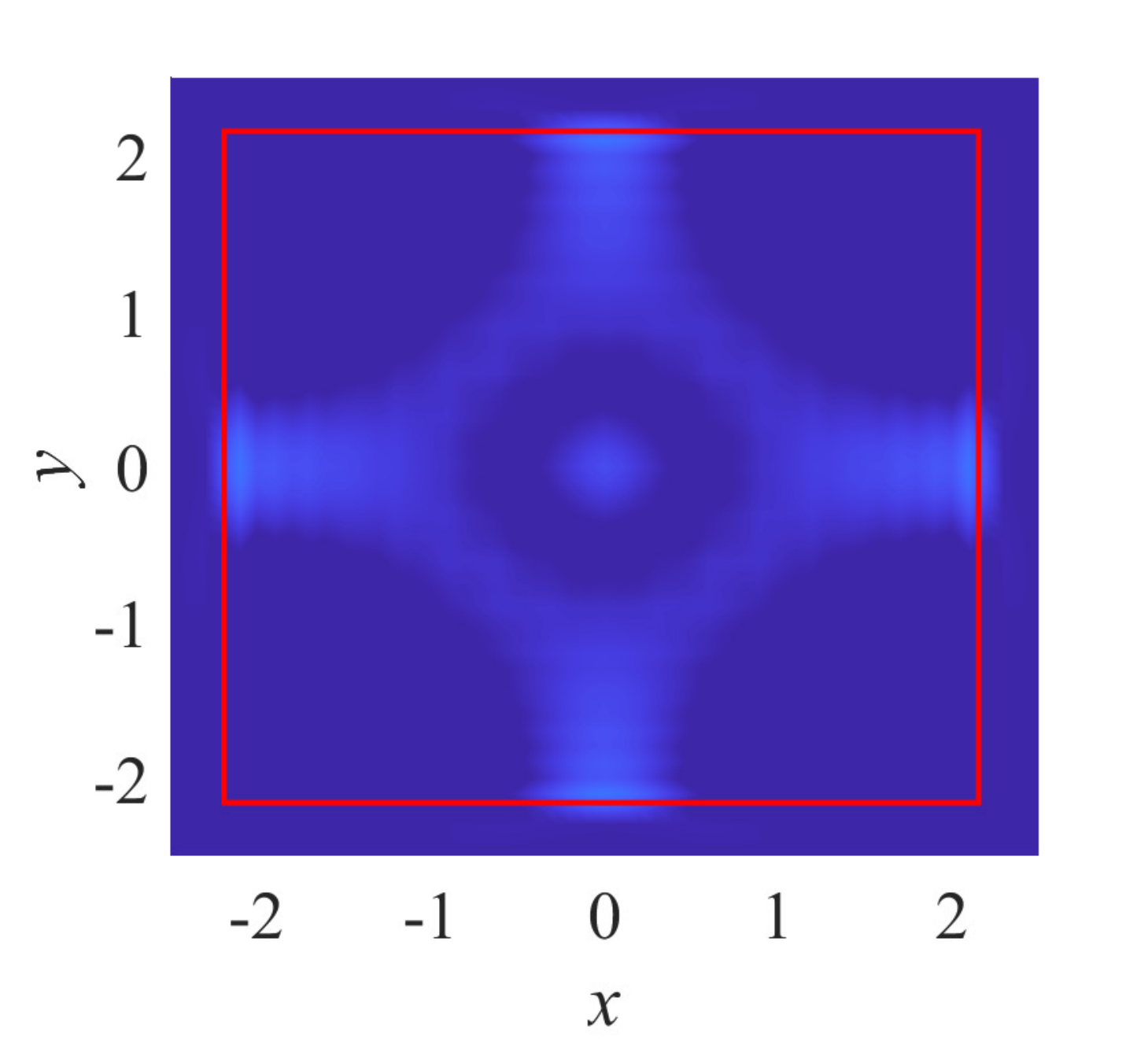}
\vspace{-0.2 in}
%%%%%%%%%%
\caption{Two-dimensional two-way wave (Example \ref{Ex: FBL two way 2D}): comparison of FBL (left), PML I (middle), and PML II (right) at $t=1, 3, 5$. The red solid rectangles indicate the boundaries of interior domain.
}\label{Fig:FBL_PML_2D}
\end{figure}

\begin{figure}
\centering
\includegraphics[clip, trim=1cm 1cm 0cm 0cm, width=0.42\linewidth]{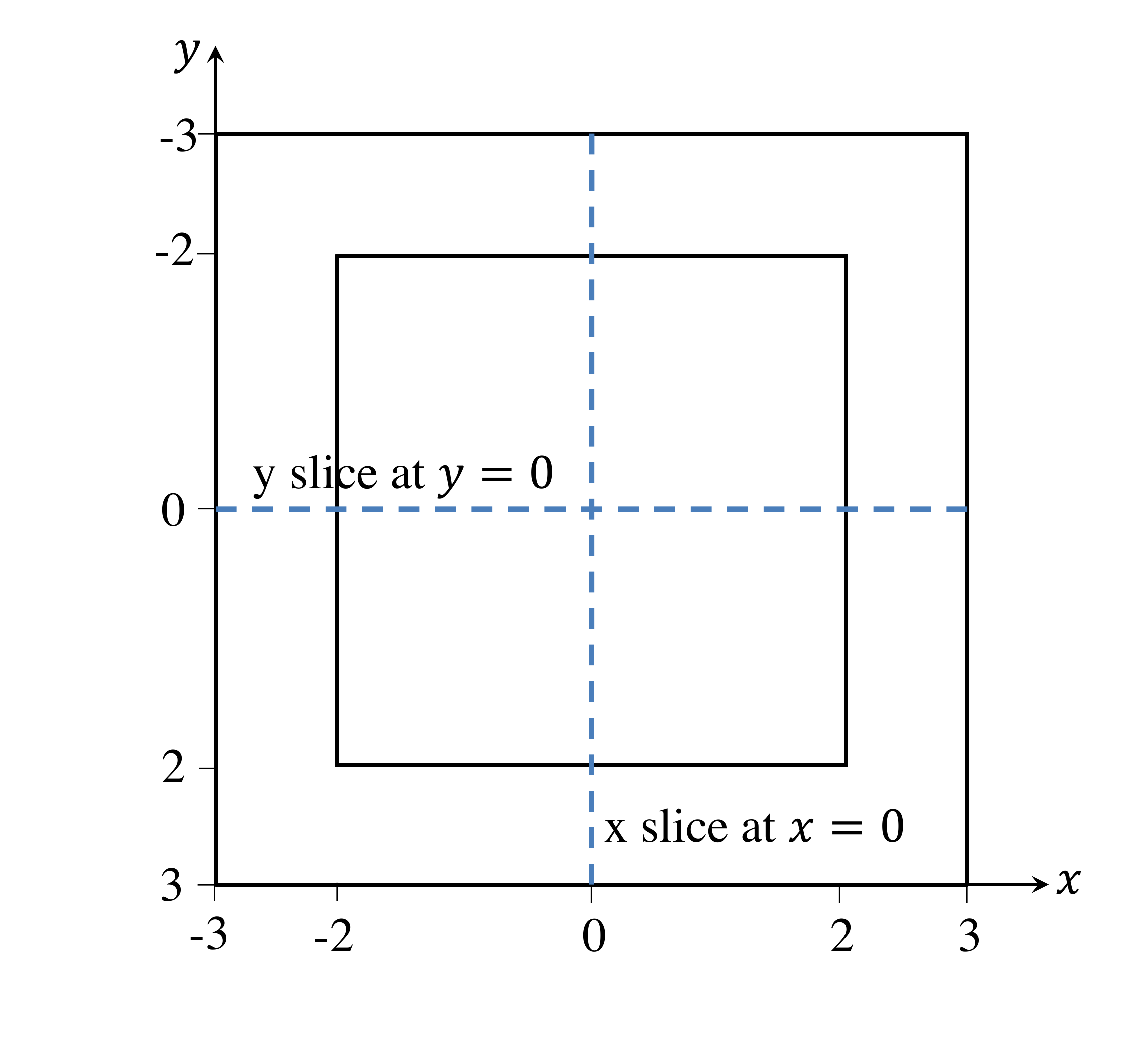}
\\
$t=1$\\[-9 pt]
\noindent\rule{4cm}{0.7pt}
\\
\includegraphics[clip, trim=0.5cm 0cm 1.7cm 0cm, width=0.42\linewidth]{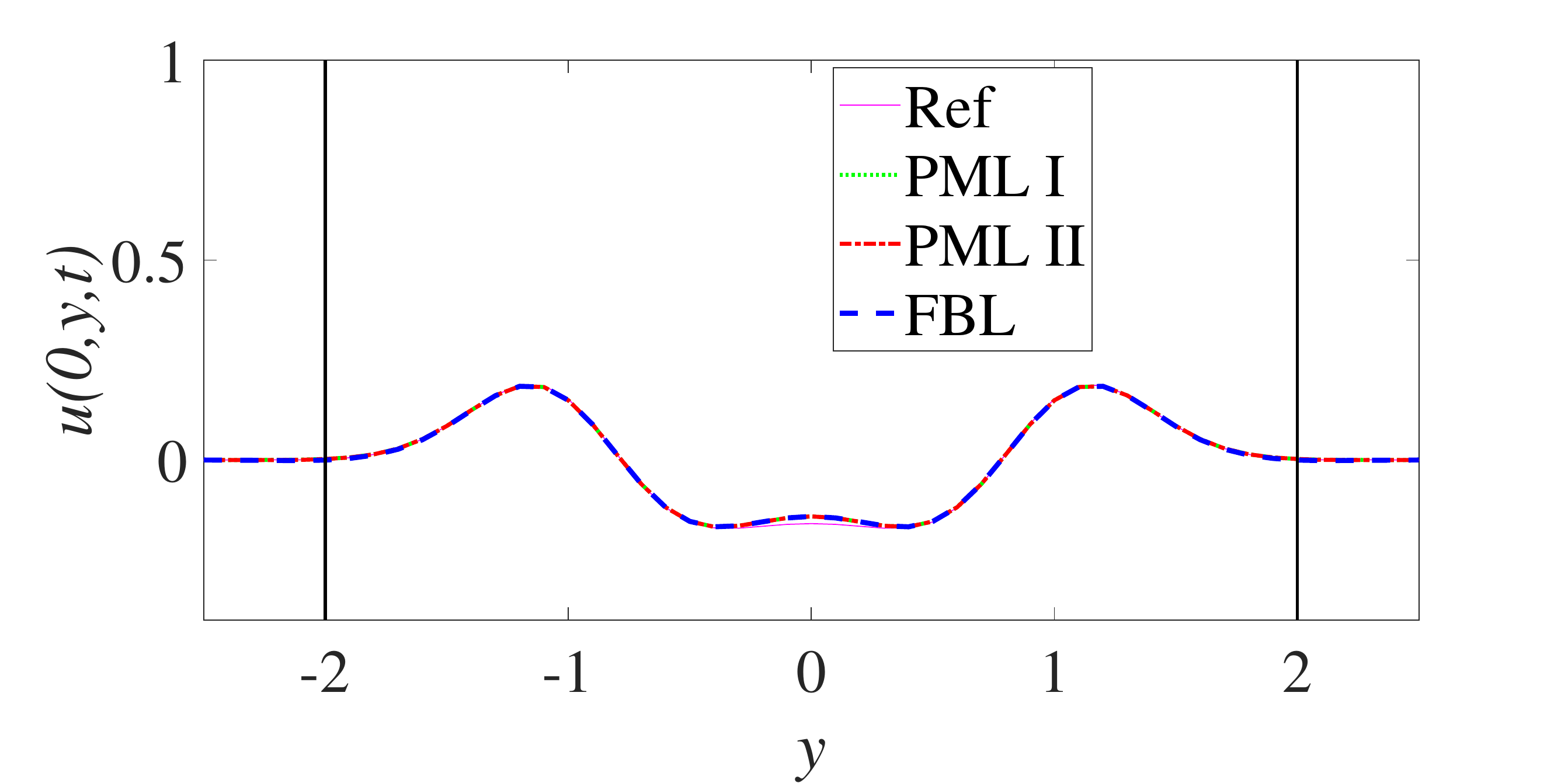}
\includegraphics[clip, trim=0.5cm 0cm 1.7cm 0cm, width=0.42\linewidth]{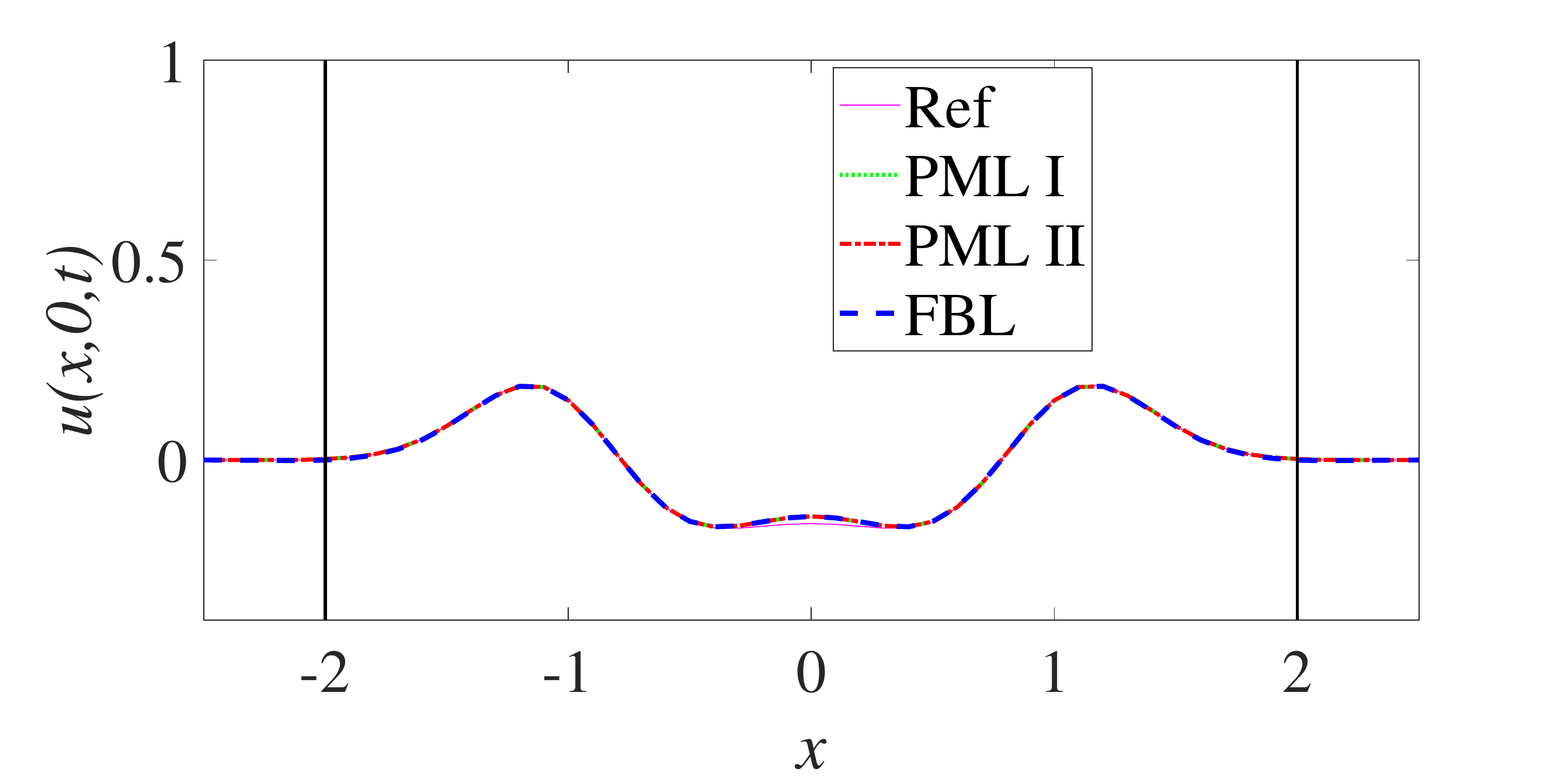}
\\
$t=3$\\[-9 pt]
\noindent\rule{4cm}{0.7pt}
\\
\includegraphics[clip, trim=0.5cm 0cm 1.7cm 0cm, width=0.42\linewidth]{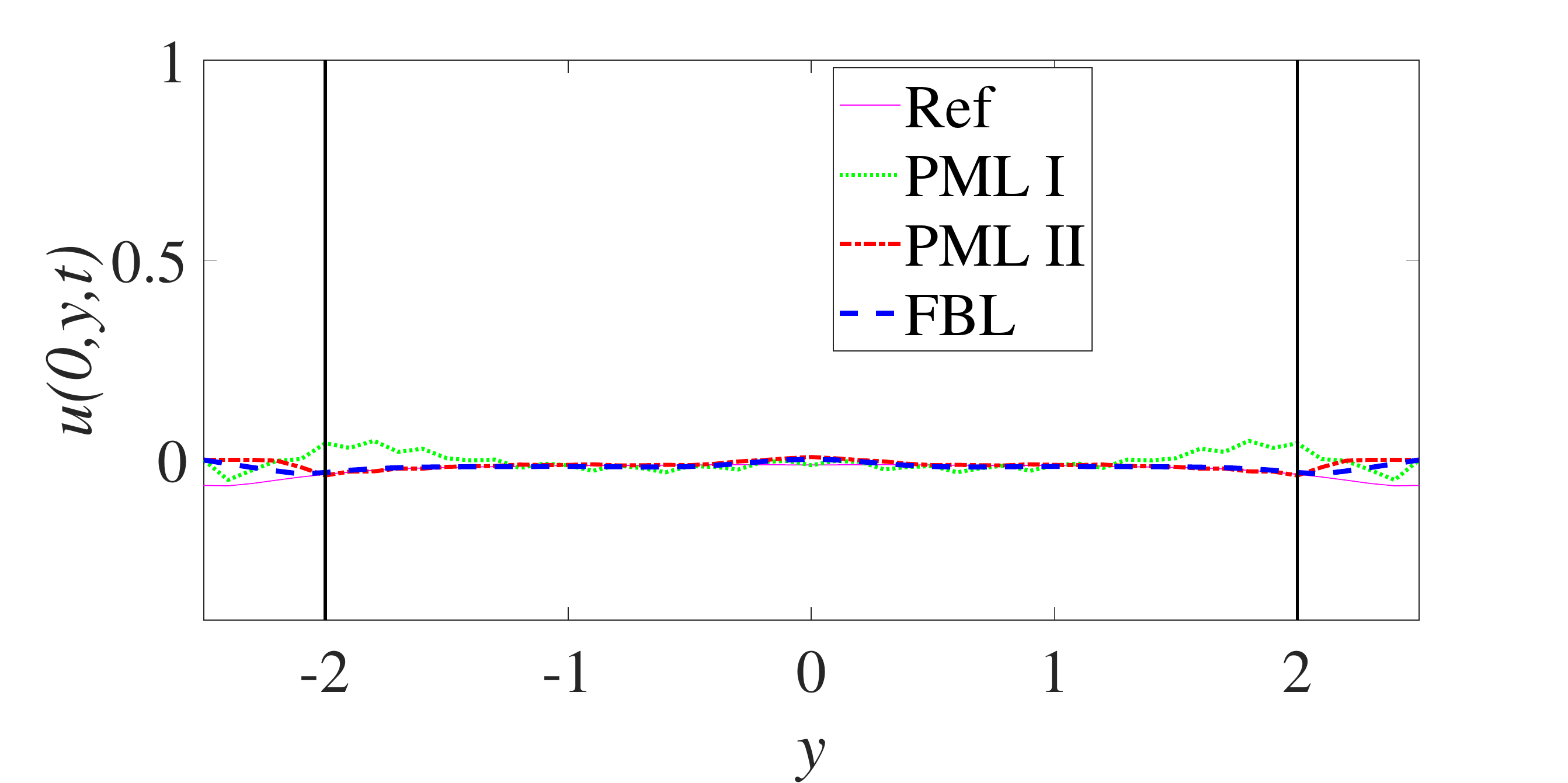}
\includegraphics[clip, trim=0.5cm 0cm 1.7cm 0cm, width=0.42\linewidth]{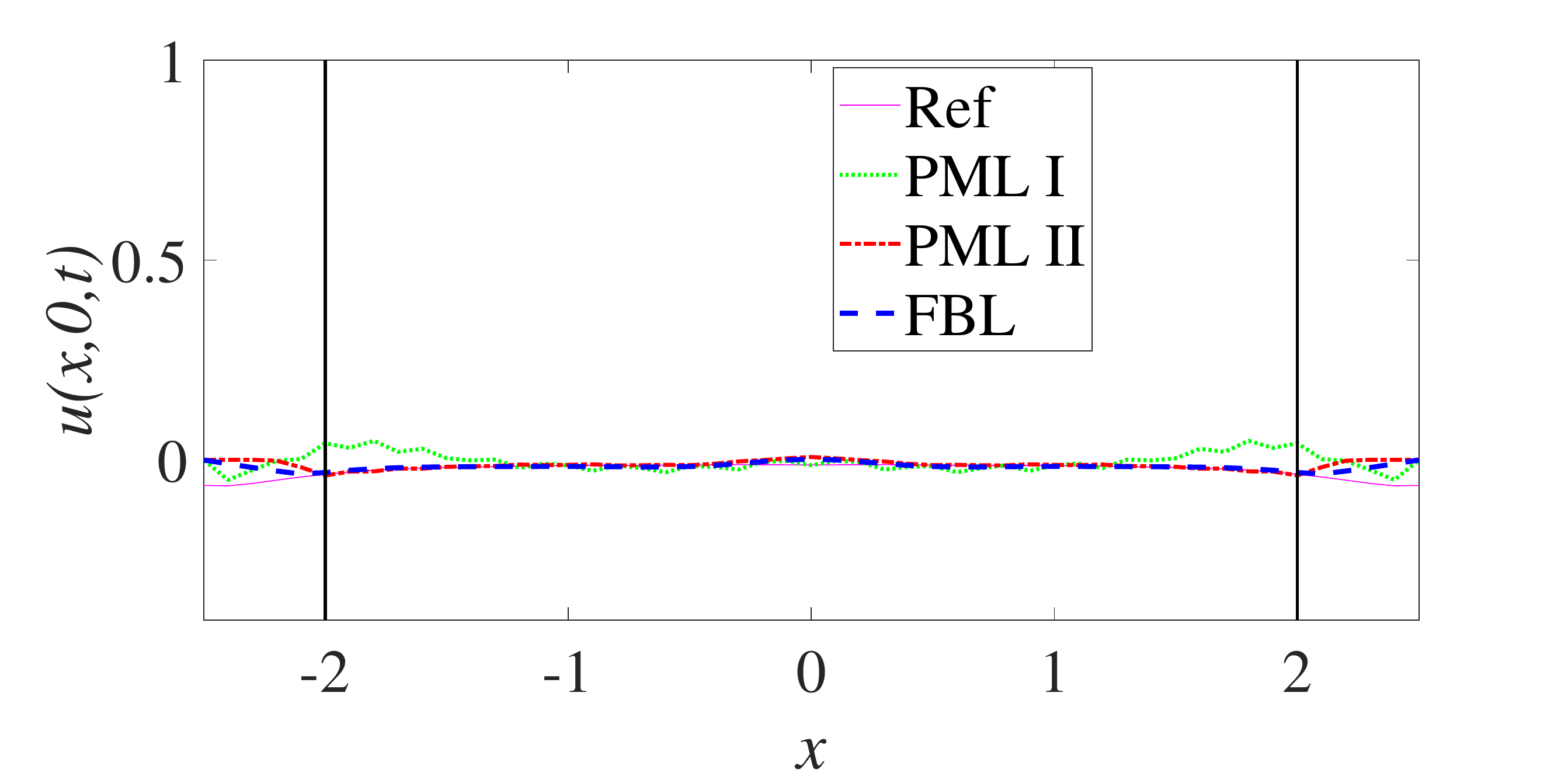}
\\
$t=5$\\[-9 pt]
\noindent\rule{4cm}{0.7pt}
\\
\includegraphics[clip, trim=0.5cm 0cm 1.7cm 0cm, width=0.42\linewidth]{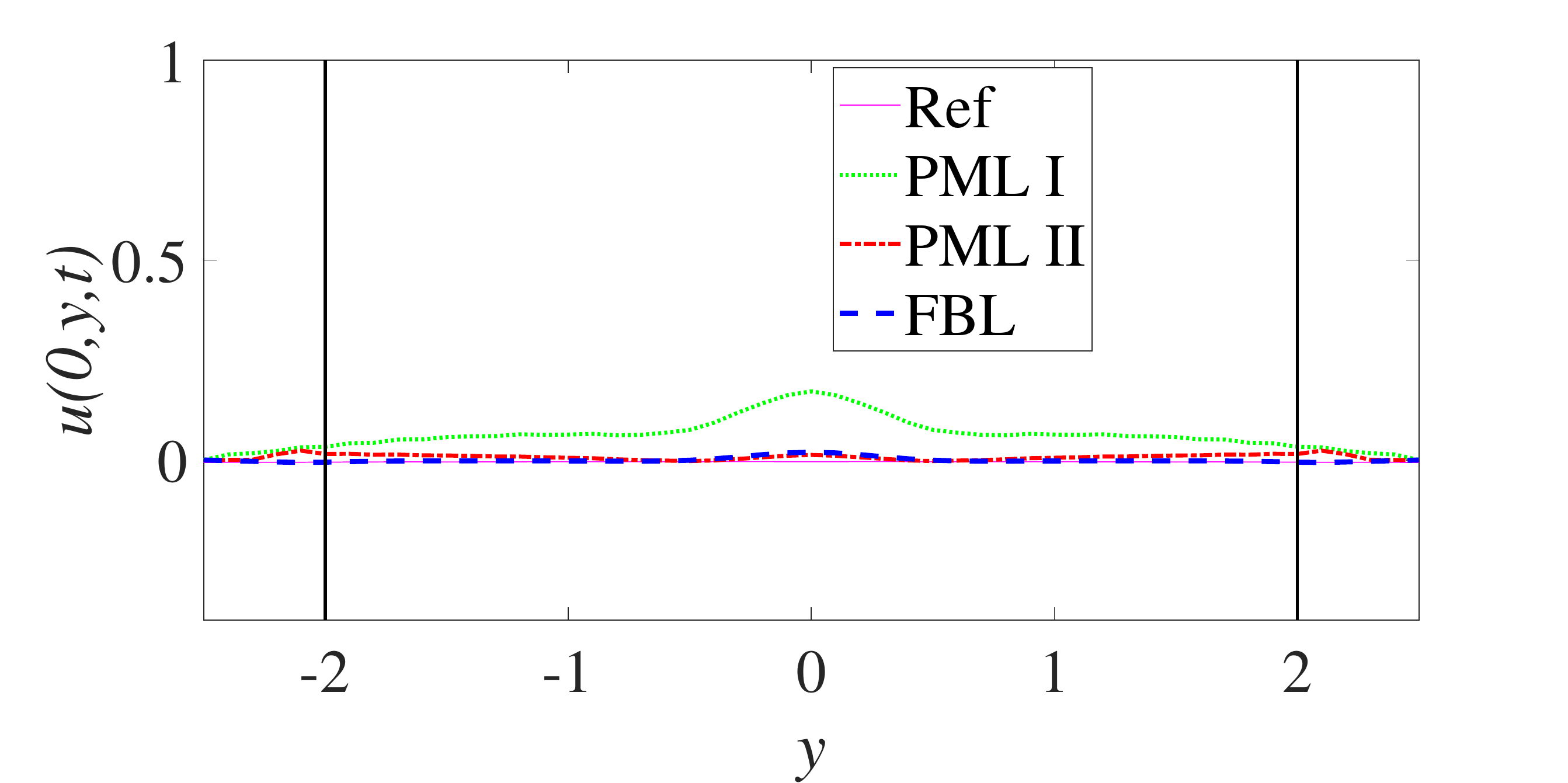}
\includegraphics[clip, trim=0.5cm 0cm 1.7cm 0cm, width=0.42\linewidth]{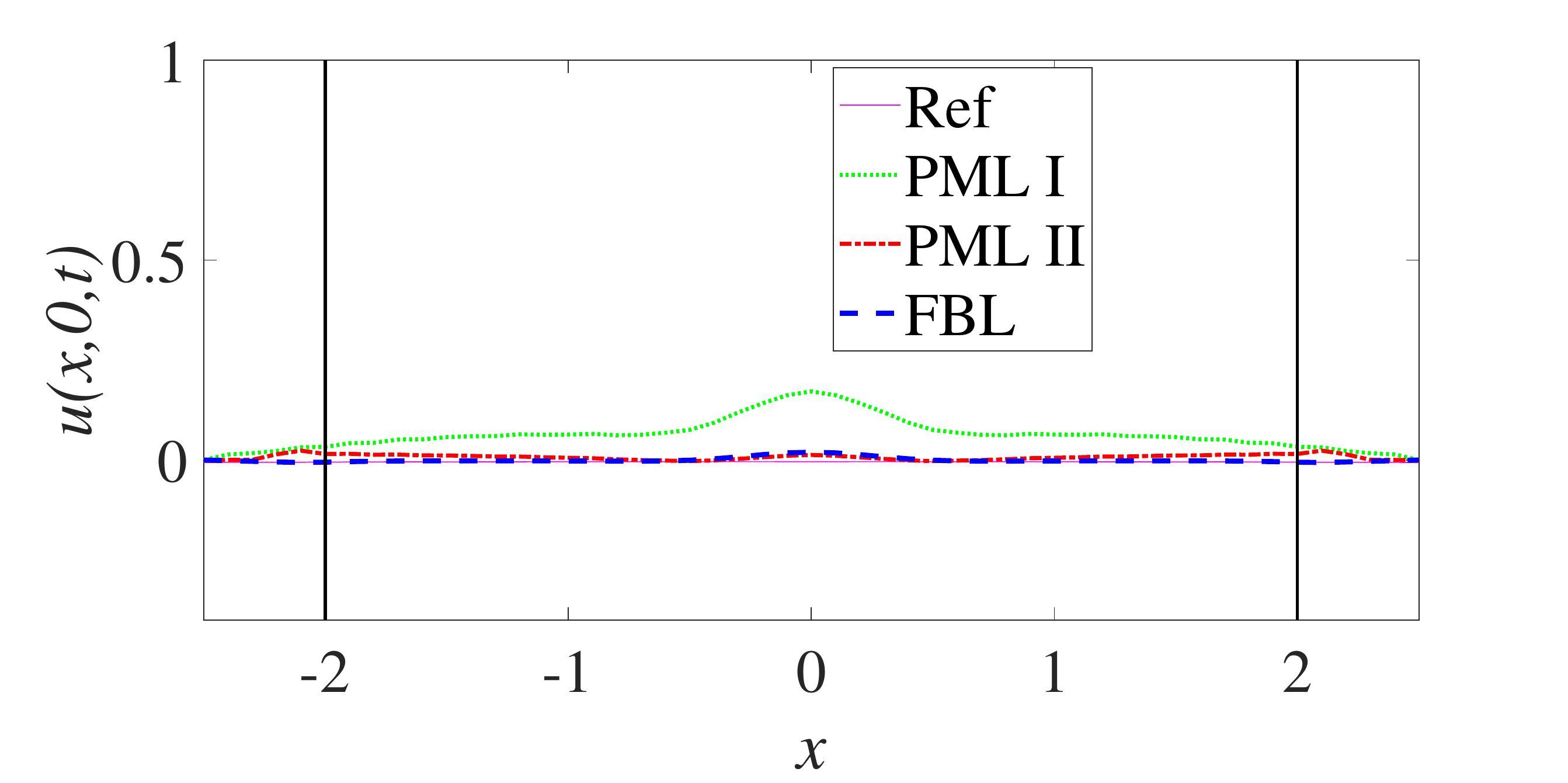}
%%%%%%%%%%
\caption{Two-dimensional two-way wave (Example \ref{Ex: FBL two way 2D}): comparison of FBL, PML I, PML II, and reference solution. Left column: $x$-slices for fixed $x=0$. Right column: $y$-slices for fixed $y=0$.
The black solid lines indicate the boundaries of interior domain.
}
\label{Fig:Slices_2D}
\end{figure}

We compare the FBL formulation with the PML formulations, namely, PML I given by \eqref{eq:PML4WaveEq_2D_I} and PML II given by \eqref{eq:PML4WaveEq_2D} in \ref{Sec:PML_2D}. Figure \ref{Fig:FBL_PML_2D} shows the FBL (left), PML I (middle), and PML II (right) solutions at different snapshots at $t=1$, $t=3$, and $t=5$. We can clearly see the corner reflections in the PML solution when the wave reaches to the boundaries of the extended domain while the FBL absorbs the wave with no reflections. We also plot the $x$- and $y$-slices of the FBL, PML I, PML II, and the reference solutions in Fig. \ref{Fig:Slices_2D}. The left and right columns give the  $x$-slices for fixed $x=0$ and $y$-slices for fixed $y=0$, respectively, at $t=1$, $t=3$, and $t=5$. It is evident that the proposed FBL performs better than the PML at absorbing the coming waves without reflection and corner issues. The reference solution is obtained by numerically solving \eqref{eq:WaveEq_2D} on a bigger domain $[-5,5]\times[-5,5]$ by spectral collocation method in space with $P_x=P_y=50$ and the central finite difference formula for time integration with $\tau=10^{-5}$. 

% In addition to no reflection and no corner problem, 
Another advantage of the proposed FBL formulation \eqref{eq:FracLayer4ScalarWaveEq_2D} is the computational efficiency with less CPU time. The FBL formulation has three equations in the coupled system \eqref{eq:FracLayer4ScalarWaveEq_2D} and thus theoretically requires less CPU time compared to the PML formulations which have more equations. We report that the CPU time for the FBL formulation \eqref{eq:FracLayer4ScalarWaveEq_2D}, the PML I in \eqref{eq:PML4WaveEq_2D_I}, and  the PML II in  \eqref{eq:PML4WaveEq_2D} are 175, 207, and 230 seconds, respectively. The numerical simulations are performed by Matlab R2019b on an Intel core i3 1.70 GHz CPU.

% The CPU time for Eq. \eqref{eq:FracLayer4ScalarWaveEq_2D} is theoretically less than that for Eq. \eqref{eq:PML4WaveEq_2D}. 

% In the present numerical simulation performed by Matlab R2019b on a Thinkpad E450s laptop with 4GB RAM and an Intel i3 core at 1.70 GHz, the CPU time for the time integration in Eq. \eqref{eq:FracLayer4ScalarWaveEq_2D} is 175 seconds while the time integration in Eq. \eqref{eq:PML4WaveEq_2D} costs 205 seconds. 

%%%%%%%%%%%%%%%
\section{Summary}\label{sec:Conclusion}
%%%%%%%%%%%%%%%
In this paper, we consider the application of variable-order space fractional differential operators in constructing absorbing boundary layers for wave propagation in finite domains. The variable-order functions are solely functions of space. We formulate our method by extending the space domain by attaching fractional buffer layers (FBLs) to its boundaries. Then, we define the corresponding variable-order functions such that we recover the original equation (or its corresponding split equation) inside the interior domain of interest
and introduce dissipation in the buffer layer. This approach is compared with PML via several numerical simulations and is shown to be effective in eliminating reflections with no corner artifacts in two-dimensional problems. As an extension of this formulation in the future, we will consider developing FBLs to absorb waves in arbitrary bounded domains for higher dimensional 
and multiphysics problems. 

%%%%%%%%%%%%%%%%%%%%%%
\section*{Acknowledgements}
%%%%%%%%%%%%%%%%%%%%%%
This work was partially supported by China Scholarship Council (No. 201906890040). This work was also supported by the Department of Energy on the PhILMs project (DE-SC0019453) and by the MURI/ARO, USA on “Fractional PDEs for Conservation Laws and Beyond: Theory, Numerics and Applications” (W911NF-15-1-0562). We would like to thank Dr. Zhiping Mao in
the Division of Applied Mathematics at Brown University for his helpful
suggestions.



\newpage
\appendix
%%%%%%%%%%%%%%%
\section{Fully discrete schemes for FBLs}\label{Sec:FullyDiscreteScheme}
%%%%%%%%%%%%%%%
We derive fully discrete schemes for the proposed FBLs. Here, we employ a spectral collocation method for spatial discretization. The time integration is based on Crank-Nicolson formula or the second-order Adams-Bashforth method on the uniformed temporal grids $\{t_{n}=n\tau\}_{n=0}^{N}$ with $\tau=\frac{T}{N}$ being the time step.

%%%%%%%%%%%%%%%
\subsection{Spatial discretization in the  one-dimensional case\label{Sec:DifferentialMatrix_SCM1D}}
%%%%%%%%%%%%%%%
Assume that $\{x_{j}\}_{j=0}^{P}$ satisfying $x_{0}=x_L-\delta$ and $x_{P}=x_R+\delta$ are collocation points in the computational domain $[x_L-\delta,x_R+\delta]$. For $j=0,1,2,\ldots,P$, we denote by $L_{j}(x)$ the $j$-th order Lagrange interpolation polynomial  satisfying $L_{j}(x_{i})=\delta_{i,j}$. Let $\{L_{j}(x)\}_{j=0}^{P}$ be the set of basis functions and define the space
\begin{equation}
\mathbb{V}_{P}={\rm span}\{L_{j}(x), \ 0\leq j\leq P\}.
\end{equation}
For any $f_{P}(x)\in\mathbb{V}_{P}$, there holds the expansion
\begin{equation}
f_{P}(x)=\sum\limits_{j=0}^{P}f_P(x_{j})L_{j}(x).
\end{equation}
Assume that $f(x)$ satisfies homogeneous boundary conditions such that $f(x_{0})=f(x_{P})=0$. Taking the left- and right-sided Riemann-Liouville derivatives of the interpolation $f_{P}(x)$ at the collocation points $\{x_{j}\}_{j=1}^{P-1}$ gives
\begin{equation}
\prescript{RL}{x_L-\delta}{\mathcal{D}}_{\mathbf{x}}^{\alpha(\mathbf{x})}\mathbf{f}_{P}=\mathbf{D}^{\bm\alpha}_{L}\mathbf{f}_{P}
\end{equation}
and
\begin{equation}
\prescript{RL}{\mathbf{x}}{\mathcal{D}}_{x_R+\delta}^{\alpha(\mathbf{x})}\mathbf{f}_{P}=\mathbf{D}^{\bm\alpha}_{R}\mathbf{f}_{P}.
\end{equation}
Here $\mathbf{x}=[x_{1},\cdots,x_{P-1}]^{\top}$, $\mathbf{f}_{P}=f_{P}(\mathbf{x})$, $\left(\mathbf{D}^{\bm\alpha}_{L}\right)_{i,j}=\prescript{RL}{x_L-\delta}{\mathcal{D}}_{x_{i}}^{\alpha(x_{i})}L_{j}(x_{i})$, and $\left(\mathbf{D}^{\bm\alpha}_{R}\right)_{i,j}=\prescript{RL}{x_{i}}{\mathcal{D}}_{x_R+\delta}^{\alpha(x_{i})}L_{j}(x_{i})$.

In order to compute the entries of the differentiation matrices, we introduce the affine transformation $\hat{x}=\frac{2(x-x_L+\delta)}{(x_R-x_L+2\delta)}-1\in[-1,1]$ and recall that the Lagrange polynomials can be expanded as
\begin{equation}
L_{j}(\hat{x})=\sum\limits_{k=0}^{M'}\ell_{k}^{j}p^{a,b}_{k}(\hat{x}),
\ a,b>-1,
\end{equation}
with $p_{k}^{a,b}(\hat{x})$ being the Jacobi polynomial and the coefficients given by
\begin{equation}
\ell_{k}^{j}
=
\left\{
\begin{array}{ll}
\frac{p_{k}^{a,b}(\hat{x}_{j})\omega_{j}}{\gamma_{k}^{a,b}},
&k=0,1,\ldots,M'-1,
\\[3pt]
\frac{p_{M'}^{a,b}(\hat{x}_{j})\omega_{j}}{\left(2+\frac{a+b+1}{M'}\right)\gamma_{k}^{a,b}},
&k=M'.
\end{array}
\right.
\end{equation}
Here $\{\hat{x}_{j}\}_{j=0}^{M'}$ and $\{\omega_{j}\}_{j=0}^{M'}$ are the nodes and weights of the Jacobi-Gauss-Lobatto quadrature and $\gamma_{k}^{a,b}=\frac{2^{a+b+1}\Gamma(k+a+1)\Gamma(k+b+1)}{(2k+a+b+1)\Gamma(k+1)\Gamma(k+a+b+1)}$. Consequently, ${}_{RL}{\rm D}^{\alpha(x_{i})}_{x_L-\delta,x_{i}}L_{j}(x)$ and ${}_{RL}{\rm D}^{\alpha(x_{i})}_{x_{i},x_R+\delta}L_{j}(x)$  can be evaluated by
\begin{equation}
\begin{aligned}
&\prescript{RL}{x_L-\delta}{\mathcal{D}}_{x_{i}}^{\alpha(x_{i})}L_{j}(x_{i})
=\left(\frac{x_R-x_L+2\delta}{2}\right)^{-\alpha(x_{i})}\sum\limits_{k=0}^{M'}\ell_{k}^{j}\,\prescript{RL}{-1}{\mathcal{D}}_{\hat{x}_{i}}^{\alpha(x_{i})}p^{a,b}_{k}(\hat{x}_{i}), 
\end{aligned}
\end{equation}
and
\begin{equation}
\begin{aligned}
&\prescript{RL}{x_{i}}{\mathcal{D}}_{x_R+\delta}^{\alpha(x_{i})}L_{j}(x_{i})
=\left(\frac{x_R-x_L+2\delta}{2}\right)^{-\alpha(x_{i})}\sum\limits_{k=0}^{M'}\ell_{k}^{j}\,\prescript{RL}{\hat{x}_{i}}{\mathcal{D}}_{1}^{\alpha(x_{i})}p^{a,b}_{k}(\hat{x}_{i}).
\end{aligned}
\end{equation}
Let $M'=P$. We have
\begin{equation}
    \left\{
    \begin{aligned}
    &\left(\mathbf{D}^{\bm\alpha}_{L}\right)_{i,j}
    =\mathbf{P}^{a,b,i}_{L}{\bm\ell}^{j},
    \\[3pt]
    &\left(\mathbf{D}^{\bm\alpha}_{R}\right)_{i,j}
    =\mathbf{P}^{a,b,i}_{R}{\bm\ell}^{j},
    \end{aligned}
    \right.
\end{equation}
where ${\bm\ell}^{j}=[\ell_{0}^{j},\cdots,\ell_{P}^{j}]^{\top}$, $\mathbf{P}_{L}^{a,b,i}=[\widehat{p}^{a,b,i}_{L,0},\cdots,\widehat{p}^{a,b,i}_{L,P}]$, 
$\mathbf{P}_{R}^{a,b,i}=[\widehat{p}^{a,b,i}_{R,0},\cdots,\widehat{p}^{a,b,i}_{R,P}]$ with
\begin{equation}
    \widehat{p}^{a,b,i}_{L,j}=\left(\frac{x_R-x_L+2\delta}{2}\right)^{-\alpha(x_i)}
    \prescript{RL}{-1}{\mathcal{D}}_{\hat{x}_{i}}^{\alpha(x_{i})}p_{j}^{a,b}(\hat{x}_{i}),
\end{equation}
and
\begin{equation}
    \widehat{p}^{a,b,i}_{R,j}=\left(\frac{x_R-x_L+2\delta}{2}\right)^{-\alpha(x_i)}
    \prescript{RL}{\hat{x}_{i}}{\mathcal{D}}_{1}^{\alpha(x_{i})}p_{j}^{a,b}(\hat{x}_{i}).
\end{equation}

It remains to evaluate Riemann-Liouville derivatives of Jacobi polynomials. Utilizing recursive relation of Jacobi polynomials, the recursive formulae for the left-sided Riemann-Liouville derivative of Jacobi polynomials are given as \cite{Zeng@2012,Zeng@Arxiv2014},
\begin{equation}
\prescript{RL}{-1}{\mathcal{D}}_{\hat{x}}^{\alpha(x)}p_{0}^{a,b}(\hat{x})
=\frac{(\hat{x}+1)^{-\alpha(x)}}{\Gamma(1-\alpha(x))},
\end{equation}
\begin{equation}
\prescript{RL}{-1}{\mathcal{D}}_{\hat{x}}^{\alpha(x)}p_{1}^{a,b}(\hat{x})
=\frac{a+b+2}{2}\frac{(\hat{x}+1)^{1-\alpha(x)}}{\Gamma(2-\alpha(x))}
-(b+1)\frac{(\hat{x}+1)^{-\alpha(x)}}{\Gamma(1-\alpha(x))},
\end{equation}
and
\begin{equation}
\begin{aligned}
&\prescript{RL}{-1}{\mathcal{D}}_{\hat{x}}^{\alpha(x)}p_{j+1}^{a,b}(\hat{x})
=\frac{A_{j}^{a,b}}{1+(n-\alpha(x))A_{j}^{a,b}\hat{C}_{j}^{a,b}}\frac{{\rm d}^{n}}{{\rm d}\hat{x}^{n}}\left(\hat{x}\,\prescript{RL}{-1}{\mathcal{J}}_{\hat{x}}^{n-\alpha(x)}p_{j}^{a,b}(\hat{x})\right)
\\[3pt]
&-\frac{B_{j}^{a,b}+(n-\alpha(x))A_{j}^{a,b}\hat{B}_{j}^{a,b}}{1+(n-\alpha(x))A_{j}^{a,b}\hat{C}_{j}^{a,b}}\prescript{RL}{-1}{\mathcal{D}}_{\hat{x}}^{\alpha(x)}p_{j}^{a,b}(\hat{x})
\\[3pt]
&-\frac{C_{j}^{a,b}+(n-\alpha(x))A_{j}^{a,b}\hat{A}_{j}^{a,b}}{1+(n-\alpha(x))A_{j}^{a,b}\hat{C}_{j}^{a,b}}\prescript{RL}{-1}{\mathcal{D}}_{\hat{x}}^{\alpha(x)}p_{j-1}^{a,b}(\hat{x})
\\[3pt]
&+\frac{(n-\alpha(x))A_{j}^{a,b}(\hat{A}^{a,b}_{j}p_{j-1}^{a,b}(-1)+\hat{B}^{a,b}_{j}p_{j}^{a,b}(-1)+\hat{C}^{a,b}_{j}p_{j+1}^{a,b}(-1))}{\Gamma(1-\alpha(x))\left[1+(n-\alpha(x))A_{j}^{a,b}\hat{C}_{j}^{a,b}\right]}(\hat{x}+1)^{-\alpha(x)}.
\end{aligned}
\end{equation}
Here the coefficients are given by
\begin{equation}
\label{eq:recurrence4Jacobi(1)}
\left\{
\begin{aligned}
{}&A_{j}^{a,b}
=\frac{(2j+a+b+1)(2j+a+b+2)}
{2(j+1)(j+a+b+1)},\\
&B_{j}^{a,b}
=\frac{(b^{2}-a^{2})(2j+a+b+1)}
{2(j+1)(j+a+b+1)(2j+a+b)},\\
&C_{j}^{a,b}
=\frac{(j+a)(j+b)(2j+a+b+2)}
{(j+1)(j+a+b+1)(2j+a+b)},
\end{aligned}
\right.
\end{equation}
and
\begin{equation}
\left\{
\begin{aligned}
%{}&\widehat{A}_{1}^{a,b}=0,\\
&\widehat{A}_{j}^{a,b}
=\frac{-2(j+a)(j+b)}
{(j+a+b)(2j+a+b)(2j+a+b+1)},\\
&\widehat{B}_{j}^{a,b}
=\frac{2(a-b)}{(2j+a+b)(2j+a+b+2)},\\
&\widehat{C}_{j}^{a,b}
=\frac{2(j+a+b+1)}{(2j+a+b+1)(2j+a+b+2)}.
\end{aligned}
\right.
\end{equation}
We can similarly derive the right-sided Riemann-Liouville derivative for Jacobi polynomials as follows,
\begin{equation}
\prescript{RL}{\hat{x}}{\mathcal{D}}_{1}^{\alpha(x)}p_{0}^{a,b}(\hat{x})
=\frac{(1-\hat{x})^{-\alpha(x)}}{\Gamma(1-\alpha(x))},
\end{equation}
\begin{equation}
\prescript{RL}{\hat{x}}{\mathcal{D}}_{1}^{\alpha(x)}p_{1}^{a,b}(\hat{x})
=-\frac{a+b+2}{2}\frac{(1-\hat{x})^{1-\alpha(x)}}{\Gamma(2-\alpha(x))}
+(a+1)\frac{(1-\hat{x})^{-\alpha(x)}}{\Gamma(1-\alpha(x))},
\end{equation}
and
\begin{equation}
\begin{aligned}
&\prescript{RL}{\hat{x}}{\mathcal{D}}_{1}^{\alpha(x)}p_{j+1}^{a,b}(\hat{x})
=\frac{(-1)^{n}A_{j}^{a,b}}{1+(n-\alpha(x))A_{j}^{a,b}\hat{C}_{j}^{a,b}}\frac{{\rm d}^{n}}{{\rm d}\hat{x}^{n}}\left(\hat{x}\,\prescript{RL}{\hat{x}}{\mathcal{J}}_{1}^{n-\alpha(x)}p_{j}^{a,b}(\hat{x})\right)
\\[3pt]
&-\frac{B_{j}^{a,b}+(n-\alpha(x))A_{j}^{a,b}\hat{B}_{j}^{a,b}}{1+(n-\alpha(x))A_{j}^{a,b}\hat{C}_{j}^{a,b}}\prescript{RL}{\hat{x}}{\mathcal{D}}_{1}^{\alpha(x)}p_{j}^{a,b}(\hat{x})
\\[3pt]
&-\frac{C_{j}^{a,b}+(n-\alpha(x))A_{j}^{a,b}\hat{A}_{j}^{a,b}}{1+(n-\alpha(x))A_{j}^{a,b}\hat{C}_{j}^{a,b}}\prescript{RL}{\hat{x}}{\mathcal{D}}_{1}^{\alpha(x)}p_{j-1}^{a,b}(\hat{x})
\\[3pt]
&+\frac{(n-\alpha(x))A_{j}^{a,b}(\hat{A}^{a,b}_{j}p_{j-1}^{a,b}(1)+\hat{B}^{a,b}_{j}p_{j}^{a,b}(1)+\hat{C}^{a,b}_{j}p_{j+1}^{a,b}(1))}{\Gamma(1-\alpha(x))\left[1+(n-\alpha(x))A_{j}^{a,b}\hat{C}_{j}^{a,b}\right]}(1-\hat{x})^{-\alpha(x)}.
\end{aligned}
\end{equation}

%%%%%%%%%%%%%%%
\subsection{Spatial discretization in the two-dimensional case\label{Sec:DifferentialMatrix_SCM_2D}}
%%%%%%%%%%%%%%%
Let $f(x,y)$ be suitably smooth with $(x,y)\in[x_L-\delta,x_R+\delta]\times[y_L-\delta,y_R+\delta]$. Assume that $f(x,y)$ can be approximated by the series
\begin{equation}
    f_{P_x,P_y}(x,y)
    =\sum\limits_{i=0}^{P_x}\sum\limits_{j=0}^{P_y}f_{i,j}L_{i}(x)L_{j}(y).
\end{equation}
Here $\{L_{i}(x)\}_{i=0}^{P_x}$ and $\{L_{j}(y)\}_{j=0}^{P_y}$ are Lagrange interpolation basis corresponding to the collocation points $(x_{i},y_{j})$ with $\{x_L-\delta=x_{0}<x_{1}<\cdots<x_{P_x-1}<x_{P_x}=x_R+\delta\}$ and $\{y_L-\delta=y_{0}<y_{1}<\cdots<y_{P_y-1}<y_{P_y}=y_R+\delta\}$.  It is evident that $f(x_{i},y_{j})\approx f_{P_x,P_y}(x_i,y_j)=f_{i,j}$. Then there hold 
\begin{equation}
\left\{
\begin{aligned}
&\prescript{RL}{x_L}{\mathcal{D}}_{x}^{\alpha(x)}f(x,y)
\approx\sum\limits_{i=0}^{P_x}\sum\limits_{j=0}^{P_y}
f_{P_x,P_y}(x_i,y_j)L_{j}(y)\,\prescript{RL}{x_L}{\mathcal{D}}_{x}^{\alpha(x)}L_{i}(x),
\\[3pt]
&\prescript{RL}{x}{\mathcal{D}}_{x_R}^{\alpha(x)}f(x,y)
\approx\sum\limits_{i=0}^{P_x}\sum\limits_{j=0}^{P_y}
f_{P_x,P_y}(x_i,y_j)L_{j}(y)\,\prescript{RL}{x}{\mathcal{D}}_{x_R}^{\alpha(x)}L_{i}(x),
\end{aligned}
\right.
\end{equation}
and
\begin{equation}
\left\{
\begin{aligned}
&\prescript{RL}{y_L}{\mathcal{D}}_{y}^{\beta(y)}f(x,y)
\approx\sum\limits_{i=0}^{P_x}\sum\limits_{j=0}^{P_y}
f_{P_x,P_y}(x_i,y_j)L_{i}(x)\,\prescript{RL}{y_L}{\mathcal{D}}_{y}^{\beta(y)}L_{j}(y),
\\[3pt]
&\prescript{RL}{y}{\mathcal{D}}_{y_R}^{\beta(y)}f(x,y)
\approx\sum\limits_{i=0}^{P_x}\sum\limits_{j=0}^{P_y}
f_{P_x,P_y}(x_i,y_j)L_{i}(x)\,\prescript{RL}{y}{\mathcal{D}}_{y_R}^{\beta(y)}L_{j}(y).
\end{aligned}
\right.
\end{equation}
Furthermore, if $f(x,y)$ satisfies homogeneous Dirichlet boundary conditions, its Riemann-Liouville derivatives at the interior collocation points $(x_{i},y_{j})$ have the following matrix forms,
\begin{equation}
\begin{aligned}
&\begin{bmatrix}
\prescript{RL}{x_L-\delta}{\mathcal{D}}_{x_{1}}^{\alpha(x_{1})}f(x_{1},y_{1})
&\cdots
&\prescript{RL}{x_L-\delta}{\mathcal{D}}_{x_{1}}^{\alpha(x_{1})}f(x_{1},y_{P_y-1})
\\[3pt]
\vdots&\ddots&\vdots
\\[3pt]
\prescript{RL}{x_L-\delta}{\mathcal{D}}_{x_{P_x-1}}^{\alpha(x_{P_x-1})}f(x_{P_x-1},y_{1})
&\cdots
&\prescript{RL}{x_L-\delta}{\mathcal{D}}_{x_{P_x-1}}^{\alpha(x_{P_x-1})}f(x_{P_x-1},y_{P_y-1})
\end{bmatrix}
\approx
\widetilde{\mathbf{D}}^{\bm\alpha}_{L}\mathbf{f},
\end{aligned}
\end{equation}
\begin{equation}
\begin{bmatrix}
\prescript{RL}{x_{1}}{\mathcal{D}}_{x_R+\delta}^{\alpha(x_{1})}f(x_{1},y_{1})
&\cdots
&\prescript{RL}{x_{1}}{\mathcal{D}}_{x_R+\delta}^{\alpha(x_{1})}f(x_{1},y_{P_y-1})
\\[3pt]
\vdots&\ddots&\vdots
\\[3pt]
\prescript{RL}{x_{P_x-1}}{\mathcal{D}}_{x_R+\delta}^{\alpha(x_{P_x-1})}f(x_{P_x-1},y_{1})
&\cdots
&\prescript{RL}{x_{P_x-1}}{\mathcal{D}}_{x_R+\delta}^{\alpha(x_{P_x-1})}f(x_{P_x-1},y_{P_y-1})
\end{bmatrix}
\approx
\widetilde{\mathbf{D}}^{\bm\alpha}_{R}\mathbf{f},
\end{equation}
\begin{equation}
\begin{aligned}
&\begin{bmatrix}
\prescript{RL}{y_L-\delta}{\mathcal{D}}_{y_{1}}^{\beta(y_{1})}f(x_{1},y_{1})
&\cdots
&\prescript{RL}{y_L-\delta}{\mathcal{D}}_{y_{P_y-1}}^{\beta(y_{P_y-1})}f(x_{1},y_{P_y-1})
\\[3pt]
\vdots&\ddots&\vdots
\\[3pt]
\prescript{RL}{y_L-\delta}{\mathcal{D}}_{y_{1}}^{\beta(y_{1})}f(x_{P_x-1},y_{1})
&\cdots
&\prescript{RL}{y_L-\delta}{\mathcal{D}}_{y_{P_y-1}}^{\beta(y_{P_y-1})}f(x_{P_x-1},y_{P_y-1})
\end{bmatrix}
\approx
\mathbf{f}\left(\overline{\mathbf{D}}^{\bm\beta}_{L}\right)^{\top},
\end{aligned}
\end{equation}
and
\begin{equation}
\begin{bmatrix}
\prescript{RL}{y_{1}}{\mathcal{D}}_{y_R+\delta}^{\beta(y_{1})}f(x_{1},y_{1})
&\cdots
&\prescript{RL}{y_{P_y-1}}{\mathcal{D}}_{y_R+\delta}^{\beta(y_{P_y-1})}f(x_{1},y_{P_y-1})
\\[3pt]
\vdots&\ddots&\vdots
\\[3pt]
\prescript{RL}{y_{1}}{\mathcal{D}}_{y_R+\delta}^{\beta(y_{1})}f(x_{P_x-1},y_{1})
&\cdots
&\prescript{RL}{y_{P_y-1}}{\mathcal{D}}_{y_R+\delta}^{\beta(y_{P_y-1})}f(x_{P_x-1},y_{P_y-1})
\end{bmatrix}%_{(M_x-1)\times(P_y-1)}
\approx
\mathbf{f}\left(\overline{\mathbf{D}}^{\bm\beta}_{R}\right)^{\top}.
\end{equation}
Here $\widetilde{\mathbf{D}}^{\bm\alpha}_{L}$, $\widetilde{\mathbf{D}}^{\bm\alpha}_{R}$, $\overline{\mathbf{D}}^{\bm\beta}_{L}$, and $\overline{\mathbf{D}}^{\bm\beta}_{R}$ are fractional differential matrices with respect to $x$ and $y$ that introduced in the section \ref{Sec:DifferentialMatrix_SCM1D}. Also, $\mathbf{f}$ is defined by
\begin{equation}
\mathbf{f}
=
\begin{bmatrix}
f_{1,1}&f_{1,2}&\cdots&f_{1,P_y-2}&f_{1,P_y-1}
\\[3pt]
f_{2,1}&f_{2,2}&\cdots&f_{2,P_y-2}&f_{2,P_y-1}
\\[3pt]
\vdots&\vdots&\ddots&\vdots&\vdots
\\[3pt]
f_{P_x-2,1}&f_{P_x-2,2}&\cdots&f_{P_x-2,P_y-2}&f_{P_x-2,P_y-1}
\\[3pt]
f_{P_x-1,1}&f_{P_x-1,2}&\cdots&f_{P_x-1,P_y-2}&f_{P_x-1,P_y-1}
\end{bmatrix}.
\end{equation}

%%%%%%%%%%%%%%%
\subsection{Fully discrete schemes}
%%%%%%%%%%%%%%%
In order to discretize Eqs. \eqref{eq:FracAdvcLeftRL} and \eqref{eq:FracAdvcRightRL} with $|V|=1$, we assume that the solution $u(x,t)$ at $t=t_{n}$ has the following Lagrange interpolation,
\begin{equation}
    u(x,t_n)
    \approx
    u_{P}(x,t_n)
    =\sum\limits_{j=0}^{P}u_{P}(x_j,t_n)L_j(x).
\end{equation}
When the spatial discretization is based on the spectral collocation method and the time discretization is based on the Crank-Nicolson formula, the fully discrete schemes for Eqs. \eqref{eq:FracAdvcLeftRL} and \eqref{eq:FracAdvcRightRL} with $|V|=1$ are given by
\begin{equation}\label{eq:FullyDiscrete4LeftAdvec}
\frac{\mathbf{u}^{n+1}-\mathbf{u}^{n}}{\tau}
=\frac{1}{2}\mathbf{D}^{\bm\alpha}_{R}\left(\mathbf{u}^{n+1}+\mathbf{u}^{n}\right),
\end{equation}
\begin{equation}\label{eq:FullyDiscrete4RightAdvec}
\frac{\mathbf{u}^{n+1}-\mathbf{u}^{n}}{\tau}
=\frac{1}{2}\mathbf{D}^{\bm\alpha}_{L}\left(\mathbf{u}^{n+1}+\mathbf{u}^{n}\right).
\end{equation}
The corresponding matrices forms are
\begin{equation}
\mathbf{u}^{n+1}=\left[\mathbf{I}-\frac{\tau}{2}\mathbf{D}^{\bm\alpha}_{R}\right]^{-1}\left[\left(\mathbf{I}+\frac{\tau}{2}\mathbf{D}^{\bm\alpha}_{R}\right)\mathbf{u}^{n}\right],
\end{equation}
\begin{equation}
\mathbf{u}^{n+1}=\left[\mathbf{I}-\frac{\tau}{2}\mathbf{D}^{\bm\alpha}_{L}\right]^{-1}\left[\left(\mathbf{I}+\frac{\tau}{2}\mathbf{D}^{\bm\alpha}_{L}\right)\mathbf{u}^{n}\right],
\end{equation}
respectively. Here $\mathbf{u}^{n}=[u_{1}^{n},\cdots,u_{P-1}^{n}]^{\top}$ with $u_{j}^{n}=u_{P}(x_{j},t_{n})$.

In order to discretize Eq. \eqref{eq:FracLayer4Wave} with $c=1$, we assume that
\begin{equation}
    V(x,t_n)
    \approx
    V_{P}(x,t_{n})
    =\sum\limits_{j=0}^{P}V_P(x_{j},t_{n})L_{j}(x),
    %=\sum\limits_{j=0}^{P}V_P(x_{j},t_{n})\sum\limits_{k=0}^{P}\ell_{k}^{j}p^{a,b}_{k}(\hat{x}),
\end{equation}
and
\begin{equation}
    W(x,t_n)
    \approx
    W_{P}(x,t_{n})
    =\sum\limits_{j=0}^{P}W_{P}(x_{j},t_{n})L_{j}(x).
\end{equation}
If the time derivative is discretized by the Crank-Nicolson formula and the space fractional derivatives are evaluated by the spectral collocation method, the corresponding fully discrete scheme for Eq. \eqref{eq:FracLayer4Wave} with $c=1$ is given by
\begin{equation}\label{eq:FullyDiscrete_1D}
\left\{
\begin{aligned}
  &\frac{\mathbf{V}^{n+1}-\mathbf{V}^{n}}{\tau}
=\frac{1}{2}\mathbf{D}^{\bm\alpha}_{R}\left(\mathbf{V}^{n+1}+\mathbf{V}^{n}\right),
\\[3pt]
&\frac{\mathbf{W}^{n+1}-\mathbf{W}^{n}}{\tau}
=\frac{1}{2}\mathbf{D}^{\bm\alpha}_{L}\left(\mathbf{W}^{n+1}+\mathbf{W}^{n}\right).
\end{aligned}
\right.
\end{equation}
Equivalently, the fully discrete scheme \eqref{eq:FullyDiscrete_1D} can be written as
\begin{equation}
\left\{
\begin{aligned}
&\mathbf{V}^{n+1}=\left[\mathbf{I}-\frac{\tau}{2}\mathbf{D}^{\bm\alpha}_{R}\right]^{-1}\left[\left(\mathbf{I}+\frac{\tau}{2}\mathbf{D}^{\bm\alpha}_{R}\right)\mathbf{V}^{n}\right],
\\[3pt]
&\mathbf{W}^{n+1}=\left[\mathbf{I}-\frac{\tau}{2}\mathbf{D}^{\bm\alpha}_{L}\right]^{-1}\left[\left(\mathbf{I}+\frac{\tau}{2}\mathbf{D}^{\bm\alpha}_{L}\right)\mathbf{W}^{n}\right].
\end{aligned}
\right.
\end{equation}
Here $\mathbf{V}^{n}=[V_{1}^{n},\cdots,V_{P-1}^{n}]^{\top}$ and $\mathbf{W}^{n}=[W_{1}^{n},\cdots,W_{P-1}^{n}]^{\top}$ with $V_{j}^{n}=V_{P}(x_{j},t_{n})$ and $W_{j}^{n}=W_{P}(x_{j},t_{n})$.

Let the time derivative be discretized by the two-step Adams–Bashforth and the space derivatives be approximated by the spectral collocation method. The fully discrete scheme for Eq. \eqref{eq:FracLayer4ScalarWaveEq_2D} with $c=1$ reads,
\begin{equation}\label{eq:FullyDiscrete_AB1}
\left\{
\begin{aligned}
\mathbf{v}^{1}
=&\mathbf{v}^{0}
+\frac{\tau}{2}\left[\widetilde{\mathbf{D}}^{\bm\alpha}_{L}+\widetilde{\mathbf{D}}^{\bm\alpha}_{R}\right]\mathbf{v}^{0}
+\frac{\tau}{2}\mathbf{v}^{0}\left[\overline{\mathbf{D}}^{\bm\beta}_{L}+\overline{\mathbf{D}}^{\bm\beta}_{R}\right]^{\top}
\\[3pt]
&+\frac{\tau}{2}\left[\widetilde{\mathbf{D}}^{\bm\alpha}_{L}-\widetilde{\mathbf{D}}^{\bm\alpha}_{R}\right]\mathbf{w}_{1}^{0}
+\frac{\tau}{2}\mathbf{w}_{2}^{0}\left[\overline{\mathbf{D}}^{\bm\beta}_{L}-\overline{\mathbf{D}}^{\bm\beta}_{R}\right]^{\top},
\\[3pt]
\mathbf{w}^{1}_{1}
=&\mathbf{w}^{0}_{1}
+\frac{\tau}{2}\left[\widetilde{\mathbf{D}}^{\bm\alpha}_{L}-\widetilde{\mathbf{D}}^{\bm\alpha}_{R}\right]\mathbf{v}_{1}^{0}
+\frac{\tau}{2}\left[\widetilde{\mathbf{D}}^{\bm\alpha}_{L}+\widetilde{\mathbf{D}}^{\bm\alpha}_{R}\right]\mathbf{w}^{0}_{1},
\\[3pt]
\mathbf{w}^{1}_{2}
=&\mathbf{w}^{0}_{2}
+\frac{\tau}{2}\mathbf{v}_{1}^{0}\left[\overline{\mathbf{D}}^{\bm\beta}_{L}-\overline{\mathbf{D}}^{\bm\beta}_{R}\right]^{\top}
+\frac{\tau}{2}\mathbf{w}^{0}_{2}\left[\overline{\mathbf{D}}^{\bm\beta}_{L}+\overline{\mathbf{D}}^{\bm\beta}_{R}\right]^{\top},
\end{aligned}
\right.
\end{equation}
and for $n\geq 1$,
\begin{equation}\label{eq:FullyDiscrete_AB2s}
\left\{
\begin{aligned}
\mathbf{v}^{n+1}
=&\mathbf{v}^{n}
+\frac{\tau}{2}\left[\widetilde{\mathbf{D}}^{\bm\alpha}_{L}+\widetilde{\mathbf{D}}^{\bm\alpha}_{R}\right]\left(\frac{3}{2}\mathbf{v}^{n}-\frac{1}{2}\mathbf{v}^{n-1}\right)
+\frac{\tau}{2}\left(\frac{3}{2}\mathbf{v}^{n}-\frac{1}{2}\mathbf{v}^{n-1}\right)\left[\overline{\mathbf{D}}^{\bm\beta}_{L}+\overline{\mathbf{D}}^{\bm\beta}_{R}\right]^{\top}
\\[3pt]
&+\frac{\tau}{2}\left[\widetilde{\mathbf{D}}^{\bm\alpha}_{L}-\widetilde{\mathbf{D}}^{\bm\alpha}_{R}\right]\left(\frac{3}{2}\mathbf{w}_{1}^{n}-\frac{1}{2}\mathbf{w}_{1}^{n-1}\right)
+\frac{\tau}{2}\left(\frac{3}{2}\mathbf{w}_{2}^{n}-\frac{1}{2}\mathbf{w}_{2}^{n-1}\right)\left[\overline{\mathbf{D}}^{\bm\beta}_{L}-\overline{\mathbf{D}}^{\bm\beta}_{R}\right]^{\top},
\\[3pt]
\mathbf{w}^{n+1}_{1}
=&\mathbf{w}^{0}_{1}
+\frac{\tau}{2}\left[\widetilde{\mathbf{D}}^{\bm\alpha}_{L}-\widetilde{\mathbf{D}}^{\bm\alpha}_{R}\right]\left(\frac{3}{2}\mathbf{v}^{n}-\frac{1}{2}\mathbf{v}^{n-1}\right)
+\frac{\tau}{2}\left[\widetilde{\mathbf{D}}^{\bm\alpha}_{L}+\widetilde{\mathbf{D}}^{\bm\alpha}_{R}\right]\left(\frac{3}{2}\mathbf{w}_{1}^{n}-\frac{1}{2}\mathbf{w}_{1}^{n-1}\right),
\\[3pt]
\mathbf{w}^{n+1}_{2}
=&\mathbf{w}^{n}_{2}
+\frac{\tau}{2}\left(\frac{3}{2}\mathbf{v}^{n}-\frac{1}{2}\mathbf{v}^{n-1}\right)\left[\overline{\mathbf{D}}^{\bm\beta}_{L}-\overline{\mathbf{D}}^{\bm\beta}_{R}\right]^{\top}
+\frac{\tau}{2}\left(\frac{3}{2}\mathbf{w}_{2}^{n}-\frac{1}{2}\mathbf{w}_{2}^{n-1}\right)\left[\overline{\mathbf{D}}^{\bm\beta}_{L}+\overline{\mathbf{D}}^{\bm\beta}_{R}\right]^{\top}.
\end{aligned}
\right.
\end{equation} 

%%%%%%%%%%%%%%%
\section{Perfectly matched layers (PMLs) in two dimensions}\label{Sec:PML_2D}
%%%%%%%%%%%%%%%

We introduce two PML formulations for the two-dimensional wave equation \eqref{eq:WaveEq_2D}. Both of the formulations are based on the idea of splitting the wave and deal with the split equation \eqref{eq:ScalarWaveEq_2D}. 

Following \cite{Johnson2007}, the PML formualtion takes the following form, 
\begin{equation}\label{eq:PML4WaveEq_2D_I}
    \text{PML I: \qquad}
    \left\{
      \begin{aligned}
        &\frac{\partial}{\partial t}v
        =c\,\frac{\partial}{\partial x}w_{1}
        +c\,\frac{\partial}{\partial y}w_{2}
        -\sigma_{x}v+\psi,
        \\[3pt]
        &\frac{\partial}{\partial t}w_{1}
        =c\,\frac{\partial}{\partial x}v
        -\sigma_{x} w_{1},
        \\[3pt]
        &\frac{\partial}{\partial t}w_{2}
        =c\,\frac{\partial}{\partial y}v, 
        \\[3pt]
        &\frac{\partial}{\partial t}\psi
        =c\,\sigma_{x}\frac{\partial}{\partial y}w_{2},
      \end{aligned}
    \right.
\end{equation}
where $v=\frac{\partial}{\partial t}u$, $w_{1}=\frac{\partial}{\partial x}u$, $w_{2}=\frac{\partial}{\partial y}u$, and $\psi(x,y,0)=0$. Here, the homogeneous Dirichlet boundary conditions are imposed to Eq. \eqref{eq:PML4WaveEq_2D_I}. The damping function $\sigma_{x}(x)$ vanishes at the interior domain and increases from zero in the layers. In the numerical simulation, $\sigma_{x}(x)$ is defined similar to $\alpha(x)$ such that $\sigma_{x}(x)=\alpha(x)-1$ with $\epsilon=0$.

Following \cite{TurkelYefet1998}, the PML for coupled two-dimensional scalar wave equation \eqref{eq:ScalarWaveEq_2D} is given as  
\begin{equation}\label{eq:PML4WaveEq_2D}
    \text{PML II: \qquad}
    \left\{
      \begin{aligned}
        &\frac{\partial}{\partial t}v
        =c\,\frac{\partial}{\partial x}w_{1}
        +c\,\frac{\partial}{\partial y}w_{2}
        -(\sigma_{x}+\sigma_{y}) v
        +c\,\sigma_{x}\frac{\partial}{\partial y}Q
        +c\,\sigma_{y}\frac{\partial}{\partial x}R,
        \\[3pt]
        &\frac{\partial}{\partial t}w_{1}
        =c\,\frac{\partial}{\partial x}v
        -\sigma_{x} w_{1},
        \\[3pt]
        &\frac{\partial}{\partial t}w_{2}
        =c\,\frac{\partial}{\partial y}v
        -\sigma_{y} w_{2},
        \\[3pt]
        &\frac{\partial}{\partial t}Q
        =w_{2},
        \\[3pt]
        &\frac{\partial}{\partial t}R
        =w_{1},
      \end{aligned}
    \right.
\end{equation}
where $v=\frac{\partial}{\partial t}u$, $w_{1}=\frac{\partial}{\partial x}u$, $w_{2}=\frac{\partial}{\partial y}u$, and $Q(x,y,0)=R(x,y,0)=0$. Here, the homogeneous Dirichlet boundary conditions are imposed. The damping functions $\sigma_{x}(x)$ and $\sigma_{y}(y)$ are defined as 
\begin{equation}
     \sigma_{x}(x)
    =\left\{
    \begin{array}{ll}
    \frac{x_L-x}{\delta}\eta, & x\in[x_L-\delta,x_L), 
    \\
    0, &x\in[x_L,x_R], 
    \\
    \frac{x-x_R}{\delta}\eta, & x\in(x_R,x_R+\delta],
    \end{array}
    \right.
\end{equation}
and 
\begin{equation}
     \sigma_{y}(y)
    =\left\{
    \begin{array}{ll}
    \frac{y_L-y}{\delta}\eta, & y\in[y_L-\delta,y_L), 
    \\
    0, &y\in[y_L,y_R], 
    \\
    \frac{y-y_{R}}{\delta}\eta, & y\in(y_R,y_R+\delta], 
    \end{array}
    \right.
\end{equation}
with $\eta\geq0$ being the attenuation constant. Here, we choose $\eta=100$ in the numerical simulation.

\end{document}